\documentclass[11pt]{article}
\usepackage{authblk}
\usepackage[T1]{fontenc}
\usepackage[utf8]{inputenc}
\usepackage[english]{babel}
\usepackage{lmodern}
\usepackage{mathtools} 
\usepackage[a4paper,DIVcalc]{typearea} 
\typearea[10mm]{16} 
\usepackage[pagewise]{lineno}
\providecommand{\abstractt}[1]
{{
  \small
  \textbf{Abstract -} #1
}}

\providecommand{\keywords}[1]
{{
  \small 
  \textbf{Keywords -} #1
}}

\providecommand{\AMSClassification}[1]
{{
  \small
  \textbf{AMS Subject classifications -} #1
}}
\newcommand{\ds}{\displaystyle}
 \usepackage{comment}
\usepackage{cases}
\usepackage{setspace}
\usepackage{amsthm,amssymb,amsmath,eqnarray}
\usepackage{dsfont}
\usepackage{latexsym}
\usepackage[pdftex]{graphicx}
\usepackage{mathdots}
\usepackage[hyperindex=true, pdftex=true, linkcolor=black,colorlinks=true,citecolor=blue,urlcolor=blue]{hyperref}
\usepackage{amssymb,amsfonts,textcomp}
\usepackage{multirow}
\usepackage{longtable}
\usepackage{xcolor}
\definecolor{LightGray}{gray}{0.9}

\usepackage{pgf, tikz}
\usepackage{authblk}

\usepackage{lmodern} 

\title{\LARGE {The Lie symmetry algebra of the Longstaff-Schwartz model}}
\author{OUKNINE Anas*}

\affil{Univ Rouen Normandie, CNRS, Normandie Univ, LMRS UMR 6085, F-76000 Rouen France\\
Email: \texttt{anas.ouknine1@univ-rouen.fr}}

\newtheorem{thm}{Theorem}[section]

\newtheorem{prop}[thm]{Proposition}
\newtheorem{lemma}[thm]{Lemma}
\newtheorem{remark}[thm]{Remark}
\newtheorem{corollary}[thm]{Corollary}
\newtheorem{exm}[thm]{Example}
\numberwithin{equation}{section}
\theoremstyle{remark}

\usepackage{cases}
\usepackage{empheq}
\usepackage{hyperref}

 \allowdisplaybreaks[1]

\begin{document}

\maketitle

\noindent \abstractt{This study uses Lie's theory of symmetries to compute the symmetry group of a class of partial differential equations parameterized by four constants: $u_{t}=-\left((a-bx)u_{x}+(d-ey)u_{y}+\frac{x}{2}u_{xx}+\frac{y}{2}u_{yy}\right)$; under the various conditions on the constants $a,b,d$ and $e$, we deduce the largest and smallest Lie algebra of symmetries, and we also determined the structure of these algebras.}
\newline

\noindent \keywords{Stochastic differential equation, Longstaff-Schwartz model, Kolmogorov backward equation, Itô formula, Lie algebra, symmetry group, prolongation formula, determining system, infinitesimal criterion of symmetry(Lie condition).}\newline

\noindent\AMSClassification{34A26, 91B28 ,60H10}

\section{Introduction}
The Longstaff-Schwartz model is a financial pricing model used to value American-style options. It was developed by Francis Longstaff and Eduardo Schwartz in 1992\cite{longstaff1992two}. The model is widely used in the field of quantitative finance and is particularly effective for valuing options on assets that have complex underlying dynamics.\\
The Longstaff-Schwartz model is base on a two-factor model based on Cox Ingersoll and Ross (CIR) process 
\begin{empheq}[left=\empheqlbrace]{align}
\nonumber d X_{t} &=(a-bX_{t}) d t+\sqrt{X_{t}} d B^{1}_{t}\label{LS}\\
d Y_{t} &=(d-eY_{t}) d t+\sqrt{Y_{t}} d B^{2}_{t}
\end{empheq}
\noindent with $B^{1}$ and $B^{2}$ are two independent Brownian motions.\\
The weighted sum of these factors is considered to form the short-term interest rate:
$$r_{t}=\alpha X_{t}+\beta Y_{t}$$
with $\alpha$ and $\beta$ nonnegative constants.\\
One of the main motivations for developing a multifactor model of the term structure is that single-factor models imply perfect correlation in the instantaneous returns on bonds of all maturities, a characteristic clearly inconsistent with reality.\\
The purpose of this study is to determine the symmetry groups of the class of (1+2) parabolic linear second-order partial differential equations:  
\begin{equation}
 u_{t}=-\mathcal{L}u \label{KBELS}
\end{equation}
where
\begin{equation}
 \nonumber\mathcal{L}u=(a-bx)u_{x}+(d-ey)u_{y}+\frac{x}{2}u_{xx}+\frac{y}{2}u_{yy}
\end{equation}
is the infinitesimal generator of equation \eqref{LS}, where $a,b,d$ and $e$ are constants.\\
\eqref{LS}.\\
It follows from the Ito formula that for each $u\in\mathcal{C}^{1,2}(\mathbb{R}^{+}\times \mathbb{R}^{2})$, $u(t,X_{t},Y_{t})$ is a martingale if and only if $u_{t}+\mathcal{L}u=0$. With a final condition at time $T$,

equation \eqref{KBELS} is known as the Kolmogorov backward equation, describing the backward time evolution of the probability density of the SDE \eqref{LS}, the density $\rho$ verifies $\rho(s,x,y)dy=\mathbb{P}(X_{s}\in dy/_{X_{T}=x})$ , $s\leq T$ and $\rho(T-t,x,y)$ is a solution of the Kolmogorov backward equation.


\noindent As we proceed in the computations,  we must distinguish between cases for the constants $a, b, d$ and $e$ to identify the largest groups of symmetries that leave the equation invariant.   \\
\indent Before starting our work, we observed in various articles that certain software programs, such as Maple, Mathematica, and Python are capable of computing Lie symmetry groups for partial differential equations. However, a challenge arises in the Longstaff Schwartz model, that involves four constant parameters in the drift function. Existing software tools lack the ability to differentiate between these constants or impose necessary conditions (if needed) to identify symmetry groups. This limitation serves as our motivation to initiate computations, seeking to determine these conditions independently so that we can systematically ascertain symmetry groups in all cases and compare our findings with those obtained using Maple.\\

\indent The remainder of this paper is organized as follows. In the next section, we recall the basic results of the Lie point symmetries for partial differential equations with an illustration of the heat equation. Section 3 is dedicated to identifying and analyzing the system of partial differential equations that determines the final structure of the symmetries of the Longstaff-Shwartz equation. Next, in Section 4, we solve the determining system and calculate the associated symmetries, and we give the class of solutions in every case. In the last section, we determine the structure of the Lie algebras linked to each Lie symmetry group and establish an isomorphism between these Lie algebras.
\section{Prerequisites and preliminaries}
\subsection{Symmetry groups}
\indent  In this section, we provide a concise overview of the theory of Lie symmetry groups for differential equations,
closely following the presentation in Craddock and Platen \cite{CraSym}.\\
The theory of Lie symmetry groups offers a systematic method to determine the Lie symmetry group connected to a given differential equation. This involves finding the infinitesimal generators of the group, which are vector fields that generate the symmetry transformations. The generators encode information regarding the symmetries of the equation and can be used to derive exact solutions or reduce the order of the equation.
Lie symmetry groups have been applied in various fields of mathematics and physics. They are particularly important in the study of partial differential equations, where they can be used to classify and solve equations in areas such as fluid dynamics, quantum mechanics, general relativity, and finance. Lie symmetry groups are also connected to other branches of mathematics, including group theory, differential geometry, and representation theory.\\
\indent In Lie symmetry theory, there are at least two methods to find symmetry groups. First, the method of Olver's prolongation, which we use in this study. Olver's book \cite{OLV86}  provides a comprehensive treatment of Lie symmetry groups and their application to differential equations. It explores various mathematical techniques, such as Lie algebras, Lie derivatives, and group actions, to characterize and analyze symmetries. The book also discusses practical algorithms and methods for finding the symmetries of specific differential equations, as well as the connection between symmetries and conservation laws.\\
Second, the method of isovectors introduced by B. Kent Harrison and Frank B. Estabrook in 1971, given a system of partial differential equations after a change in variables and/or possible unknowns, they can be expressed as the vanishing of a family of differential forms. An isovector is then a vector field of all variables preserving the differential ideal generated by the forms; see \cite{harrison1971geometric} for more details on this method.\\
The symmetry of a partial differential equation is a modification or transformation that maps the solutions of an equation to new solutions. More precisely, symmetries allow the construction of a complex solution from simple solutions. If $\mathcal{C}_{\Delta}$ denotes the set of all solutions of the PDE 
\begin{equation}
    \Delta(x,u^{(n)})=0\label{Delta}
\end{equation}
Then we can see the symmetry $\mathcal{S}$ as an automorphism of  $\mathcal{C}_{\Delta}$ \,i.e\, $\mathcal{S}\,:\,\mathcal{C}_{\Delta}\longrightarrow \mathcal{C}_{\Delta}$\\
means if  $u\in\mathcal{C}_{\Delta}$ it implies that $\mathcal{S}u\in\mathcal{C}_{\Delta}.$\\
In the late XIX century, Lie developed a method to determine all symmetry groups for a system of differential equations, providing an essential and powerful tool for the analysis of differential equations. Olver's book \cite{OLV86} provides a comprehensive description of this theory, along with the works of Bluman and Kumei \cite{blu}, \cite{ovsiannikov2014group} by Ovsiannikov, and the studies of Lescot and Zambrini \cite{LESZA,Les12}, who explored the symmetries of the Hamilton-Jacobi-Bellman equation and the Black-Scholes equation used in finance using the isovectors method \cite{harrison1971geometric}.\\
The theory of symmetries is summarized by Lie's theorem, which will be discussed later. For the motivation of our work, we consider a PDE of order $n$ involving  $p$ independent variables and $q$ dependent variables, defined on a connected subset $M\subset X\times U$. The PDE takes the form (\ref{Delta}) where $\Delta(x,y)$ is a $\mathcal{C}^{\infty}$ function on $M\subset X\times U$ where $X$ and $U$ contain independent and dependent variables, respectively.
$$u^{(n)}_{J}=\frac{\partial^{k}u}{\partial x^{j_{1}}_{1}...\partial x^{j_{k}}_{k}}$$
Here, $J = (j_{1}, ... ,j_{k})$ is a multi-index with $j_{i}\in \mathbb{N}$ for $i\in\{1,2,..,k\}$,  
the order of such a multi-index, which we denote by $\#J=k\leq n$, indicates how many derivatives are being taken.
For more details, we refer to Olver's book \cite{OLV86} Chapter 2, which we describe below.\\
We begin by introducing the notion of vector fields in the form of 
\begin{equation}
    \mathbf{v}=\sum_{k=1}^{p}\xi_{k}(x,u)\frac{\partial}{\partial x_{k}}+\sum_{i=1}^{q}\phi_{i}(x,u)\frac{\partial}{\partial u_{i}},\label{Vector}
\end{equation}
defined on $M\subset X\times U.$\\
We can view $\textbf{v}$ as a first-order differential operator: The vector field (\ref{Vector}) is reffered to as the infinitesimal generator of a one parameter local Lie group  called the flow of $\textbf{v}$, that transforms elements $(x,u)\in M\subset X\times U$. This group is denoted as $\mathcal{G}$. Our goal is to determine the conditions for $\xi_{k}$ and $\phi_{i}$ which will guarantee that $\mathcal{G}$ will be a symmetry group of \eqref{Delta}.\\
Before stating Lie's theorem, we must introduce the notion of the \textbf{nth\,\,prolongation} of \textbf{v}, denoted $pr^{(n)}(\textbf{v})$ is defined on the  space $M^{(n)} \subset X \times U^{(n)}$ 
where  $X\times U^{(n)}$is the space of independent
variables, dependent variables and the derivatives of the dependent
variables up to order $n$, we call it the n-th order \textbf{jet space} of the underlying
space $X \times U$.

\noindent First, we define the nth prolongation of group $\mathcal{G}$, which is the natural extension of the action of $\mathcal{G}$ from $(x,u)$ to all derivatives of u up to order $n$, denoted by $pr^{n}\mathcal{G}$ \cite[p.~100]{OLV86}.
The infinitesimal generator of $pr^{n}\mathcal{G}$ is called the nth prolongation of \textbf{v}, and is denoted $pr^{n}\textbf{v}$. It is possible to obtain an explicit formula for $pr^{n}\textbf{v}$. 
\begin{thm}{\textbf{Olver}}\cite[p.~113]{OLV86}\\
 Let $\mathbf{v}$ be a vector field of the form (\ref{Vector}) on an open subset $M\subset X\times U$. Then, the $n$-th prolongation of $\mathbf{v}$ is the vector field: 
$$
pr^{(n)} \mathbf{v}=\mathbf{v}+\sum_{\alpha=1}^q \sum_J \phi_\alpha^J\left(x, u^{(n)}\right) \frac{\partial}{\partial u_J^\alpha}
$$
defined on the corresponding jet space $M^{(n)} \subset X \times U^{(n)}$, the second summation is over all multi-indices $J=\left(j_1, \ldots, j_k\right)$, with $1 \leqslant j_k \leqslant p, 1 \leqslant k \leqslant n$. The coefficient function $\phi_\alpha^J$ of $\mathrm{pr}^{(n)} \mathbf{v}$ is given by the following formula:
$$
\phi_\alpha^J\left(x, u^{(n)}\right)=D_J\left(\phi_\alpha-\sum_{i=1}^p \xi^i u_i^\alpha\right)+\sum_{i=1}^p \xi^i u_{J, i}^\alpha,
$$
where $u_i^\alpha=\partial u^\alpha / \partial x^i$, and $u_{J, i}^\alpha=\partial u_J^\alpha / \partial x^i=\frac{\partial^{k+1}u^{\alpha}}{\partial x^{i}\partial x^{j_{1}}...\partial x^{j_{k}}}$.\\
Here $D_{J}$ denotes the total derivative operator $D$.
\end{thm}
\begin{prop}
    Given $P(x,u^{(n)})$ a smooth function of $x,u$ and derivatives
of $u$ up to order $n$, defined on an open subset $M^{(n)} \subset X \times U^{(n)}$, the total derivative of $P$ with respect to $x^{i}$ has the form $$D_{i}P=\frac{\partial P}{\partial x^{i}}+\sum_{\alpha=1}^{p}\sum_{J} u_{J,i}^{\alpha}\frac{\partial P}{\partial u_{J}^{\alpha}}$$
\end{prop}
\noindent For example, suppose that $q=1$(one dependent variable u) and $p=2$ (2 independent variables x,t). Let $J=(1,1)$ be a multi-index; then, $u_{J}=u_{xt}=\frac{\partial^{2}u}{\partial x\partial t}\,,\, u_{J,t}=u_{xtt}$. It is clear that $\phi^{x}$ can be written for $\phi^{J}$ when $J=(1,0)$, $\phi^{xx}=\phi^{J}$ when $J=(2,0)$ etc.  \\
Now, we will present a version of Lie's theorem, which we call the infinitesimal criterion of symmetry, which is the key to the theory of Lie symmetry groups for differential equations. This allows the necessary and sufficient conditions for the vector field (\ref{Vector}) to generate symmetries of the differential equation. 
\begin{thm}{\textbf{Lie}\cite[p.~165]{OLV86}\label{Liecondition}}
Let $\Delta(x,u^{(n)})=0$ 
be an n-th order partial differential equation as defined above. A local group of transformations $\mathcal{G}$ acting on an open subset $M\subset X\times U$ is a symmetry group of the system if and only if $pr^{n}\mathbf{v}[\Delta(x,u^{(n)})]=0$
whenever $\Delta(x,u^{(n)})=0$, for every infinitesimal generator of $\mathcal{G}$.  
\end{thm}
\noindent When Theorem \eqref{Liecondition} is applied to a partial differential equation, it leads us to what we refer to as the determining system for the functions $\xi_{k}$ and $\phi$. These determining systems are often solved using heuristic methods. Finally, we derive vector fields that generate symmetry groups. What is particularly interesting about these vector fields is that they form Lie algebra under the typical Lie bracket. We obtainede the following results owing to Lie.
\begin{corollary}\label{Liealg}
Let $\Delta(x,u^{(n)})=0$
be an n-th order partial differential equation.The set of all infinitesimal symmetries of this system
forms a Lie subalgebra of the algebra of $\mathcal{C}^{\infty}$-vector fields on M. Moreover, if this subalgebra is finite-dimensional, the symmetry group of the system is a local Lie group of
transformations acting on M.
\end{corollary}
\noindent The procedure to find the group transformation, which is generated by the infinitesimal symmetries that we already find, is known as \textbf{exponentiating} the vector field. To exponentiate the infinitesimal symmetry \textbf{v}$_{k}$, we solve the system of first-order ordinary differential equations: 
\begin{align*}
    \frac{d\tilde{x}}{d\varepsilon}&=\xi(\tilde{x},\tilde{t},\tilde{y})\,\,,\,\, \frac{d\tilde{t}}{d\varepsilon}=\tau(\tilde{x},\tilde{t},\tilde{y})\,\,,\,\,\frac{d\tilde{u}}{d\varepsilon}=\phi(\tilde{x},\tilde{t},\tilde{y}).\\ 
\end{align*}
with initial conditions 
$\tilde{x}(0)=x\,,\,\tilde{t}(0)=t\,,\,\tilde{u}(0)=u.$\\
We define $u^{\varepsilon,v}(t,x):=\Tilde{u}_{\varepsilon}(\tilde{t}_{-\varepsilon},\tilde{x}_{-\varepsilon}).$
\begin{thm}{\label{newSol}}
    For every infinitesimal generator that satisfies Lie's condition \eqref{Liecondition}, we have that $u^{\varepsilon , v}$ is a solution of  PDE \eqref{Delta} for every $\varepsilon$ sufficiently small.
\end{thm}
\noindent For more details on the transformed solution of the PDE, see \cite[p.~94]{OLV86}.\\
\begin{remark}
\begin{itemize}
    \item 
    $u^{\varepsilon ,v}$ is the action of the symmetry generated by v on u. 
    \item 
    $\varepsilon$ is the parameter of the symmetry group.
\end{itemize}   
\end{remark}
\begin{exm}[Heat equation]

We illustrate here the example of the well known one dimensional heat equation, 
\begin{equation}
    u_{xx}=u_{t}\label{heat}.
\end{equation}
Here p=2 and q=1, so  $X=\mathbb{R}^{2}$ and $U=\mathbb{R}$. \\
To find the symmetries of \eqref{heat}, we set $$v=\xi(x,t,u)\frac{\partial}{\partial x}+\tau(x,t,u)\frac{\partial}{\partial t}+\phi(x,t,u)\frac{\partial}{\partial u}$$
and compute the second prolongation of $\mathbf{v}$.Thanks to Lie's theorem $\mathbf{v}$ generates a symmetry group of the heat equation if and only if 
$$pr^{2}\mathbf{v}(u_{xx}-u_{t})=0$$
whenever $u_{xx}-u_{t}=0$,\\
where $$pr^{2}\mathbf{v}=\mathbf{v}+\phi^{x}\frac{\partial}{\partial u_{x}}+\phi^{t}\frac{\partial}{\partial u_{t}}+\phi^{xx}\frac{\partial}{\partial u_{xx}}+\phi^{xt}\frac{\partial}{\partial u_{xt}}+\phi^{tt}\frac{\partial}{\partial u_{tt}}.$$
This condition provides the determining system for $\xi$, $\tau$ and $\phi$. For more details see Olver \cite[p.~120]{OLV86}. \\
Solving the determining system gives the expression of $\xi$, $\tau$ and $\phi$, then we determine the Lie algebra spanned by the infinitesimal symmetries. Thus, the basis of the Lie algebra of symmetry is 
\begin{align*}
    v_{1}&=\frac{\partial}{\partial x}\\
    v_{2}&=\frac{\partial}{\partial t}\\
    v_{3}&=u\frac{\partial}{\partial u},\\
    v_{4}&=x\frac{\partial}{\partial x}+2t\frac{\partial}{\partial t}-\frac{u}{2}\frac{\partial}{\partial u}
    \\v_{5}&=2t\frac{\partial}{\partial x}-xu\frac{\partial}{\partial u},\\
     v_{6}&=4xt\frac{\partial}{\partial x}+4t^{2}\frac{\partial}{\partial t}-(x^{2}+2t)u\frac{\partial}{\partial u}.\\
\end{align*}
Additionally, the infinite dimensional Lie subalgebra spanned by $v_{\beta}=\beta(x,t)\frac{\partial}{\partial u}$ where $\beta$ is an arbitrary solution of the heat equation, which corresponds to the linearity of the equation.\\

\noindent Exponentiating $v_{k}$ and using theorem \ref{newSol} gives the symmetry transformations below:
\begin{align}
    u^{\varepsilon,v_{1}}(t,x)&=u(x-\varepsilon,t)\label{eq1},\\
    u^{\varepsilon,v_{2}}(x,t)&=u(x,t-\varepsilon)\label{eq2},\\
   u^{\varepsilon,v_{3}}(x,t)&=e^{\varepsilon}u(x,t)\label{eq3},\\
    u^{\varepsilon,v_{4}}(x,t)&=e^{-\frac{\varepsilon}{2}}u(e^{-\varepsilon}x,e^{-2\varepsilon}t)\label{eq4},\\
    u^{\varepsilon,v_{5}}(x,t)&=e^{-\varepsilon x+\varepsilon^{2}t}u((x-2\varepsilon t),t)\label{eq5},\\
    u^{\varepsilon,v_{6}}(x,t)&=\frac{1}{\sqrt{1+4\varepsilon t}}\exp\left(\frac{-\varepsilon x^{2}}{1+4\varepsilon t}\right)u\left(\frac{x}{1+4\varepsilon t},\frac{t}{1+4\varepsilon t}\right).\label{eq6}
\end{align}
This implies that if $u(x,t)$ is a solution of the heat equation, then the right-hand side of equations ((\ref{eq1}),...,(\ref{eq6}))  is also a solution for $\varepsilon$ to be sufficiently small. For an application, let us examine the infinitesimal symmetry (\ref{eq6}). When $u(x,t)=1$ is a solution of the heat equation, then by symmetry: $$u^{\varepsilon,v_{6}}(x,t)=\frac{1}{\sqrt{1+4\varepsilon t}}\exp\left(\frac{-\varepsilon x^{2}}{1+4\varepsilon t}\right)\,\,\text{is also a solution.}$$
Allow $t\longrightarrow t-\frac{1}{4\varepsilon}$ and $\varepsilon =\pi$ then $g(x,t)=\frac{1}{\sqrt{4\pi t}}\exp\left(\frac{- x^{2}}{4t}\right)$
which is the fundamental solution for the heat equation that we found from $u=1$ by a group transformation.
\end{exm}
\noindent Lescot and Zambrini \cite{LESZA} found the same symmetries of the heat equation using the isovectors method \cite{harrison1971geometric}. 
\subsection{A Lie alebgra}
We define $\mathcal{M}=<E,F,H,X,Y,Z>$ with $[E,F]=H$, 
$[H,F]=-2F$, $[H,E]=2E$, $[X,Z]=0$, $[Y,Z]=0$ and $[X,Y]=Z$. $\mathcal{M}$ is a Lie algebra. From the commutator relations, we have that, $<E,F,H>\cong \mathfrak{sl}_{2}$ and $<X,Y,Z>\cong \mathfrak{h}_{3}$, the Heisenberg algebra of dimension 3 
 \small
\begin{longtable}{|c|c|c|c|c|c|c|} 
\caption{Lie brackets of the Algebra of heat equation.}\\
\hline
$[.,.]$ & $E$ & $F$&$H$&$X$&$Y$&$Z$\\
\hline
$E$&0&$H$&$-2E$&0&$X$&0\\ 
$F$&$-H$&0&$2F$&$Y$&$0$&0\\
$H$&$2E$&$-2F$&0&$X$&$-Y$&$0$\\
$X$&0&$-Y$&$-X$&0&$Z$&0\\
$Y$&$-X$&0&$Y$&$-Z$&0&0\\
$Z$&0&0&0&0&0&0\\
\hline
\end{longtable}
\noindent From the commutator table, we notice that $<X,Y,Z>$ is an ideal of $\mathcal{M}$, and $<E,F,H>$ is a subalgebra of $\mathcal{M}$, therefore, $\mathcal{M}$ is a semidirect product $\mathfrak{sl}_{2}\ltimes \mathfrak{h}_{3}$. More precisely, from \cite{lescot2014solving} and Theorem 4.1, in \cite{lescot2014symmetries} $ \mathcal{M}\cong \mathcal{H}_{0,0}$, where $\mathcal{H}_{0,0}$ is the Lie symmetry algebra (of dimension 6) of the heat equation $$\frac{\partial u}{\partial t}=-\frac{\partial^{2} u}{\partial q^{2}}.$$
Also, from the commutator table, one can sees that $\mathcal{N}:=\{E,F,H,Z\}$ is a subalgebra of $\mathcal{M}$, and\\ $\mathcal{N}\cong \mathfrak{sl}_{2}\times \mathbb{R}$. Using Theorem (9.1) in \cite{lescot2014symmetries}, for each $C\ne 0$ and $D$, $\mathcal{N}$ is isomorphic to the symmetry algebra of the equation $$\frac{\partial u}{\partial t}=-\frac{\partial^{2} u}{\partial q^{2}}+\left(\frac{C}{q^{2}}+Dq^{2}\right)u.$$

\section{Equations Defining Infinitesimal Symmetries.}
In this section, we aim to determine the Lie algebra of symmetries called $\mathcal{V}_{a,b,d,e}$ for the PDE:
\begin{equation}
    \Delta(x,y,t,u^{(2)})=(a-bx)u_{x}+(d-ey)u_{y}+\frac{x}{2}u_{xx}+\frac{y}{2}u_{yy}+u_{t}=0.\label{I}
\end{equation}
Here we have three
independent variables $x,y$ and  $t$, and one dependent variable $u$, resulting in $p = 3$ and
$q = 1$ in our notation. Equation (\ref{I}) is a second-order PDE with $n = 2$. We can identify the set of its solutions with the linear subvariety in $X \times U^{(2)}$ determined by the
vanishing of $\Delta(x,y,t,u^{(2)})$.\\

\noindent We look for vector fields on $X\times U$  of the form  $$ v=\xi(x,y,t,u)\frac{\partial}{\partial x}+\gamma(x,y,t,u)\frac{\partial}{\partial y}+\tau(x,y,t,u)\frac{\partial}{\partial t}+\phi(x,y,t,u)\frac{\partial}{\partial u}.$$
Our goal is to determine all coefficient
functions $\xi$, $\gamma$ , $\tau$ and $\phi$ such that the corresponding one-parameter group $exp(\varepsilon v)$
is a symmetry group of (\ref{I}). Otherwise, we search for symmetry groups that leave the subvariety $\mathcal{S}_{\Delta}=\{(x,y,t,u^{(2)}):\Delta(x,y,t,u^{(2)})=0\}$ invariant.\\
According to Theorem  (\ref{Liealg}), we
need to calculate the second prolongation 
$$pr^{2}\mathbf{v}=\mathbf{v}+\phi^{x}\frac{\partial}{\partial u_{x}}+\phi^{y}\frac{\partial}{\partial u_{y}}+\phi^{t}\frac{\partial}{\partial u_{t}}+\phi^{xx}\frac{\partial}{\partial u_{xx}}+\phi^{yy}\frac{\partial}{\partial u_{yy}}+\phi^{tt}\frac{\partial}{\partial_{u_{tt}}}+\phi^{xy}\frac{\partial}{\partial u_{xy}}+\phi^{xt}\frac{\partial}{\partial u_{xt}}+\phi^{yt}\frac{\partial}{\partial u_{yt}}$$
of $\mathbf{v}$ whose coefficients are 
\begin{doublespace}
\begin{align*}\ds
\phi^{x}&=\phi_{x}+(\phi_{u}-\xi_{x})u_{x}-\xi_{u}u_{x}^{2}-\gamma_{x}u_{y}-\tau_{x}u_{t}-\gamma_{u}u_{x}u_{y}-\tau_{u}u_{x}u_{t},\\
\phi^{y}&=\phi_{y}+(\phi_{u}-\gamma_{y})u_{y}-\gamma_{u}u_{y}^{2}-\xi_{y}u_{x}-\xi_{u}u_{x}u_{y}-\tau_{y}u_{t}-\tau_{u}u_{t}u_{y},\\
\phi^{t}&=\phi_{t}+(\phi_{u}-\tau_{t})u_{t}-\xi_{t}u_{x}-\xi_{u}u_{t}u_{x}-\gamma_{t}u_{y}-\gamma_{u}u_{t}u_{y}-\tau_{u}u_{t}^{2},\\
\phi^{xx}&=\phi_{xx}+(2\phi_{xu}-\xi_{xx})u_{x}+(\phi_{u}-2\xi_{x})u_{xx}+(\phi_{uu}-2\xi_{xu})u_{x}^{2}-3\xi_{u}u_{x}u_{xx}
-2\gamma_{xu}u_{x}u_{y}-2\gamma_{x}u_{xy}-2\gamma_{u}u_{x}u_{xy}\\
&-\xi_{uu}u_{x}^{3}-\gamma_{xx}u_{y}-\gamma_{u}u_{xx}u_{y}-\gamma_{uu}u_{x}^{2}u_{y}
-2\tau_{x}u_{xt}-2\tau_{u}u_{xt}u_{x}-\tau_{xx}u_{t}-2\tau_{xu}u_{x}u_{t}-\tau_{u}u_{xx}u_{t}-\tau_{uu}u_{x}^{2}u_{t} ,\\
\phi^{yy}&=\phi_{yy}+(2\phi_{yu}-\gamma_{yy})u_{y}+(\phi_{u}-2\gamma_{y})u_{yy}+(\phi_{uu}-2\gamma_{yu})u_{y}^{2}-3\gamma_{u}u_{y}u_{yy}
-2\xi_{y}u_{xy}-2\xi_{u}u_{xy}u_{y}-\xi_{yy}u_{x}\\
&-2\xi_{yu}u_{x}u_{y}
-\xi_{x}u_{yy}-\xi_{uu}u_{x}u_{y}^{2}-\gamma_{uu}u_{y}^{3}-2\tau_{y}u_{yt}-2\tau_{u}u_{y}u_{yt}
-\tau_{yy}u_{t}-2\tau_{uy}u_{t}u_{y}-\tau_{u}u_{t}u_{y}-\tau_{uu}u_{t}u_{y}^{2}.
\end{align*}
\end{doublespace}
\begin{remark}
    In this context, we refrained from calculating $\phi^{tt}, \phi^{yt}, \phi^{xt}$and $\phi^{xy}$ because the terms $u_{tt}, u_{xy}, u_{yt}$ and $u_{xt}$ are not present in Equation (\ref{I}). Consequently, in the prolongation formula, the derivatives of $\Delta$  with respect to $u_{tt}, u_{xy}, u_{yt}$ and $u_{xt}$ vanish.
\end{remark}
\noindent Applying $pr^{(2)}\mathbf{v}$ to $\Delta(x,y,t,u^{(2)})$,
$\textbf{pr}^{(2)}\mathbf{v}(\Delta)=0$ \text{knowing that} $\Delta=0 ,$ with
\begin{align*}\ds
\textbf{pr}^{(2)}\mathbf{v}(\Delta)&=-b\xi u_{x}+\frac{\xi}{2}u_{xx}-e\gamma u_{y}+\frac{\gamma}{2}u_{yy}\\
+(a-bx)&\left[\phi_{x}+(\phi_{u}-\xi_{x})u_{x}-\xi_{u}u_{x}^{2}-\gamma_{x}u_{y}-\tau_{x}u_{t}-\gamma_{u}u_{x}u_{y}-\tau_{u}u_{x}u_{t}\right] \\
+(d-ey)&\left[\phi_{y}+(\phi_{u}-\gamma_{y})u_{y}-\gamma_{u}u_{y}^{2}-\xi_{y}u_{x}-\xi_{u}u_{x}u_{y}-\tau_{y}u_{t}-\tau_{u}u_{t}u_{y}\right]\\
+\frac{x}{2}\Bigg[\phi_{xx}+&(2\phi_{xu}-\xi_{xx})u_{x}+(\phi_{u}-2\xi_{x})u_{xx}+(\phi_{uu}-2\xi_{xu})u_{x}^{2}-3\xi_{u}u_{x}u_{xx}-2\gamma_{xu}u_{x}u_{y}-2\gamma_{x}u_{xy}-2\gamma_{u}u_{x}u_{xy}\\
-\xi_{uu}u_{x}^{3}-&\gamma_{xx}u_{y}-\gamma_{u}u_{xx}u_{y}-\gamma_{uu}u_{x}^{2}u_{y}-2\tau_{x}u_{xt}-2\tau_{u}u_{xt}u_{x}
-\tau_{xx}u_{t}-2\tau_{xu}u_{x}u_{t}-\tau_{u}u_{xx}u_{t}-\tau_{uu}u_{x}^{2}u_{t}\Bigg]\\
+\frac{y}{2}\Bigg[\phi_{yy}+&(2\phi_{yu}-\gamma_{yy})u_{y}+(\phi_{u}-2\gamma_{y})u_{yy}+(\phi_{uu}-2\gamma_{yu})u_{y}^{2}
-3\gamma_{u}u_{y}u_{yy}
-2\xi_{y}u_{xy}-2\xi_{u}u_{xy}u_{y}-\xi_{yy}u_{x}\\
-2\xi_{yu}u_{x}u_{y}&-\xi_{x}u_{yy}-\xi_{uu}u_{x}u_{y}^{2}-\gamma_{uu}u_{y}^{3}-2\tau_{y}u_{yt}
-2\tau_{u}u_{y}u_{yt}-\tau_{yy}u_{t}-2\tau_{uy}u_{t}u_{y}-\tau_{u}u_{t}u_{y}-\tau_{uu}u_{t}u_{y}^{2}\Bigg]=0 
\end{align*}
and replace $ u_{t}$ with $-(a-bx)u_{x}-(d-ey)u_{y}-\frac{x}{2}u_{xx}-\frac{y}{2}u_{yy}$ wherever it occurs.
The infinitesimal criterion of symmetry in (\ref{Liecondition})
results in a polynomial equation for the derivatives
of $u$. The coefficients of this polynomial equation are constructed from $\xi$, $\gamma$, $\tau$ and $\phi$ along with their derivatives. However, because $\xi$, $\gamma$, $\tau$ and $\phi$ depend only on $x, y, t$ and $u$, the coefficients must all be zero. Hence, the 
determining system is:\\

\begin{tabular}{>{\centering\arraybackslash}m{1.5cm}>{\centering\arraybackslash}p{13.5cm}}
        \hline
        Monomials  & {\centering Coefficients  }\\
        \hline
         $u_{x}$  &\begin{spacing}{-0.5}
 \begin{eqnarray} \nonumber
             &&-b\xi-(a-bx)\xi_{x}+(a-bx)^{2}\tau_{x}-(d-ey)\xi_{y}+(a-bx)(d-ey)\tau_{y}-\xi_{t}\\[-2pt]
            &&+(a-bx)\tau_{t}+\frac{x}{2}(2\phi_{xu}-\xi_{xx})+\frac{x}{2}(a-bx)\tau_{xx}-\frac{y}{2}\xi_{yy}+\frac{y}{2}(a-bx)\tau_{yy}=0
            \quad\quad\label{a} 
        \end{eqnarray}\end{spacing}
        \\\hline 
        $u_{y}$ &\begin{spacing}{-0.5}
        \begin{eqnarray} \nonumber
        &&-e\gamma-(d-ey)\gamma_{y}+(d-ey)^{2}\tau_{y}-(a-bx)\gamma_{x}+(a-bx)(d-ey)\tau_{x}-\gamma_{t}\\[-2pt]
        &&+(d-ey)\tau_{t}+\frac{y}{2}(2\phi_{yu}-\gamma_{yy})+\frac{y}{2}(d-ey)\tau_{yy}-\frac{x}{2}\gamma_{xx}+\frac{x}{2}(d-ey)\tau_{xx}=0\quad \quad
        \label{b}
         \end{eqnarray}\end{spacing}
        \\\hline
         $u_{x}u_{y}$ &\begin{spacing}{-0.5}
 \begin{eqnarray} 
             &&-x\gamma_{xu}+x(d-ey)\tau_{xu}-y\xi_{yu}-y(a-bx)\tau_{uy}=0
         \end{eqnarray}\end{spacing}
        \\\hline
         $u_{x}^{2}$ &\begin{spacing}{-0.5}
        \begin{eqnarray} 
         &&-(a-bx)\xi_{u}+(a-bx)^{2}\tau_{u}
        +\frac{x}{2}(\phi_{uu}-2\xi_{xu})+x(a-bx)\tau_{xu}=0\hspace{1cm} 
        \label{c}
        \end{eqnarray}\end{spacing}
             \\\hline 
        $u_{y}^{2}$ &\begin{spacing}{-0.5}
          \begin{eqnarray} 
             &&-(d-ey)\gamma_{u}+(d-ey)^{2}\tau_{u}+\frac{y}{2}(\phi_{uu}-2\gamma_{yu})+y(d-ey)\tau_{uy}=0\hspace{1cm} 
        \label{c'}  
        \end{eqnarray}\end{spacing}
         \\\hline
         $u_{xx}$ & \begin{spacing}{-0.5}
        \begin{eqnarray} 
         &&\frac{\xi}{2}+\frac{x}{2}(a-bx)\tau_{x}+\frac{y}{2}(d-ey)\tau_{y}+\frac{x}{2}(\phi_{u}-2\xi_{x})+\frac{x^{2}}{4}\tau_{xx}+\frac{y^{2}}{4}\tau_{yy}=0 \quad\label{d}         \end{eqnarray}\end{spacing}
          \\\hline
          $u_{yy}$ & \begin{spacing}{-0.5}
        \begin{eqnarray}
          &&\frac{\gamma}{2}+\frac{y}{2}(d-ey)\tau_{y}+\frac{y}{2}(a-bx)\tau_{x}+\frac{y}{2}(\phi_{u}-2\gamma_{y})+\frac{y^{2}}{4}\tau_{yy}+\frac{x^{2}}{4}\tau_{xx}=0\label{e}      \quad \quad   \end{eqnarray}\end{spacing}
           \\\hline
          $u_{x}u_{xx}$ & \begin{spacing}{-0.5}
        \begin{eqnarray}
        &&-x\xi_{u}+\frac{x^{2}}{2}\tau_{xu}=0\quad\quad          \end{eqnarray}\end{spacing}
          \\\hline
          $u_{y}u_{yy}$ & \begin{spacing}{-0.5}
        \begin{eqnarray} 
        &&-y\gamma_{u}+\frac{y^{2}}{2}\tau_{yu}=0\quad\quad\end{eqnarray}\end{spacing}
          \\\hline
          $u_{yy}u_{x}$ &  \begin{spacing}{-0.5}
        \begin{eqnarray}
          &&\frac{xy}{2}\tau_{xu}=0\quad\quad\end{eqnarray}\end{spacing}
         \\\hline
         $u_{xx}u_{y}$ &  \begin{spacing}{-0.5}
        \begin{eqnarray}
        \frac{xy}{2}\tau_{uy}=0
        \end{eqnarray}\end{spacing}
         \\ \hline
         $u_{xx}u_{yy}$ & \begin{spacing}{-0.5}
        \begin{eqnarray}
        -xy\tau_{u}=0
        \end{eqnarray}\end{spacing}
         \\\hline
         $u_{yy}^{2}$ &\begin{spacing}{-0.5}
        \begin{eqnarray}
        \frac{y^{2}}{4}\tau_{u}=0
         \end{eqnarray}\end{spacing}
          \\\hline
         $u_{xx}^{2}$ & \begin{spacing}{-0.5}
        \begin{eqnarray}
          \frac{x^{2}}{4}\tau_{u}=0
          \end{eqnarray}\end{spacing}
         \\\hline      
         $u_{xy}$ & \begin{spacing}{-0.5}
        \begin{eqnarray}
         x\gamma_{x}+y\xi_{y}=0\hspace{2cm}
         \label{f}
         \end{eqnarray}\end{spacing}
         \\\hline
         $u_{x}u_{xy}$ & \begin{spacing}{-0.5}
        \begin{eqnarray}
         -x\gamma_{u}=0\hspace{2cm}
         \label{g}
         \end{eqnarray}\end{spacing}
         \\\hline
          $u_{x}^{2}u_{y}$ & \begin{spacing}{-0.5}
        \begin{eqnarray}
         -\frac{x}{2}\gamma_{uu}+\frac{x}{2}(d-ey)\tau_{uu}=0
         \end{eqnarray}\end{spacing}
         \\\hline        
         $u_{x}^{3}$ &  \begin{spacing}{-0.5}
        \begin{eqnarray}
         -\frac{x}{2}\xi_{uu}-\frac{x}{2}(a-bx)\tau_{uu}=0
         \end{eqnarray}\end{spacing}
         \\\hline     
         $u_{y}^{3}$ & \begin{spacing}{-0.5}
        \begin{eqnarray}
         -\frac{y}{2}\gamma_{uu}-\frac{y}{2}(d-ey)\tau_{uu}=0
         \end{eqnarray}\end{spacing}
         \\\hline
         $u_{xt}$ &\begin{spacing}{-0.5}
        \begin{eqnarray}
        -x\tau_{x}=0\hspace{2cm}
        \label{h}
        \end{eqnarray}\end{spacing}
         \\\hline
         $u_{xt}u_{x}$ &\begin{spacing}{-0.5}
        \begin{eqnarray}
         -x\tau_{u}=0\hspace{3cm}
         \label{i}
         \end{eqnarray}\end{spacing}
          \\\hline
         \end{tabular}

         \begin{tabular}{>{\centering\arraybackslash}m{1.5cm}>{\centering\arraybackslash}p{13.5cm}}
        \hline
          
        Monomes  & {\centering Coefficients  }\\
        \hline
          $u_{x}u_{y}^{2}$ & \begin{spacing}{-0.5}
        \begin{eqnarray}
        -\frac{y}{2}\xi_{uu}+\frac{y}{2}(a-bx)\tau_{uu}=0
        \end{eqnarray}\end{spacing}
         \\\hline
           $u_{x}^{2}u_{t}$ & \begin{spacing}{-0.5}
        \begin{eqnarray}
        \frac{x}{2}\tau_{uu}=0
        \end{eqnarray}\end{spacing}
         \\\hline
         $u_{y}^{2}u_{t}$ & \begin{spacing}{-0.5}
        \begin{eqnarray}
        \frac{y}{2}\tau_{uu}=0
         \end{eqnarray}\end{spacing}
         \\\hline
          $u_{xy}u_{y}$ &  \begin{spacing}{-0.5}
        \begin{eqnarray}
         -y\xi_{u}=0\hspace{1cm}
         \label{j}
         \end{eqnarray}\end{spacing}
         \\\hline
          $u_{yt}$ & \begin{spacing}{-0.5}
        \begin{eqnarray}
         y\tau_{y}=0\label{(3.27)}
         \label{k}
         \end{eqnarray}\end{spacing}
         \\\hline
          $u_{y}u_{yt}$ &\begin{spacing}{-0.5}
        \begin{eqnarray}
         \hspace{1cm}-y\tau_{u}=0
         \end{eqnarray}\end{spacing}
         \\\hline
         ....... & \begin{spacing}{-0.5}
        \begin{eqnarray}
         \frac{y}{2}(d-ey)\tau_{uu}=0
         \end{eqnarray}\end{spacing}
         \\\hline
         $ u_{y}^{2}u_{yy}$& \begin{spacing}{-0.5}
        \begin{eqnarray}
         \frac{y^{2}}{4}\tau_{uu}=0
          \end{eqnarray}\end{spacing}
         \\\hline
        $u_{x}^{2}u_{xx}$& \begin{spacing}{-0.5}
        \begin{eqnarray}
        \frac{x^{2}}{4}\tau_{uu}=0
        \end{eqnarray}\end{spacing}
         \\\hline
         constant term & \begin{spacing}{-0.5}
        \begin{eqnarray}
        (a-bx)\phi_{x}+(d-ey)\phi_{y}+\frac{x}{2}\phi_{xx}+\frac{y}{2}\phi_{yy}+\phi_{t}=0\label{(3.32)}\hspace{1cm}
        \end{eqnarray}\end{spacing}
          \\\hline
\end{tabular}
\bigskip
\noindent First, Equations (\ref{h}),(\ref{i}) and (\ref{(3.27)}) require that $\tau$ be solely a function of  $t$. Subsequently, Equations(\ref{g}) and (\ref{j}) demonstrate that  $\xi$ and  $\gamma$ are independent of $u$.\\
Furthermore, according to (\ref{c}) and (\ref{c'}), $ \phi$ is affine in $u$, implying that $$\phi(x,y,t,u)=\lambda(x,y,t)u+\mu(x,y,t)\label{phi is affine on u}$$
for certain functions $\lambda$ and $\mu$.
Substituting $\phi$ into equation  (\ref{(3.32)}) reveals that $\lambda$ and  $\mu$ satisfy equation (\ref{I}).\\
Since $\phi$ is affine on $u$, we define $\mathcal{L}_{a,b,d,e}=\{v\in \mathcal{V}_{a,b,d,e}\,\,|\,\,\mu_{v}=0\}$, a subspace of $\mathcal{V}_{a,b,d,e}$ where $\mu_{v}:=\mu$.
\begin{lemma}
    $\mathcal{L}_{a,b,d,e}$ is a subalgebra of $\mathcal{V}_{a,b,d,e}$.
\end{lemma}
\begin{proof}
    Consider $v\,,v^{\prime} \in \mathcal{L}_{a,b,d,e}$ implying $\mu_{v}=\mu_{v^{\prime}}=0$. We aim to show that $[v,v^{\prime}]\in \mathcal{L}_{a,b,d,e}$.
    By the definition of $\mathcal{L}_{a,b,d,e}$, it suffices to prove that $\mu_{[v,v^{\prime}]}=0$. Expressing $\phi$ as \\
    $\phi_{v}=v(u)=\lambda_{v}u+\mu_{v}\,\,\text{we have}\,\,\phi_{[v,v^{\prime}]}=\lambda_{[v,v^{\prime}]}u+\mu_{[v,v^{\prime}]}$.\\
    We compute as follows:
    \begin{align*}
    \phi_{[v,v^{\prime}]}=[v,v^{\prime}](u) &=v(v^{\prime}(u))-v^{\prime}(v(u))\\
        &=v(\lambda_{v^{\prime}}u+\mu_{v^{\prime}})-v^{\prime}(\lambda_{v}u+\mu_{v})\\
        &=(v(\lambda_{v^{\prime}})-v^{\prime}(\lambda_{v}))u+(v(\mu_{v^{\prime}})-v^{\prime}(\mu_{v})+\lambda_{v^{\prime}}\mu_{v}-\lambda_{v}\mu_{v^{\prime}})  
    \end{align*}
    Since $\xi\,,\,\gamma\,,\,\tau$ do not depend on $x\,,y\,,t$, the terms $v(\lambda_{v^{\prime}})-v^{\prime}(\lambda_{v})$ and $v(\mu_{v^{\prime}})-v^{\prime}(\mu_{v})+\lambda_{v^{\prime}}\mu_{v}-\lambda_{v}\mu_{v^{\prime}}$ depend only on $x,y,t$. Thus, $\mu_{[v,v^{\prime}]}=v(\mu_{v^{\prime}})-v^{\prime}(\mu_{v})+\lambda_{v^{\prime}}\mu_{v}-\lambda_{v}\mu_{v^{\prime}}$ and $\lambda_{[v,v^{\prime}]}=v(\lambda_{v^{\prime}})-v^{\prime}(\lambda_{v})$.
     Hence, $\mu_{[v,v^{\prime}]}=0$, providing that $[v,v^{\prime}]\in \mathcal{L}_{a,b,d,e}.$
\end{proof}
\begin{remark}
      $\mathcal{L}_{a,b,d,e}$ represent the space of non trivial symmetries, that is, symmetries that do not arise from the linearity of the equation.
\end{remark}
\noindent Now we seek symmetries of the form:$$v=\xi(x,y,t)\frac{\partial}{\partial x}+\gamma(x,y,t)\frac{\partial}{\partial y}+\tau(t)\frac{\partial}{\partial t}+(\lambda(x,y,t)u+\mu(x,y,t))\frac{\partial}{\partial u}.$$
Because $\tau$ does not depend on $x,y$ or $u$, conditions (\ref{d}) and (\ref{e}) become:
$$
\left\{\begin{aligned}
\frac{\xi}{2}+\frac{x}{2}\tau_{t}-x\xi_{x}&=0,\\
\frac{\gamma}{2}+\frac{y}{2}\tau_{t}-y\gamma_{y}&=0;\\
\end{aligned}\right.\\
$$

\noindent Solving those two equations gives: 
$$
\left\{\begin{aligned}
\xi&=\tau_{t}x+\sqrt{x}f(t,y),\\
\gamma&=\tau_{t}y+\sqrt{y}g(x,t),
\end{aligned}\right.
$$
for certain functions $f$ and $g$.\\
Substituting $\xi$ and $\gamma$ into equation  (\ref{f})  gives  $$2(x\gamma_{x}+y\xi_{y})=2x(\sqrt{y}g_{x})+2y(\sqrt{x}f_{y})=2\sqrt{x}g_{x}+2\sqrt{y}f_{y}=0,$$
thus there exists a function  $h$ such that $-2\sqrt{x}g_{x}=2\sqrt{y}f_{y}=h(t).$\\

$$
\text{Then}\,\,\left\{\begin{aligned}
g_{x}&=-\frac{1}{2\sqrt{x}}h(t),\\
f_{y}&=\frac{1}{2\sqrt{y}}h(t),
\end{aligned}\right.
$$
$$
\text{so}\,\,\left\{\begin{aligned}
g(x,t)&=-\sqrt{x}h(t)+k(t),\\
f(y,t)&=\sqrt{y}h(t)+l(t),
\end{aligned}\right.
$$
for certain functions  $k$ and  $l$. We remind that the unknown functions are  $h(t), l(t), k(t), \tau(t)$.\\
$$
\text{We have }\,\,\left\{\begin{aligned}
\xi&=\tau_{t}x+\sqrt{x}(\sqrt{y}h(t)+l(t)),\\
\gamma&=\tau_{t}y+\sqrt{y}(-\sqrt{x}h(t)+k(t)).
\end{aligned}\right.
$$
Replacing  $\xi$ and $\gamma$ in the  equations  (\ref{a}) and  (\ref{b}) gives 
{\small
$$
\left\{\begin{aligned}
&-b\left(\tau_{t}x+\sqrt{x}(\sqrt{y}h+l)\right)-(a-bx)\left(\frac{h}{2}\frac{\sqrt{y}}{\sqrt{x}}+\frac{l}{2}\frac{1}{\sqrt{x}}\right)-(d-ey)\left(\frac{h}{2}\frac{\sqrt{x}}{\sqrt{y}}\right)\\
&-\tau_{tt}x-h_{t}\sqrt{x}\sqrt{y}
-\sqrt{x}l_{t}+x\lambda_{x}+\frac{x}{2}\left(\frac{h}{4}\sqrt{y}x^{-\frac{3}{2}}+\frac{l}{4}x^{-\frac{3}{2}}\right)+\frac{y}{2}\left(\frac{h}{4}\sqrt{x}y^{-\frac{3}{2}}\right)=0,\\\\
&-e\left(\tau_{t}y+\sqrt{y}(-\sqrt{x}h+k)\right)-(d-ey)\left(-\frac{h}{2}\frac{\sqrt{x}}{\sqrt{y}}+\frac{k}{2}\frac{1}{\sqrt{y}}\right)-(a-bx)\left(-\frac{h}{2}\frac{\sqrt{y}}{\sqrt{x}}\right)\\&-\tau_{tt}y+\sqrt{x}\sqrt{y}h_{t}-k_{t}\sqrt{y}+y\lambda_{y}-\frac{y}{2}\left(\frac{h}{4}\sqrt{x}y^{-\frac{3}{2}}-\frac{k}{4}y^{-\frac{3}{2}}\right)-\frac{x}{2}\left(\frac{h}{4}\sqrt{y}x^{-\frac{3}{2}}\right)=0,
\end{aligned}\right.
$$
}
{\small
simplifying we have 
\begin{empheq}[left=\empheqlbrace]{align}
\nonumber\lambda(x,y,t)&=(\tau_{tt}+b\tau_{t})x+2\left(h_{t}+\frac{h}{2}(b-e)\right)\sqrt{x}\sqrt{y}+h\left(\frac{1}{4}-a\right)\frac{\sqrt{y}}{\sqrt{x}}+h\left(d-\frac{1}{4}\right)\frac{\sqrt{x}}{\sqrt{y}}\label{*}
\\
&+2\left(l_{t}+\frac{bl}{2}\right)\sqrt{x}+l\left(\frac{1}{4}-a\right)\frac{1}{\sqrt{x}}+\theta(y,t),
\\\nonumber\\
\nonumber\lambda(x,y,t)&=(\tau_{tt}+e\tau_{t})y+2\left(-h_{t}+\frac{h}{2}(b-e)\right)\sqrt{x}\sqrt{y}+h\left(\frac{1}{4}-a\right)\frac{\sqrt{y}}{\sqrt{x}}+h\left(d-\frac{1}{4}\right)\frac{\sqrt{x}}{\sqrt{y}}\label{**}
\\
&+2\left(k_{t}+\frac{ek}{2}\right)\sqrt{y}
+k\left(\frac{1}{4}-d\right)\frac{1}{\sqrt{y}}+\beta(x,t),
\end{empheq}
}
for some functions $\theta$ and $\beta$.\\
For a fixed $t$ and $y$, utilizing \eqref{*}, we observe that $\sqrt{x}\lambda(x,y,t)\underset{x\longrightarrow 0}\longrightarrow 0$ and\\ 
    $\sqrt{x}\lambda(x,y,t)\underset{x\longrightarrow 0}\longrightarrow h(\frac{1}{4}-a)\sqrt{y}+l(\frac{1}{4}-a)$. From the uniqueness of the limit, we obtain \\$h(\frac{1}{4}-a)\sqrt{y}+l(\frac{1}{4}-a)=0$. Similarly, for a fixed $t$ and $x$, we find that $h(\frac{1}{4}-d)\sqrt{x}+k(\frac{1}{4}-d)=0$.\\
    Thus, for every $x$ and $y$: 
    \begin{align*}
        h(\frac{1}{4}-a)\sqrt{y}+l(\frac{1}{4}-a)&=0\\
        h(\frac{1}{4}-d)\sqrt{x}+k(\frac{1}{4}-d)&=0.
    \end{align*} 
    Since $h,k$ and $l$ depend only on $t$, we have \begin{align}
        h(\frac{1}{4}-a)&=h(\frac{1}{4}-d)=l(\frac{1}{4}-a)=k(\frac{1}{4}-d)=0\,\,\text{for every}\,\,t.\label{h,k,l}
    \end{align}

\begin{remark}

    \begin{itemize}
    \item If $a\ne \frac{1}{4}$ :\\
    From \ref{h,k,l} we deduce that $h=l=0$.
    \item If $d\ne \frac{1}{4}$ :\\
    From \ref{h,k,l} we deduce that $h=k=0$.
\end{itemize}
\end{remark}
\noindent Thus equations \eqref{*} and \eqref{**} becomes 
{\small

\begin{empheq}[left=\empheqlbrace]{align}
\lambda(x,y,t)&=(\tau_{tt}+b\tau_{t})x+2\left(h_{t}+\frac{h}{2}(b-e)\right)\sqrt{x}\sqrt{y}+2\left(l_{t}+\frac{bl}{2}\right)\sqrt{x}+\theta(y,t),\label{*'}\\
\lambda(x,y,t)&=(\tau_{tt}+e\tau_{t})y+2\left(-h_{t}+\frac{h}{2}(b-e)\right)\sqrt{x}\sqrt{y}+2\left(k_{t}+\frac{ek}{2}\right)\sqrt{y}+\beta(x,t),\label{**'}
\end{empheq}
}

 

\noindent From \eqref{*'} and \eqref{**'}  we have 

\begin{align*}
\lim\limits_{x \to 0}&\Bigg[(\tau_{tt}+b\tau_{t})x+2\left(h_{t}+\frac{h}{2}(b-e)\right)\sqrt{x}\sqrt{y}
+2\left(l_{t}+\frac{bl}{2}\right)\sqrt{x}+\theta(y,t)\Bigg]\\
=\lim\limits_{x \to 0}&\Bigg[(\tau_{tt}+e\tau_{t})y+2\left(-h_{t}+\frac{h}{2}(b-e)\right)\sqrt{x}\sqrt{y}+2\left(k_{t}+\frac{ek}{2}\right)\sqrt{y}
+\beta(x,t)\Bigg]
\end{align*}
and 

\begin{align*}
\lim\limits_{y \to 0}&\Bigg[(\tau_{tt}+b\tau_{t})x+2\left(h_{t}+\frac{h}{2}(b-e)\right)\sqrt{x}\sqrt{y}+2\left(l_{t}+\frac{bl}{2}\right)\sqrt{x}+\theta(y,t)\Bigg]\\
=\lim\limits_{y \to 0}&\Bigg[(\tau_{tt}+e\tau_{t})y+2\left(-h_{t}+\frac{h}{2}(b-e)\right)\sqrt{x}\sqrt{y}
+2\left(k_{t}+\frac{ek}{2}\right)\sqrt{y}+\beta(x,t)\Bigg].
\end{align*}
Therefore \begin{empheq}[left=\empheqlbrace]{align}
\nonumber\theta(y,t)&=(\tau_{tt}+e\tau_{t})y+2\left(k_{t}+\frac{ek}{2}\right)\sqrt{y}+s(t),\\
\beta(x,t)&=(\tau_{tt}+b\tau_{t})x+2\left(l_{t}+\frac{bl}{2}\right)\sqrt{x}+s(t)\label{theta et beta}.
\end{empheq}

\noindent where $s(t):=\beta(0,t)=\theta(0,t).$

{\small
\begin{remark}\label{h_{t}=0}
    Using \eqref{*'}, \eqref{**'}, and \eqref{theta et beta}, we obtain: 
    \begin{align*}
        \lambda(x,y,t)&=(\tau_{tt}+b\tau_{t})x+2\left(h_{t}+\frac{h}{2}(b-e)\right)\sqrt{x}\sqrt{y}+2\left(l_{t}+\frac{bl}{2}\right)\sqrt{x}+(\tau_{tt}+e\tau_{t})y+2\left(k_{t}+\frac{ek}{2}\right)\sqrt{y}+s(t)
        \\
        &=(\tau_{tt}+e\tau_{t})y+2\left(-h_{t}+\frac{h}{2}(b-e)\right)\sqrt{x}\sqrt{y}+2\left(k_{t}+\frac{ek}{2}\right)\sqrt{y}+(\tau_{tt}+b\tau_{t})x+2\left(l_{t}+\frac{bl}{2}\right)\sqrt{x}+s(t).
    \end{align*}
    So, $4h_{t}\sqrt{x}\sqrt{y}=0$ for every $x$ and $y$. Therefore, $h_{t}=0.$
    \end{remark}}
    
\noindent Because $\lambda$ satisfies (\ref{I}), and from equation \eqref{h,k,l} and remark \ref{h_{t}=0}, then
{\small
\begin{doublespace}
\begin{align*}\ds
    &(a-bx)\left[\tau_{tt}+b\tau_{t}+\frac{h}{2}\left(b-e\right)\frac{\sqrt{y}}{\sqrt{x}}+\left(l_{t}+\frac{bl}{2}\right)\frac{1}{\sqrt{x}}\right]+(d-ey)\Bigg[\tau_{tt}+e\tau_{t}+\frac{h}{2}\left(b-e\right)\frac{\sqrt{x}}{\sqrt{y}}+\left(k_{t}+\frac{ek}{2}\right)\frac{1}{\sqrt{y}}\Bigg]\\
    &+\frac{x}{2}\Bigg[-\frac{h}{4}\left(b-e\right)\sqrt{y}x^{-\frac{3}{2}}-\frac{1}{2}\left(l_{t}+\frac{bl}{2}\right)x^{-\frac{3}{2}}\Bigg]+\frac{y}{2}\left[-\frac{h}{4}\left(b-e\right)\sqrt{x}y^{-\frac{3}{2}}-\frac{1}{2}\left(k_{t}+\frac{ek}{2}\right)y^{-\frac{3}{2}}\right]\\
    &+\Bigg[(\tau_{ttt}+b\tau_{tt})x+2\left(h_{tt}+\frac{h_{t}}{2}(b-e)\right)\sqrt{x}\sqrt{y}+2\left(l_{tt}+\frac{bl_{t}}{2}\right)\sqrt{x}+(\tau_{ttt}+e\tau_{tt})y
    +2\left(k_{tt}+\frac{ek_{t}}{2}\right)\sqrt{y}+s_{t}\Bigg]=0.
\end{align*}
\end{doublespace}}
{\small
\begin{doublespace}
\noindent
Therefore
\begin{align*}
&(\tau_{ttt}-b^{2}\tau_{t})x+(\tau_{ttt}-e^{2}\tau_{t})y-\frac{h}{2}(b^{2}-e^{2})\sqrt{x}\sqrt{y}
+\left(2l_{tt}-b^{2}\frac{l}{2}\right)\sqrt{x}+\left(2k_{tt}-e^{2}\frac{k}{2}\right)\sqrt{y}\\
&+((a+d)\tau_{tt}+(ab+de)\tau_{t}+s_{t})=0\,\,\text{for every}\,\,x\,\, \text{and} \,\,y.
\end{align*}
\end{doublespace}
}
\noindent Then
{\small
\begin{subequations}\label{eq:system}
\begin{empheq}[left=\empheqlbrace]{align}
&h_{t}=0,\\
&(b^{2}-e^{2})\tau_{t}=0\label{b^{2}},\\
&(b^{2}-e^{2})h=0\label{eb},\\
&l_{tt}-\frac{b^{2}}{4}l=0\label{l},\\
&k_{tt}-\frac{e^{2}}{4}k=0\label{k},\\
&(a+d)\tau_{tt}+(ab+de)\tau_{t}+s_{t}=0.\label{s}
\end{empheq}
\end{subequations}
}
and $\mu$ a solution of equation \eqref{I}.\\
Recall that the unknowns are $\tau,\,h,\,l,\,k,\,s,\,$ and
$\mu$. Thus, the coefficients of the vector field are given by: 
{\small
\begin{align*}
    \xi(t,x,y)&=\tau_{t}x+\sqrt{x}(\sqrt{y}h(t)+l(t)),\\
    \gamma(t,x,y)&=\tau_{t}y+\sqrt{y}(-\sqrt{x}h(t)+k(t)),\\
    \phi(t,x,y,u)&=\left[(\tau_{tt}+b\tau_{t})x+h(b-e)\sqrt{x}\sqrt{y}+2\left(l_{t}+\frac{bl}{2}\right)\sqrt{x}+(\tau_{tt}+e\tau_{t})y+2\left(k_{t}+\frac{ek}{2}\right)\sqrt{y}+s(t)\right]u+\mu(t,x,y).
\end{align*}
}
This provides the final structure of the symmetry group. To proceed further, it is essential to distinguish between cases based on conditions involving the constants $a,b,d$ and $e$. \\


\section{Computing the Infinitesimal Symmetries.}
Now, we determine the Lie algebras of symmetries and compute their Lie brackets. From equations \ref{h,k,l} and \ref{eb}, we split our analysis into $(4\times 4)$ cases:
$$(b=e,\,\,b\ne 0),\,\,(b=-e,\,\,b\ne 0),\,\,(b=e=0),\,\,(b\ne\pm e)$$
and $$(a=\frac{1}{4},d=\frac{1}{4}),\,(a=\frac{1}{4},d\ne \frac{1}{4}),\,(a\ne\frac{1}{4},d= \frac{1}{4}),\,(a\ne\frac{1}{4},d\ne \frac{1}{4}).$$
\begin{remark}
We observe that the determining equations \eqref{b^{2}}, \eqref{eb}, \eqref{l}, \eqref{k}, and \eqref{s} 
always have three particular solutions.
\begin{align*}
    \xi&=\gamma=\tau=\mu=0,\,\, \lambda=cte \,\,\,\text{means}\,\,\, h=l=k=\tau=0,\,\,s=cte,\\
    \xi&=\gamma=\lambda=\mu=0,\,\, \tau=cte \,\,\,\text{means}\,\,\, h=l=k=s=0,\,\,\tau=cte,\\
     \xi&=\gamma=\tau=0 ,\,\, \lambda=0 \,\,\text{and}\,\, \mu\,\, \text{is a solution of the PDE}\,\,\eqref{I}\,\text{means}\,\,\, h=l=k=\tau=0,\,\,s=0\\
     &\text{and} \,\, \mu\,\, \text{is a solution} ,\\
\end{align*}
These solutions provides us with multiples of the following symmetries:
\begin{align*}
v_{u}&=u\frac{\partial}{\partial u},\quad v_{t}=\frac{\partial}{\partial t},\,\,v_{\mu}=\mu(x,y,t)\frac{\partial}{\partial u} \,\,\text{where}\,\, \mu \,\,\text{is a solution of} \,\,\eqref{I}.
\end{align*}
$v_{\mu}$ corresponds to the linear superposition principle, which generates infinite-dimensional subalgebra. $v_{u}$ illustrates that we can multiply the solutions by constants, and $v_{t}$ reflects the invariance of the equation under translation in time.

 \end{remark}
    \subsection{Case 1 ($a=\frac{1}{4},d= \frac{1}{4}$)}
\subsubsection{Subcase 1.1 ($b\ne\pm e\,,\,a=\frac{1}{4},d= \frac{1}{4}$)}\label{case1.1}
From equations \eqref{b^{2}}, \eqref{eb}, \eqref{l}, \eqref{k}, and \eqref{s}, we obtain: 
$$
\tau_{t}=0,\,\,h=0,\,\,l_{tt}-b^{2}\frac{l}{4}=0,
\,\,k_{tt}-e^{2}\frac{k}{4}=0\,\,\text{and}\,\,
s_{t}=0,
$$

\noindent solving these equations yields:
$$
\left\{\begin{aligned}
\xi&=\sqrt{x}\left(c_{3}\exp(\frac{b}{2}t)+c_{4}\exp(-\frac{b}{2}t)\right),\\
\gamma&=\sqrt{y}\left(c_{5}\exp(\frac{e}{2}t)+c_{6}\exp(-\frac{e}{2}t)\right),\\
\tau&=c_{1},\\
\phi&=\left(2b\exp(\frac{b}{2}t)c_{3}\sqrt{x}+2e\exp(\frac{e}{2}t)c_{5}\sqrt{y}+c_{2} \right)u+\mu(x,y,t).
\end{aligned}\right.
$$

\noindent where $c_{1},c_{2},c_{3},c_{4},c_{5},c_{6}$ are arbitrary constants and  $\mu(x,y,t)$ is an arbitrary solution of the PDE (\ref{I}). For $c_{i}= 1$ and $c_{j}= 0 \,\,\text{if}\,\, i\ne j$, so the basis of the Lie algebra is:
{\small
\begin{doublespace}
\begin{align*}\ds
    v_{1}&=\frac{\partial}{\partial t}\,\,,\,\, v_{2}=u\frac{\partial}{\partial u}\,\,,\,\,v_{3}=\sqrt{x}\exp(\frac{b}{2}t)\frac{\partial}{\partial x}+2b\exp(\frac{b}{2}t)\sqrt{x}u\frac{\partial}{\partial u},\\
    v_{4}&=\sqrt{x}\exp(-\frac{b}{2}t)\frac{\partial}{\partial x}
     v_{5}=\sqrt{y}\exp(\frac{e}{2}t)\frac{\partial}{\partial y}+2e\exp(\frac{e}{2}t)\sqrt{y}u\frac{\partial}{\partial u},\\
     v_{6}&=\sqrt{y}\exp(-\frac{e}{2}t)\frac{\partial}{\partial y}.
\end{align*}
\end{doublespace}
}
\noindent The commutation relations between these vector fields are given by the commutator table of $\mathcal{V}_{a,b,d,e}$ as follows:\\
{\small
\begin{longtable}{|c|c|c|c|c|c|c|c|}\caption{Commutator table($b\ne\pm e\,,\,a=\frac{1}{4},d= \frac{1}{4}$).}\\
\hline
 $[v_{i},v_{j}]$& $v_{1}$ & $v_{2}$&$v_{3}$&$v_{4}$&$v_{5}$&$v_{6}$& $v_{\mu}$ \endhead
\hline
$v_{1}$&0&0&$-\frac{b}{2}v_{3}$&$-\frac{e}{2}v_{4}$&$\frac{b}{2}v_{5}$&$\frac{e}{2}v_{6}$&$v_{\mu_{t}}$\\
$v_{2}$&0&0&0&0&0&0&$-v_{\mu}$\\
$v_{3}$&$\frac{b}{2}v_{3}$&0&0&0&$bv_{2}$&0&$v_{\mu_{1}}$\\
$v_{4}$&$\frac{e}{2}v_{4}$&0&0&0&0&$ev_{2}$&$v_{\mu_{2}}$\\
$v_{5}$&$-\frac{b}{2}v_{5}$&0&$-bv_{2}$&0&0&0&$v_{\mu_{3}}$\\
$v_{6}$&$-\frac{e}{2}v_{6}$&0&0&$-ev_{2}$&0&0&$v_{\mu_{4}}$\\
$v_{\mu_{t}}$&$-v_{\mu}$&$v_{\mu}$&$-v_{\mu_{1}}$&$-v_{\mu_{2}}$&$-v_{\mu_{3}}$&$-v_{\mu_{4}}$&0\\ 
\hline
\end{longtable}}

\begin{align*}
    \text{where}\,\,\,\mu_{1}&=\exp(-\frac{b}{2}t)\sqrt{x}\mu_{x}\,\,,\,\, \mu_{2}=\exp(-\frac{e}{2}t)\sqrt{y}\mu_{y},\\
    \mu_{3}&=-2\exp(\frac{b}{2}t)\sqrt{x}b\mu+\exp(\frac{b}{2}t)\sqrt{x}\mu_{x}\,\,,\,\, \mu_{4}=-2\exp(\frac{e}{2}t)\sqrt{y}e\mu+\exp(\frac{e}{2}t)\sqrt{y}\mu_{y}.\\
\end{align*}

\begin{remark}
    In all sebsequent cases, the whole table will correspond to $\mathcal{V}_{a,b,d,e}$ and the table without the last line and column to $\mathcal{L}_{a,b,d,e}$.\label{nu }
\end{remark}
\noindent Corollary (\ref{Liealg}) guarantees that the set of all infinitesimal symmetries is a Lie algebra. We can conclude that if $\mu(x,y,t)$ is any solution of equation (\ref{I}), so are $\mu_{1},\mu_{2},\mu_{3},\mu_{4}$; this is easily checked directly. The one-parameter groups $G_{i}$ generated by the $v_{i}$ is given as follows: This gives the transformed point 
$\exp \left(\varepsilon v_{i}\right)(x,y, t, u)=(\tilde{x},\tilde{y}, \tilde{t}, \tilde{u})$

\begin{doublespace}
\begin{align*}\ds
G_{1}: & (x,y, t+\varepsilon, u), \\
G_{2}: & (x,y, t, \exp(\varepsilon)u), \\
G_{3}: & \left(\frac{\left(\varepsilon \exp(\frac{bt}{2})+2\sqrt{x}\right)^{2}}{4},y, t,u\exp\left(^{\frac{ \exp(\frac{bt}{2})b\varepsilon\left(\varepsilon \exp(\frac{bt}{2})+4\sqrt{x}\right)}{2}}\right)\right), \\
G_{4}: & \left(\frac{\left(2\exp(\frac{bt}{2})\sqrt{x}+\varepsilon\right)^{2}\exp(-bt)}{4},y, t,u\right), \\
G_{5}: & \left(x,\frac{\left(\varepsilon \exp(\frac{et}{2})+2\sqrt{y}\right)^{2}}{4}, t,u\exp\left(^{\frac{ \exp(\frac{et}{2})e\varepsilon\left(\varepsilon \exp(\frac{et}{2})+4\sqrt{y}\right)}{2}}\right)\right), \\
G_{6}: & \left(x,\frac{\left(2\exp(\frac{et}{2})\sqrt{y}+\varepsilon\right)^{2}\exp(-et)}{4}, t,u\right).
\end{align*}
\end{doublespace}

\noindent Because each one-parameter group $G_{i}$ is a symmetry group, using Theorem \ref{newSol}, if $u(x,y, t)$ is a solution of (\ref{I}), so are the $u^{\varepsilon,v_{i}},\,1\leq i\leq 6$, given by:

\begin{doublespace}
\begin{align*}\ds
u^{\varepsilon,v_{1}}(x,y,t)&=u(x,y, t-\varepsilon), \\
u^{\varepsilon,v_{2}}(x,y,t)&=\exp(\varepsilon) u(x,y, t), \\
u^{\varepsilon,v_{3}}(x,y,t)&=\exp\left({\frac{b\varepsilon\exp(\frac{bt}{2})\left(4\sqrt{x}-\varepsilon \exp(\frac{bt}{2})\right)}{2}}\right) u\left(\frac{\left(2\sqrt{x}-\epsilon \exp(\frac{bt}{2})\right)^{2}}{4},y,t\right), \\
u^{\varepsilon,v_{4}}(x,y,t)&=u\left(\frac{\left(2\exp(\frac{bt}{2})\sqrt{x}-\varepsilon\right)^{2}\exp(-bt)}{4},y,t\right), \\
u^{\varepsilon,v_{5}}(x,y,t)&=\exp\left(^{\frac{e\varepsilon \exp(\frac{et}{2})\left(4\sqrt{y}-\varepsilon \exp(\frac{et}{2})\right)}{2}}\right)u\left(\frac{\left(2\sqrt{y}-\epsilon \exp(\frac{et}{2})\right)^{2}}{4},y,t\right),\\
u^{\varepsilon,v_{6}}(x,y,t)&=u\left(x,\frac{\left(2\exp(\frac{et}{2})\sqrt{y}-\epsilon\right)^{2}\exp(-et)}{4},t\right). 
\end{align*}
\end{doublespace}

 \subsubsection{Subcase 1.2 ($b=e\,,\,b\ne 0\,,\,a=\frac{1}{4},d=\frac{1}{4}$)}
{\small
$$
\left\{\begin{aligned}
\xi&=(c_{1}b\exp(bt)-c_{2}b\exp(-bt))x+\sqrt{x}\left(c_{4}\sqrt{y}+c_{5}\exp(\frac{b}{2}t)+c_{6}\exp(-\frac{b}{2}t)\right),\\
\gamma&=(c_{1}b\exp(bt)-c_{2}b\exp(-bt))y+\sqrt{y}\left(-c_{4}\sqrt{x}+c_{7}\exp(\frac{b}{2}t)+c_{8}\exp(-\frac{b}{2}t)\right),\\
\tau&=c_{1}\exp(bt)+c_{2}\exp(-bt)+c_{3},\\
\phi&=\left((2b^{2}c_{1}\exp(bt))(x+y)+2bc_{5}\exp(\frac{b}{2}t)\sqrt{x}+2bc_{7}\exp(\frac{b}{2}t)\sqrt{y}-bc_{1}\exp(bt)-\frac{1}{2}bc_{3}+c_{9}\right)u+\mu(x,y,t).
\end{aligned}\right.
$$
}
{\small
\begin{doublespace}
\begin{align*}\ds
    v_{1}&=b\exp(bt)x\frac{\partial}{\partial x}+b\exp(bt)y\frac{\partial}{\partial y}+\exp(bt)\frac{\partial}{\partial t}+2b\exp(bt)(b(x+y)-\frac{1}{2})u\frac{\partial}{\partial u},\\
     v_{2}&=-b\exp(-bt)x\frac{\partial}{\partial x}-b\exp(-bt)y\frac{\partial}{\partial y}+\exp(-bt)\frac{\partial}{\partial t},\,\,v_{3}=\frac{\partial}{\partial t}-\frac{b}{2}u\frac{\partial}{\partial u}\,\,,\,\, v_{4}=\sqrt{x}\sqrt{y}\frac{\partial}{\partial x}-\sqrt{x}\sqrt{y}\frac{\partial}{\partial y},\\
      v_{5}&=\sqrt{x}\exp(\frac{b}{2}t)\frac{\partial}{\partial x}+2b\sqrt{x}\exp(\frac{b}{2}t)u\frac{\partial}{\partial u}\,\,,\,\,v_{6}=\sqrt{x}\exp(-\frac{b}{2}t)\frac{\partial}{\partial x},\\
       v_{7}&=\sqrt{y}\exp(\frac{b}{2}t)\frac{\partial}{\partial y}+2b\sqrt{y}\exp(\frac{b}{2}t)u\frac{\partial}{\partial u}\,\,,\,\, v_{8}=\sqrt{y}\exp(-\frac{b}{2}t)\frac{\partial}{\partial y}\,\,,\,\, v_{9}=u\frac{\partial}{\partial u}.
\end{align*}
\end{doublespace}
}
{\small
\begin{longtable}{|c|c|c|c|c|c|c|c|c|c|c|} 
\caption{Commutator table($b=e\,\,,\,b\ne0\,,\,a=\frac{1}{4},d= \frac{1}{4}$).}\\
\hline
$[v_{i},v_{j}]$ & $v_{1}$ & $v_{2}$&$v_{3}$&$v_{4}$& $v_{5}$&$v_{6}$&$v_{7}$&$v_{8}$&$v_{9}$&$v_{\mu}$ \endhead
\hline
$v_{1}$&0&$-2bv_{3}$&$-bv_{1}$&0&0&$-bv_{5}$&0&$-bv_{7}$&0&$v_{\mu_{1}}$\\ 
$v_{2}$&$2bv_{3}$&0&$bv_{2}$&0&$bv_{6}$&0&$bv_{8}$&0&0&$v_{\mu_{2}}$\\
$v_{3}$&$bv_{1}$&$-bv_{2}$&0&0&$\frac{b}{2}v_{5}$&$-\frac{b}{2}v_{6}$&$\frac{b}{2}v_{7}$&$-\frac{b}{2}v_{8}$&0&$v_{\mu_{t}}$\\
$v_{4}$&0&0&0&0&$\frac{1}{2}v_{7}$&$\frac{1}{2}v_{8}$&$-\frac{1}{2}v_{5}$&$-\frac{1}{2}v_{6}$&0&$v_{\mu_{3}}$\\
$v_{5}$&0&$-bv_{6}$&$-\frac{b}{2}v_{5}$&$-\frac{1}{2}v_{7}$&0&$-bv_{9}$&0&0&0&$v_{\mu_{4}}$\\ 
$v_{6}$&$bv_{5}$&0&$\frac{b}{2}v_{6}$&$-\frac{1}{2}v_{8}$&$bv_{9}$&0&0&0&0&$v_{\mu_{5}}$\\ 
$v_{7}$&0&$-bv_{8}$&$-\frac{b}{2}v_{7}$&$\frac{1}{2}v_{5}$&0&0&0&$-bv_{9}$&0&$v_{\mu^{6}}$\\ 
$v_{8}$&$bv_{7}$&0&$\frac{b}{2}v_{8}$&$\frac{1}{2}v_{6}$&0&0&$bv_{9}$&0&0&$v_{\mu_{7}}$\\ 
$v_{9}$&0&0&0&0&0&0&0&0&0&$-v_{\mu}$\\ 
$v_{\mu}$&$-v_{\mu_{1}}$&$-v_{\mu_{2}}$&$-v_{\mu_{t}}$&$-v_{\mu_{3}}$&$-v_{\mu_{4}}$&$-v_{\mu_{5}}$&$-v_{\mu_{6}}$&$-v_{\mu_{7}}$&$v_{\mu}$&0\\ 
\hline
\end{longtable}}
{\small
\begin{doublespace}
\begin{align*}\ds   
\mu_{1}&=b\exp(bt)x\mu_{x}+b\exp(bt)y\mu_{y}+\exp(bt)\mu_{t}+b\exp(bt)(1-2b(x+y)))\mu,\\
\mu_{2}&=-b\exp(-bt)x\mu_{x}-b\exp(-bt)y\mu_{y}+\exp(bt)\mu_{t},\\
\mu_{3}&=\sqrt{x}\sqrt{y}\mu_{x}-\sqrt{x}\sqrt{y}\mu_{y}\,\,,\,\, \mu_{4}=\sqrt{x}\exp(\frac{b}{2}t)\mu_{x}-2b\exp(\frac{b}{2}t)\sqrt{x}\mu,\\
 \mu_{5}&=\sqrt{x}\exp(-\frac{b}{2}t)\mu_{x}\,\,,\,\, \mu_{6}=\sqrt{y}\exp(\frac{b}{2}t)\mu_{y}-2b\exp(\frac{b}{2}t)\sqrt{y}\mu\,\,,\,\,\mu_{7}=\sqrt{y}\exp(-\frac{b}{2}t)\mu_{y}.
\end{align*}
\end{doublespace}
}
\noindent The transformed point by a one-parameter group $\exp \left(\varepsilon v_{i}\right)(x,y, t, u)=$ $(\tilde{x},\tilde{y}, \tilde{t}, \tilde{u})$ are:
{\small
\begin{doublespace}
\begin{align*}\ds
G_{1}: & \left(\frac{x}{1-b\varepsilon \exp(bt)},\frac{y}{1-b\varepsilon \exp(bt)}, t-\frac{\ln(1-b\varepsilon \exp(bt))}{b}, (1-b\varepsilon \exp(bt))u\exp\left(\frac{2b^{2}\varepsilon \exp(bt)(x+y)}{1-b\varepsilon \exp(bt)}\right)\right), \\
G_{2}: & \left(\frac{x}{1+b\varepsilon \exp(-bt)},\frac{y}{1+b\varepsilon \exp(-bt)},t+ \frac{\ln(1+b\varepsilon \exp(-bt))}{b},u\right), \\
G_{3}: & (x,y,t+\varepsilon,u), \\
G_{4}: & \left(\sqrt{x}\sqrt{y}\sin(\varepsilon)+(\frac{x-y}{2})\cos(\varepsilon)+(\frac{x+y}{2}),-\sqrt{x}\sqrt{y}\sin(\varepsilon)-(\frac{x-y}{2})\cos(\varepsilon)-(\frac{x+y}{2}),t,u\right),\\
G_{5}: & \left(\frac{\left(\varepsilon \exp(\frac{bt}{2})+2\sqrt{x}\right)^{2}}{4},y, t,u\exp\left(\frac{ \exp(\frac{bt}{2})b\varepsilon\left(\varepsilon \exp(\frac{bt}{2})+4\sqrt{x}\right)}{2}\right)\right), \\
G_{6}: & \left(\frac{\left(2\exp(\frac{bt}{2})\sqrt{x}+\varepsilon\right)^{2}\exp(-bt)}{4},y, t,u\right), \\
G_{7}: & \left(x,\frac{\left(\varepsilon \exp(\frac{bt}{2})+2\sqrt{y}\right)^{2}}{4}, t,u\exp\left(\frac{ \exp(\frac{bt}{2})b\varepsilon\left(\varepsilon \exp(\frac{bt}{2})+4\sqrt{y}\right)}{2}\right)\right), \\
G_{8}: & \left(x,\frac{\left(2\exp(\frac{bt}{2})\sqrt{y}+\varepsilon\right)^{2}\exp(-bt)}{4},t,u\right), \\
G_{9}: & (x,y,t,ue^{\varepsilon}).
\end{align*}
\end{doublespace}
}
\noindent Since each one parameter group $G_{i}$ is a symmetry group, using Theorem \ref{newSol} if $u(x,y, t)$ is a solution of (\ref{I}), so are the $u^{\varepsilon,v_{i}},\,1\leq i\leq 9$ given by:
{\small
\begin{doublespace}
\begin{align*}\ds
u^{\varepsilon,v_{1}}(x,y,t)&=\frac{1}{1+b\epsilon \exp(bt)}\exp\left(2b^{2}\varepsilon \exp(bt)(x+y)\right)u\left(\frac{x}{1+b\epsilon \exp(bt)},\frac{x}{1+b\epsilon \exp(bt)},t-\frac{\ln(1+\varepsilon b \exp(bt) )}{b}\right), \\
u^{\varepsilon,v_{2}}(x,y,t)&=u\left(\frac{x}{1-b\varepsilon \exp(-bt)},\frac{y}{1-b\varepsilon \exp(-bt)},t+\frac{\ln(1-b\varepsilon \exp(-bt))}{b}\right), \\
u^{\varepsilon,v_{3}}(x,y,t)&=u(x,y,t-\varepsilon), \\
u^{\varepsilon,v_{4}}(x,y,t)&=u\left(-\sqrt{x}\sqrt{y}\sin(\varepsilon)+(\frac{x-y}{2})\cos(\varepsilon)+(\frac{x+y}{2}),\sqrt{x}\sqrt{y}\sin(\varepsilon)-(\frac{x-y}{2})\cos(\varepsilon)-(\frac{x+y}{2}),t\right),\\
u^{\varepsilon,v_{5}}(x,y,t)&=\exp\left(\frac{b\epsilon \exp(\frac{bt}{2})(4\sqrt{x}-\epsilon \exp(\frac{bt}{2}))}{2}\right)u\left(\frac{\left(2\sqrt{x}-\epsilon \exp(\frac{bt}{2})\right)^{2}}{4},y,t\right),  \\
u^{\varepsilon,v_{6}}(x,y,t)&=u\left(\frac{\left(2\exp(\frac{bt}{2})\sqrt{x}-\epsilon\right)^{2}}{4}\exp(-bt),y,t\right), \\
u^{\varepsilon,v_{7}}(x,y,t)&=\exp\left(\frac{b\epsilon \exp(\frac{bt}{2})(4\sqrt{y}-\epsilon \exp(\frac{bt}{2}))}{2}\right)u\left(\frac{\left(2\sqrt{y}-\epsilon \exp(\frac{bt}{2})\right)^{2}}{4},y,t\right), \\
u^{\varepsilon,v_{8}}(x,y,t)&=u\left(x,\frac{\left(2\exp(\frac{bt}{2})\sqrt{y}-\epsilon\right)^{2}\exp(-bt)}{4},t\right), \\
u^{\varepsilon,v_{9}}(x,y,t)&=\exp(\varepsilon) u(x,y,t). 
\end{align*}
\end{doublespace}
}
\subsubsection{Subcase 1.3$(b=-e\,,\,b\ne 0\,,\,a=\frac{1}{4},d=\frac{1}{4}$)}
{\small
\begin{align*}
    v_{1}&=b\exp(bt)x\frac{\partial}{\partial x}+b\exp(bt)y\frac{\partial}{\partial y}+\exp(bt)\frac{\partial}{\partial t}+\frac{b\exp(bt)(4bx-1)}{2}u\frac{\partial}{\partial u},\\
     v_{2}&=-b\exp(-bt)x\frac{\partial}{\partial x}-b\exp(-bt)y\frac{\partial}{\partial y}+\exp(-bt)\frac{\partial}{\partial t}+\frac{b\exp(-bt)(4by+1)}{2}u\frac{\partial}{\partial u},\\
       v_{3}&=\frac{\partial}{\partial t}\,\,,\,\,v_{4}=\sqrt{x}\sqrt{y}\frac{\partial}{\partial x}-\sqrt{x}\sqrt{y}\frac{\partial}{\partial y}+2b\sqrt{x}\sqrt{y}u\frac{\partial}{\partial u},\\
       v_{5}&=\sqrt{x}\exp(\frac{b}{2}t)\frac{\partial}{\partial x}+2b\sqrt{x}\exp(\frac{b}{2}t)u\frac{\partial}{\partial u}\\
     v_{6}&=\sqrt{x}\exp(-\frac{b}{2}t)\frac{\partial}{\partial x},\,\,v_{7}=\sqrt{y}\exp(-\frac{b}{2}t)\frac{\partial}{\partial y}-2b\sqrt{y}\exp(-\frac{b}{2}t)u\frac{\partial}{\partial u},\\
      v_{8}&=\sqrt{y}\exp(\frac{b}{2}t)\frac{\partial}{\partial y}\,\,,\,\, v_{9}=u\frac{\partial}{\partial u}.
\end{align*}
}
\begin{longtable}{|c|c|c|c|c|c|c|c|c|c|c|} 
\caption{Commutator table($b=-e\,,\,b\ne0\,,\,a=\frac{1}{4},d= \frac{1}{4}$).}\\
\hline
 $[v_{i},v_{j}]$& $v_{1}$ & $v_{2}$&$v_{3}$&$v_{4}$& $v_{5}$&$v_{6}$&$v_{7}$&$v_{8}$&$v_{9}$&$v_{\mu}$ \endhead
\hline
$v_{1}$&0&$-2bv_{3}$&$-bv_{1}$&0&0&$-bv_{5}$&$-bv_{8}$&0&0&$v_{\mu_{1}}$\\ 
$v_{2}$&$2bv_{3}$&0&$bv_{2}$&0&$bv_{6}$&0&0&$bv_{7}$&0&$v_{\mu_{2}}$\\
$v_{3}$&$bv_{1}$&$-bv_{2}$&0&0&$\frac{b}{2}v_{5}$&$-\frac{b}{2}v_{6}$&$-\frac{b}{2}v_{7}$&$\frac{b}{2}v_{8}$&0&$v_{\mu_{t}}$\\
$v_{4}$&0&0&0&0&$\frac{1}{2}v_{8}$&$\frac{1}{2}v_{7}$&$-\frac{b}{2}v_{6}$&$-\frac{1}{2}v_{5}$&0&$v_{\mu_{3}}$\\
$v_{5}$&0&$-bv_{6}$&$-\frac{b}{2}v_{5}$&$-\frac{1}{2}v_{8}$&0&$-bv_{9}$&0&0&0&$v_{\mu_{4}}$\\ 
$v_{6}$&$bv_{5}$&0&$\frac{b}{2}v_{6}$&$-\frac{1}{2}v_{7}$&$bv_{9}$&0&0&0&0&$v_{\mu_{5}}$\\ 
$v_{7}$&$bv_{8}$&0&$\frac{b}{2}v_{7}$&$\frac{1}{2}v_{6}$&0&0&0&$bv_{9}$&0&$v_{\mu_{6}}$\\ 
$v_{8}$&0&$-bv_{7}$&$-\frac{b}{2}v_{8}$&$\frac{1}{2}v_{5}$&0&0&$-bv_{9}$&0&0&$v_{\mu_{7}}$\\ 
$v_{9}$&0&0&0&0&0&0&0&0&0&$-v_{\mu}$\\ 
$v_{\mu}$&$-v_{\mu_{1}}$&$-v_{\mu_{2}}$&$-v_{\mu_{t}}$&$-v_{\mu_{3}}$&$-v_{\mu_{4}}$&$-v_{\mu_{5}}$&$-v_{\mu_{6}}$&$-v_{\mu_{7}}$&$v_{\mu}$&0\\ 
\hline
\end{longtable}
{\small
\begin{doublespace}
\begin{align*}\ds
     \mu_{1}&=b\exp(bt)x\mu_{x}+b\exp(bt)y\mu_{y}+\exp(bt)\mu_{t}+\frac{1}{2}b(1-4b)\exp(bt)\mu,\\
     \mu_{2}&=-b\exp(-bt)x\mu_{x}-b\exp(-bt)y\mu_{y}+\exp(-bt)\mu_{t}-\frac{1}{2}b(1+by))\exp(bt)\mu,\\
     \mu_{3}&=\sqrt{x}\sqrt{y}\mu_{x}-\sqrt{x}\sqrt{y}\mu_{y}-2b\sqrt{x}\sqrt{y}\mu\,\,,\,\, \mu_{4}=\sqrt{x}\exp(\frac{b}{2}t)\mu_{x}-2b\exp(\frac{b}{2}t)\sqrt{x}\mu,\\
      \mu_{5}&=\sqrt{x}\exp(-\frac{b}{2}t)\mu_{x}\,\,,\,\,\mu_{6}=\sqrt{y}\exp(-\frac{b}{2}t)\mu_{y}-2b\exp(-\frac{b}{2}t)\sqrt{y}\mu\,\,,\,\, \mu_{7}=\sqrt{y}\exp(-\frac{b}{2}t)\mu_{y}.   
\end{align*}
\end{doublespace}
}
\noindent The transformed point by a one-parameter group $\exp \left(\varepsilon v_{i}\right)(x,y, t, u)=$ $(\tilde{x},\tilde{y}, \tilde{t}, \tilde{u})$ are:
{\small
\begin{doublespace}
\begin{align*}\ds
\nonumber G_{1}: & \left(\frac{x}{1-b\varepsilon \exp(bt)},\frac{y}{1-b\varepsilon \exp(bt)}, t-\frac{\ln(1-b\varepsilon \exp(bt))}{b}, \sqrt{1-b\varepsilon \exp(bt)}u\exp\left(\frac{2b^{2}\varepsilon \exp(bt)x}{1-b\varepsilon \exp(bt)}\right)\right), \\
\nonumber G_{2}: & ~ \Bigg(\frac{x}{1+b\varepsilon \exp(-bt)},\frac{y}{1+b\varepsilon \exp(-bt)}, t+\frac{\ln(1-b\varepsilon \exp(-bt))}{b},\sqrt{1+b\varepsilon \exp(-bt)}u\exp\left(\frac{2b^{2}\varepsilon \exp(-bt)y}{1+b\varepsilon \exp(-bt)}\right)\Bigg), \\
\nonumber G_{3}: & ~(x,y,t+\varepsilon,u), \\
\nonumber G_{4}: &~ \Bigg(\sqrt{x}\sqrt{y}\sin(\varepsilon)+(\frac{x-y}{2})\cos(\varepsilon)+(\frac{x+y}{2}),-\sqrt{x}\sqrt{y}\sin(\varepsilon)-(\frac{x-y}{2})\cos(\varepsilon)+(\frac{x+y}{2}),t,\\\nonumber
& u\exp\left(b(2\sqrt{x}\sqrt{y}\sin(\varepsilon)+(x-y)\cos(\varepsilon)-(x-y))\right)\Bigg), \\
\nonumber G_{5}: & \left(\frac{\left(\varepsilon \exp(\frac{bt}{2})+2\sqrt{x}\right)^{2}}{4},y, t,u\exp\left(\frac{ \exp(\frac{bt}{2})b\varepsilon\left(\varepsilon \exp(\frac{bt}{2})+4\sqrt{x}\right)}{2}\right)\right), \\
\nonumber G_{6}: & \left(\frac{\left(2\exp(\frac{bt}{2})\sqrt{x}+\varepsilon\right)^{2}\exp(-bt)}{4},y, t,u\right), \\
\nonumber G_{7}: & \left(x,\frac{\left(\varepsilon \exp(\frac{et}{2})+2\sqrt{y}\right)^{2}}{4}, t,u\exp\left(\frac{ \exp(\frac{et}{2})e\varepsilon\left(\varepsilon \exp(\frac{et}{2})+4\sqrt{y}\right)}{2}\right)\right), \\
\nonumber G_{8}: & \left(x,\frac{\left(2\sqrt{y}+\varepsilon \exp(\frac{bt}{2})\right)^{2}}{4},t,u\right), \\
\nonumber G_{9}: & (x,y,t,u\exp(\varepsilon)).
\end{align*}
\end{doublespace}
}
\noindent Since each one parameter group $G_{i}$ is a symmetry group, using Theorem \ref{newSol} if $u(x,y, t)$ is a solution of (\ref{I}), so are the $u^{\varepsilon,v_{i}},\,1\leq i\leq 9$ given by: 
{\small
\begin{doublespace}
\begin{align*}\ds
u^{\varepsilon,v_{1}}(x,y,t)&=\frac{\exp\left(\frac{2b^{2}\varepsilon x \exp(bt) }{1+b\varepsilon \exp(bt)}\right)}{\sqrt{1+b\varepsilon \exp(bt)}}u\left(\frac{x}{1+b\varepsilon \exp(bt)},\frac{y}{1+b\varepsilon \exp(bt)},t-\frac{\ln(1+\varepsilon b \exp(bt) )}{b}\right), \\
u^{\varepsilon,v_{2}}(x,y,t)&=\sqrt{1+\frac{b\epsilon \exp(-bt)}{1+b\varepsilon \exp(-bt)}}\exp\left(\frac{2b^{2}\varepsilon \exp(-bt)y}{ (1-b\varepsilon \exp(-bt))(1+2b\epsilon \exp(-bt))}\right)\\
&u\left(\frac{x}{(1-\varepsilon b \exp(-bt))},\frac{y}{(1-\varepsilon b \exp(-bt))},t+\frac{\ln(1+b\varepsilon \exp(-bt))}{b}\right), \\
u^{\varepsilon,v_{3}}(x,y,t)&=u(x,y,t-\varepsilon), \\
u^{\varepsilon,v_{4}}(x,y,t)&=
\exp~ \Bigg[b~ \Bigg(\sqrt{\left(-\sqrt{xy}\sin(\varepsilon)+\frac{(x-y)}{2}\cos(\varepsilon)+\frac{(x+y)}{2}\right)\left(\sqrt{xy}\sin(\varepsilon)-\frac{(x-y)}{2}\cos(\varepsilon)+\frac{(x+y)}{2}\right)}\sin(\varepsilon)\\
  &+(-2\sqrt{x}\sqrt{y}\sin(\varepsilon)+(x-y)\cos(\varepsilon))(\cos(\varepsilon)-1)\Bigg)\Bigg]\\
  &u\left(-\sqrt{x}\sqrt{y}\sin(\varepsilon)+(\frac{x-y}{2})\cos(\varepsilon)+(\frac{x+y}{2}),\sqrt{x}\sqrt{y}\sin(\varepsilon)-(\frac{x-y}{2})\cos(\varepsilon)-(\frac{x+y}{2}),t\right),\\
u^{\varepsilon,v_{5}}(x,y,t)&=\exp\left(\frac{b\varepsilon \exp(\frac{bt}{2})\left(4\sqrt{x}-\varepsilon \exp(\frac{bt}{2})\right)}{2}\right)u\left(\frac{\left(2\sqrt{x}-\varepsilon \exp(\frac{bt}{2})\right)^{2}}{4},y,t\right), \\
u^{\varepsilon,v_{6}}(x,y,t)&=u\left(\frac{\left(2\exp(\frac{bt}{2})\sqrt{x}-\epsilon\right)^{2}}{4}\exp(-bt),y,t\right), \\
u^{\varepsilon,v_{7}}(x,y,t)&=\exp\left(-\frac{b\varepsilon e^{-\frac{bt}{2}}\left(4\sqrt{y}-\varepsilon \exp(-\frac{bt}{2})\right)}{2}\right)u\left(x,\frac{\left(2\sqrt{y}-\varepsilon \exp(-\frac{bt}{2})\right)^{2}}{4},t\right), \\
u^{\varepsilon,v_{8}}(x,y,t)&=u\left(x,\frac{\left(2\sqrt{y}-\epsilon \exp(\frac{bt}{2})\right)^{2}}{4},t\right), \\
u^{\varepsilon,v_{9}}(x,y,t)&=\exp(\varepsilon) u(x,y, t).
\end{align*}
\end{doublespace}
}
\subsubsection{Subcase 1.4$(b=e=0\,,\,a=\frac{1}{4},d=\frac{1}{4})$}
{\small
\begin{doublespace}
\begin{align*}
    v_{1}&=2tx\frac{\partial}{\partial x}+2ty\frac{\partial}{\partial y}+t^{2}\frac{\partial}{\partial t}+(2(x+y)-\frac{t}{2})u\frac{\partial}{\partial u},\\
    v_{2}&=x\frac{\partial}{\partial x}+y\frac{\partial}{\partial y}+t\frac{\partial}{\partial t},\\
    v_{3}&=\frac{\partial}{\partial t},\,\,v_{4}=\sqrt{x}\sqrt{y}\frac{\partial}{\partial x}-\sqrt{x}\sqrt{y}\frac{\partial}{\partial y},\,\,
    v_{5}=\sqrt{x}t\frac{\partial}{\partial x}+2\sqrt{x}u\frac{\partial}{\partial u},\\
    v_{6}&=\sqrt{x}\frac{\partial}{\partial x},\,\,v_{7}=\sqrt{y}t\frac{\partial}{\partial y}+2\sqrt{y}u\frac{\partial}{\partial u},\,\, v_{8}=\sqrt{y}\frac{\partial}{\partial y}\,\,,\,\,v_{9}=u\frac{\partial}{\partial u}.
\end{align*}
\end{doublespace}
}
\begin{remark}
    $\{v_{1},...,v_{9}\}$ is the largest Lie algebra of symmetries that we found.
\end{remark}

{\small
\begin{longtable}{|c|c|c|c|c|c|c|c|c|c|c|} 
\caption{Commutator table($b=e=0\,,\,a=\frac{1}{4},d= \frac{1}{4}$).}\\
\hline
 $[v_{i},v_{j}]$& $v_{1}$ & $v_{2}$&$v_{3}$&$v_{4}$& $v_{5}$&$v_{6}$&$v_{7}$&$v_{8}$&$v_{9}$&$v_{\mu}$ \endhead
\hline
$v_{1}$&0&$-v_{1}$&$-2v_{2}+\frac{1}{2}v_{9}$&0&0&$-v_{5}$&$0$&$-v_{7}$&0&$v_{\mu_{1}}$\\ 
$v_{2}$&$v_{1}$&0&$-v_{3}$&0&$\frac{1}{2}v_{5}$&$-\frac{1}{2}v_{6}$&$\frac{1}{2}v_{7}$&$-\frac{1}{2}v_{8}$&0&$v_{\mu_{2}}$\\
$v_{3}$&$2v_{2}-\frac{1}{2}v_{9}$&$v_{3}$&0&0&$v_{6}$&$0$&$v_{8}$&$0$&0&$v_{\mu_{3}}$\\
$v_{4}$&0&0&0&0&$\frac{1}{2}v_{7}$&$\frac{1}{2}v_{8}$&$-\frac{1}{2}v_{5}$&$-\frac{1}{2}v_{6}$&0&$v_{\mu_{3}}$\\
$v_{5}$&$0$&$-\frac{1}{2}v_{5}$&$-v_{6}$&$-\frac{1}{2}v_{7}$&0&$-v_{9}$&0&0&0&$v_{\mu_{4}}$\\ 
$v_{6}$&$0$&$\frac{1}{2}v_{6}$&$0$&$-\frac{1}{2}v_{8}$&$v_{9}$&0&0&0&0&$v_{\mu_{5}}$\\ 
$v_{7}$&$0$&$-\frac{1}{2}v_{7}$&$-v_{8}$&$\frac{1}{2}v_{5}$&0&0&0&$-v_{9}$&0&$v_{\mu_{6}}$\\ 
$v_{8}$&$-v_{7}$&$\frac{1}{2}v_{8}$&$0$&$\frac{1}{2}v_{6}$&0&0&$v_{9}$&0&0&$v_{\mu_{7}}$\\ 
$v_{9}$&0&0&0&0&0&0&0&0&0&$-v_{\mu}$\\ 
$v_{\mu}$&$-v_{\mu_{1}}$&$-v_{\mu_{2}}$&$-v_{\mu_{3}}$&$-v_{\mu_{3}}$&$-v_{\mu_{4}}$&$-v_{\mu_{5}}$&$-v_{\mu_{6}}$&$-v_{\mu_{7}}$&$v_{\mu}$&0\\ 
\hline
\end{longtable}}
{\small
\begin{doublespace}
\begin{align*}\ds
     \mu_{1}&=2tx\mu_{x}+2ty\mu_{y}+t^{2}\mu_{t}+(\frac{t}{2}-2(x+y))\mu,\,\,\mu_{2}=x\mu_{x}+y\mu_{y}+t\mu_{t},\,\,\mu_{3}=\mu_{t},\\
       \mu_{4}&=\sqrt{x}\sqrt{y}\mu_{x}-\sqrt{x}\sqrt{y}\mu_{y},\,\,\mu_{5}=\sqrt{x}t\mu_{x}-2\sqrt{x}\mu\,\,,\,\,\mu_{6}=\sqrt{x}\mu_{x}\,\,,\,\, \mu_{7}=\sqrt{y}t\mu_{y}-2\sqrt{y}\mu.
\end{align*}
\end{doublespace}
}
\noindent The transformed pointby a one-parameter group $\exp \left(\varepsilon v_{i}\right)(x,y, t, u)=$ $(\tilde{x},\tilde{y}, \tilde{t}, \tilde{u})$ are:
{\small
\begin{doublespace}
\begin{align*}\ds
\nonumber G_{1}: & \left(\frac{x}{(1-t\varepsilon)^{2}},\frac{y}{(1-t\varepsilon)^{2}},\frac{t}{1-t\varepsilon},u\exp\left(\frac{2\varepsilon(x+y)}{1-t\varepsilon}\right)\sqrt{1-t\varepsilon}\right), \\
\nonumber G_{2}: & (\exp(\varepsilon)x,\exp(\varepsilon)y,\exp(\varepsilon)t,u), \\
\nonumber G_{3}: & ~(x,y,t+\varepsilon,u), \\
\nonumber G_{4}: &\left(\sqrt{x}\sqrt{y}\sin(\varepsilon)+\left(\frac{x-y}{2}\right)\cos(\varepsilon)+\left(\frac{x+y}{2}\right),-\sqrt{x}\sqrt{y}\sin(\varepsilon)-\left(\frac{x-y}{2}\right)\cos(\varepsilon)+\left(\frac{x+y}{2}\right),t,u\right), \\
\nonumber G_{5}: & \left(x+t\sqrt{x}\varepsilon +\frac{t^{2}\varepsilon^{2}}{4},y,t,u\exp\left(\frac{\varepsilon(t\varepsilon+4\sqrt{x})}{2}\right)\right), \\
\nonumber G_{6}: & \left(\sqrt{x}\varepsilon+\frac{\varepsilon^{2}}{4}+x,y,t,u\right), \\
\nonumber G_{7}: & \left(x,y+t\sqrt{y}\varepsilon+\frac{t^{2}}{4}\varepsilon^{2},t,u\exp\left(\frac{\varepsilon(t\varepsilon+4\sqrt{y})}{2}\right)\right), \\
\nonumber G_{8}: & (x,\sqrt{y}\varepsilon+\frac{\varepsilon^{2}}{4}+y,t,u),\\
\nonumber G_{9}: & (x,y,t,\exp(\varepsilon)u).
\end{align*}
\end{doublespace}
}
\noindent Since each one parameter group $G_{i}$ is a symmetry group, using Theorem \ref{newSol} if $u(x,y, t)$ is a solution of (\ref{I}), so are the $u^{\varepsilon,v_{i}},\,1\leq i\leq 9$ given by: 
{\small
\begin{doublespace}
\begin{align*}\ds
u^{\varepsilon,v_{1}}(x,y,t)&=\frac{\exp\left(\frac{2\varepsilon(x+y)}{(1+t\varepsilon)^{3}}\right)}{\sqrt{1+t\varepsilon}}u\left(\frac{x}{(1+t\varepsilon)^{2}},\frac{y}{(1+t\varepsilon)^{2}},\frac{t}{(1+t\varepsilon)}\right), \\
u^{\varepsilon,v_{2}}(x,y,t)&=u\left(\exp(-\varepsilon)x,\exp(-\varepsilon)y,\exp(-\varepsilon)t\right),\\
u^{\varepsilon,v_{3}}(x,y,t)&=u(x,y,t-\varepsilon), \\
u^{\varepsilon,v_{4}}(x,y,t)&=u~\Bigg(-\sqrt{x}\sqrt{y}\sin(\varepsilon)+\left(\frac{x-y}{2}\right)\cos(\varepsilon)+\left(\frac{x+y}{2}\right),\sqrt{x}\sqrt{y}\sin(\varepsilon)-\left(\frac{x-y}{2}\right)\cos(\varepsilon)-\left(\frac{x+y}{2}\right)
,t\Bigg),\\
u^{\varepsilon,v_{5}}(x,y,t)&=\exp\left(\frac{\varepsilon \left(\varepsilon t+4\sqrt{x-t\sqrt{x}\varepsilon+\frac{t^{2}\varepsilon^{2}}{2}}\right)}{2}\right)u\left(x-t\sqrt{x}\varepsilon+\frac{t^{2}\varepsilon^{2}}{2},y,t\right), \\
u^{\varepsilon,v_{6}}(x,y,t)&=u\left(\sqrt{x}\varepsilon-\frac{\varepsilon^{2}}{4}+x,y,t\right), \\
u^{\varepsilon,v_{7}}(x,y,t)&=\exp\left(\frac{\varepsilon \left(\varepsilon t+4\sqrt{y-t\sqrt{y}\varepsilon+\frac{t^{2}\varepsilon^{2}}{2}}\right)}{2}\right)u\left(x,y-t\sqrt{y}\varepsilon+\frac{t^{2}\varepsilon^{2}}{2},t\right), \\
u^{\varepsilon,v_{8}}(x,y,t)&=u\left(-\sqrt{x}\varepsilon+\frac{\varepsilon^{2}}{4}+x,y,t\right), \\
u^{\varepsilon,v_{9}}(x,y,t)&=\exp(\varepsilon) u(x,y, t).
\end{align*}
\end{doublespace}
}
\subsection{Case 2 ($a\ne\frac{1}{4},d\ne\frac{1}{4}$)}
Because $a\ne\frac{1}{4},\,\,d\ne\frac{1}{4}$ and from remark \eqref{h,k,l}, then $h=k=l=0$ in all the following subcases.
\subsubsection{Subcase 2.1 ($b\ne\pm e\,,\,a\ne\frac{1}{4},d\ne \frac{1}{4}$)}

$$
\tau=c_{1},\,\,\xi=\gamma=0,\,\,
\phi=c_{2}u+\mu(x,y,t).
$$

\begin{align*}
    v_{1}&=\frac{\partial}{\partial t}\,\,,\,\, v_{2}=u\frac{\partial}{\partial u}
\end{align*}
\begin{remark}
$\{v_{1},v_{2}\}$ is the smallest Lie algebra of infinitesimal symmetries that we found.  
\end{remark}
{\small
\begin{longtable}{|c|c|c|c|} 
\caption{Commutator table($b\ne\pm e\,,\,a\ne\frac{1}{4},d\ne \frac{1}{4}$).}\\
\hline
$[v_{i},v_{j}]$ & $v_{1}$ & $v_{2}$& $v_{\mu}$ \endhead
\hline
$v_{1}$&$0$&$0$&$v_{\mu_{t}}$\\ 
$v_{2}$&$0$&$0$&$v_{\mu}$\\ 
$v_{\mu}$&$-v_{\mu_{t}}$&$-v_{\mu}$&$0$\\
\hline
\end{longtable}}
\noindent The transformed pointby a one-parameter group $\exp \left(\varepsilon v_{i}\right)(x,y, t, u)=$ $(\tilde{x},\tilde{y}, \tilde{t}, \tilde{u})$ are:

$$
G_{1}:  (x,y, t+\varepsilon, u),\,\,
G_{2}: (x,y, t, \exp(\varepsilon)u). 
$$

\noindent Since each one parameter group $G_{i}$ is a symmetry group, using Theorem \ref{newSol} if $u(x,y, t)$ is a solution of (\ref{I}), so are the $u^{\varepsilon,v_{i}},\,1\leq i\leq 2$ given by:

$$
u^{\varepsilon,v_{1}}(x,y,t)=f(x,y, t-\varepsilon),\,\,
u^{\varepsilon,v_{2}}(x,y,t)=\exp(\varepsilon) f(x,y, t).
$$

\subsubsection{Subcase 2.2($b=e\,,b\ne 0\,,\,a\ne\frac{1}{4},d\ne\frac{1}{4}$)}
{\small
\begin{doublespace}
\begin{align*}\ds
    v_{1}&=b\exp(bt)x\frac{\partial}{\partial x}+b\exp(bt)y\frac{\partial}{\partial y}+b\exp(bt)2(b(x+y)-(a+d))u\frac{\partial}{\partial u},\\
     v_{2}&=-b\exp(-bt)x\frac{\partial}{\partial x}-b\exp(-bt)y\frac{\partial}{\partial y}+\exp(-bt)\frac{\partial}{\partial t}, v_{3}=\frac{\partial}{\partial t}-b(a+d)u\frac{\partial}{\partial u}\,\,,\,\, v_{4}=u\frac{\partial}{\partial u}.
\end{align*}
\end{doublespace}
}
{\small
\begin{longtable}{|c|c|c|c|c|c|} 
\caption{Commutator table($b=e\,,\,b\ne0\,,\,a\ne\frac{1}{4},d\ne \frac{1}{4}$).}\\
\hline
$[v_{i},v_{j}]$ & $v_{1}$ & $v_{2}$&$v_{3}$&$v_{4}$& $v_{\mu}$ \endhead
\hline
$v_{1}$&0&$-2bv_{3}$&$-bv_{1}$&0&$v_{\mu_{1}}$\\ 
$v_{2}$&$2bv_{3}$&0&$bv_{2}$&0&$v_{\mu_{2}}$\\
$v_{3}$&$bv_{1}$&$-bv_{2}$&0&0&$v_{\mu_{3}}$\\
$v_{4}$&0&0&0&0&$-v_{\mu}$\\
$v_{\mu}$&$-v_{\mu_{1}}$&$-v_{\mu_{2}}$&$-v_{\mu_{3}}$&$v_{\mu}$&0\\
\hline
\end{longtable}
\begin{doublespace}
\begin{align*}\ds
    \mu_{1}&=b\exp(bt)x\mu_{x}+b\exp(bt)y\mu_{y}+\exp(bt)\mu_{t}+b(a+d-b(x+y))\exp(bt)\mu,\\
    \mu_{2}&=-b\exp(-bt)x\mu_{x}-b\exp(-bt)y\mu_{y}+\exp(-bt)\mu_{t}\,\,,\,\,\mu_{3}=b(a+d)\mu+\mu_{t}.
\end{align*}
\end{doublespace}}

\noindent The transformed pointby a one-parameter group $\exp \left(\varepsilon v_{i}\right)(x,y, t, u)=$ $(\tilde{x},\tilde{y}, \tilde{t}, \tilde{u})$ are:
{\small
\begin{doublespace}
\begin{align*}\ds
G_{1}: & ~ \Bigg(\frac{x}{1-\varepsilon b \exp(bt)},\frac{y}{1-\varepsilon b \exp(bt)},t-\frac{\ln(1-\varepsilon b \exp(bt))}{b},(1-\varepsilon b \exp(bt))^{2(d+a)}u\exp\left(\frac{2b^{2}\exp(bt)\varepsilon(x+y)}{1-\varepsilon b \exp(bt)}\right)\Bigg), \\
G_{2}: & \left(\frac{x}{1+\varepsilon b \exp(-bt)},\frac{y}{1+\varepsilon b \exp(-bt)},t+\frac{\ln(1+b\varepsilon \exp(-bt))}{b},u\right), \\
G_{3}: & (x,y, t+\varepsilon, \exp\left(-\varepsilon b(a+d)\right)u), \\
G_{4}: & (x, y,t, \exp(\varepsilon)u).
\end{align*}
\end{doublespace}
}
\noindent Since each one parameter group $G_{i}$ is a symmetry group, using Theorem \ref{newSol} if $u(x,y, t)$ is a solution of (\ref{I}), so are the $u^{\varepsilon,v_{i}},\,1\leq i\leq 4$ given by:

\begin{doublespace}
\begin{align*}\ds
u^{\varepsilon,v_{1}}(x,y,t)&=\frac{\exp\left(\frac{2b^{2}\epsilon \exp(bt)(x+y)}{1+b\epsilon \exp(bt)}\right)}{(1+b\epsilon \exp(bt))^{2(a+d)}}
u\left(\frac{x}{1+\varepsilon b \exp(bt)},\frac{y}{1+\varepsilon b \exp(bt)},t-\frac{\ln(1+b\epsilon \exp(bt))}{b}\right), \\
u^{\varepsilon,v_{2}}(x,y,t)&=u\left(\frac{x}{1-\varepsilon b \exp(-bt)},\frac{y}{1-\varepsilon b \exp(-bt)},t+\frac{\ln(1-\exp(-bt)b\varepsilon)}{b}\right), \\
u^{\varepsilon,v_{3}}(x,y,t)&=\exp(-\varepsilon b(a+d)) u(x,y, t-\varepsilon), \\
u^{\varepsilon,v_{4}}(x,y,t)&=\exp(\varepsilon) u(x,y, t).
\end{align*}
\end{doublespace}

\subsubsection{ Subcase 2.3 $(b=-e\,,\,b\ne 0\,,\,a\ne\frac{1}{4},d\ne\frac{1}{4})$}
{\small
\begin{doublespace}
\begin{align*}\ds
    v_{1}&=b\exp(bt)x\frac{\partial}{\partial x}+b\exp(bt)y\frac{\partial}{\partial y}+\exp(bt)\frac{\partial}{\partial t}-2b(-bx+a)\exp(bt)u\frac{\partial}{\partial u},\\
     v_{2}&=-b\exp(-bt)x\frac{\partial}{\partial x}-b\exp(-bt)y\frac{\partial}{\partial y}+\exp(-bt)\frac{\partial}{\partial t}+2b(by+d)\exp(-bt)u\frac{\partial}{\partial u},\\
       v_{3}&=\frac{\partial}{\partial t}+b(db-ab)u\frac{\partial}{\partial u},\,\,v_{4}=u\frac{\partial}{\partial u}.     
\end{align*}
\end{doublespace}
}
{\small
\begin{longtable}{|c|c|c|c|c|c|} 
\caption{Commutator table($b=-e\,,\,b\ne0\,,\,a\ne\frac{1}{4},d\ne \frac{1}{4}$).}\\
\hline
$[v_{i},v_{j}]$ & $v_{1}$ & $v_{2}$&$v_{3}$&$v_{4}$& $v_{\mu}$ \endhead
\hline
$v_{1}$&0&$-2bv_{3}$&$-bv_{1}$&0&$v_{\mu_{1}}$\\ 
$v_{2}$&$2bv_{3}$&0&$2bv_{2}$&0&$v_{\mu_{2}}$\\
$v_{3}$&$bv_{1}$&$-bv_{2}$&0&0&$v_{\mu^{3}}$\\
$v_{4}$&0&0&0&0&$-v_{\mu}$\\
$v_{\mu}$&$-v_{\mu_{1}}$&$-v_{\mu_{2}}$&$-v_{\mu_{3}}$&$v_{\mu}$&0\\ 
\hline
\end{longtable}}
{\small
\begin{doublespace}
\begin{align*}\ds
    \mu_{1}&=b\exp(bt)x\mu_{x}+b\exp(bt)y\mu_{y}+\exp(bt)\mu_{t}+2b(a-bx)\exp(bt)\mu,\\
    \mu_{2}&=-b\exp(-bt)x\mu_{x}-b\exp(-bt)y\mu_{y}+\exp(-bt)\mu_{t}-2b(d+bx)\exp(-bt)\mu,\\
     \mu_{3}&=b^{2}(a-d)\mu+\mu_{t}.
\end{align*}
\end{doublespace}
}
\noindent The transformed pointby a one-parameter group $\exp \left(\varepsilon v_{i}\right)(x,y, t, u)=$ $(\tilde{x},\tilde{y}, \tilde{t}, \tilde{u})$ are:
{\small
\begin{doublespace}
    \begin{align*}\ds
G_{1}: & \left(\frac{x}{1-b\varepsilon \exp(bt)},\frac{y}{1-b\varepsilon \exp(bt)}, t-\frac{\ln(1-b\varepsilon \exp(bt))}{b}, \frac{u\exp\left(\frac{2b^{2}\varepsilon \exp(bt)x}{1-b\varepsilon \exp(bt)}\right)}{1-b\varepsilon \exp(bt)}\right), \\
G_{2}: & ~ \Bigg(\frac{x}{1+b\varepsilon \exp(-bt)},\frac{y}{1+b\varepsilon \exp(-bt)},t+\frac{\ln(1+b\varepsilon \exp(-bt))}{b} ,\frac{\exp\left(-\frac{2b\left(\exp(-bt)dt+b\varepsilon(dt-y)\right)}{b\varepsilon+\exp(-bt)}\right) u}{(b\varepsilon +\exp(-bt))^{-2d}}\Bigg), \\
G_{3}: & \left( x,y,  t+\varepsilon,b(db - ab)\varepsilon + u(x, y, t)\right), \\
G_{4}: & \left(x,y, t, \exp(\varepsilon) u\right).
 \end{align*}
\end{doublespace}
}
\noindent Since each one parameter group $G_{i}$ is a symmetry group, using Theorem \ref{newSol} if $u(x,y, t)$ is a solution of (\ref{I}), so are the $u^{\varepsilon,v_{i}},\,1\leq i\leq 4$ given by:
{\small
\begin{doublespace}
 \begin{align*}
u^{\varepsilon,v_{1}}(x,y,t)&=(1+b\varepsilon \exp(bt))\exp\left(\frac{2b^{2}\varepsilon \exp(bt)x}{1+b\varepsilon \exp(bt)}\right)u\left(\frac{x}{1+b\varepsilon \exp(bt)},\frac{y}{1+b\varepsilon \exp(bt)},t-\frac{\ln(1+b\varepsilon \exp(bt))}{b}\right), \\
u^{\varepsilon,v_{2}}(x,y,t)&=\frac{\exp\left(\frac{2b^{2}\varepsilon \exp(-bt)y}{1-b\varepsilon \exp(-bt)}\right)}{(1-b\varepsilon \exp(-bt))^{2d}}u\left(\frac{x}{1-b\varepsilon \exp(-bt)},\frac{y}{1-b\varepsilon \exp(-bt)}, t+\frac{\ln(1-b\varepsilon \exp(-bt))}{b}\right), \\
u^{\varepsilon,v_{3}}(x,y,t)&=\exp(\varepsilon b(d-a))u(x,y, t-\varepsilon), \\
u^{\varepsilon,v_{4}}(x,y,t)&=\exp(\varepsilon)u\left( x,y,t\right)
\end{align*}
\end{doublespace}
}
\subsubsection{Subcase 2.4$(b=e=0\,,\,a\ne\frac{1}{4},d\ne\frac{1}{4})$}
{\small
\begin{align*}
    v_{1}&=2tx\frac{\partial}{\partial x}+2ty\frac{\partial}{\partial y}+t^{2}\frac{\partial}{\partial t}+(2(x+y)-2(a+d)t)u\frac{\partial}{\partial u},\,\, v_{2}=x\frac{\partial}{\partial x}+y\frac{\partial}{\partial y}+t\frac{\partial}{\partial t}\,\,,\,\, v_{3}=\frac{\partial}{\partial t}\,\,,\,\, v_{4}=u\frac{\partial}{\partial u}.
\end{align*}
}
\begin{longtable}{|c|c|c|c|c|c|} 
\caption{Commutator table($b=e=0\,,\,a\ne\frac{1}{4},d\ne \frac{1}{4}$).}\\
\hline
$[v_{i},v_{j}]$ & $v_{1}$ & $v_{2}$&$v_{3}$&$v_{4}$& $v_{\mu}$ \endhead
\hline
$v_{1}$&0&$-v_{1}$&$-2v_{2}+2(a+d)v_{4}$&0&$v_{\mu_{1}}$\\ 
$v_{2}$&$v_{1}$&0&$-v_{3}$&0&$v_{\mu_{2}}$\\
$v_{3}$&$2v_{2}-2(a+d)v_{4}$&$v_{3}$&$0$&0&$v_{\mu_{3}}$\\
$v_{4}$&0&0&0&0&$-v_{\mu}$\\
$v_{\mu}$&$-v_{\mu_{1}}$&$-v_{\mu_{2}}$&$-v_{\mu_{3}}$&$v_{\mu}$&0\\
\hline
\end{longtable}
{\small
\begin{doublespace}
\begin{align*}\ds
     \mu_{1}&=2tx\mu_{x}+2ty\mu_{y}+t^{2}\mu_{t}+(2t(a+d)-2(x+y))\mu,\,\,\mu_{2}=x\mu_{x}+y\mu_{y}+t\mu_{t}\,\,,\,\,\mu_{3}=\mu_{t}.    
\end{align*}
\end{doublespace}
}
\noindent The transformed pointby a one-parameter group $\exp \left(\varepsilon v_{i}\right)(x,y, t, u)=$ $(\tilde{x},\tilde{y}, \tilde{t}, \tilde{u})$ are:
{\small
\begin{doublespace}
\begin{align*}\ds
\nonumber G_{1}: &  \left(\frac{x}{(1-t\varepsilon)^{2}},\frac{y}{(1-t\varepsilon)^{2}},\frac{t}{1-t\varepsilon},u\exp\left(\frac{2\varepsilon(x+y)}{1-t\varepsilon}\right)(1-t\varepsilon)^{a+d}\right), \\
\nonumber G_{2}: & (\exp(\varepsilon)x,\exp(\varepsilon)y,\exp(\varepsilon)t,u), \\
\nonumber G_{3}: & ~(x,y,t+\varepsilon,u), \\
\nonumber G_{4}: & (x,y,t,\exp(\varepsilon)u).
\end{align*}
\end{doublespace}
}
\noindent Since each one parameter group $G_{i}$ is a symmetry group, using Theorem \ref{newSol} if $u(x,y, t)$ is a solution of (\ref{I}), so are the $u^{\varepsilon,v_{i}},\,1\leq i\leq 9$ given by: 

\begin{doublespace}
\begin{align*}\ds
u^{\varepsilon,v_{1}}(x,y,t)&=\frac{\exp\left(\frac{2\varepsilon(x+y)}{(1+t\varepsilon)^{3}}\right)}{(1+t\varepsilon)^{a+d}}u\left(\frac{x}{(1+t\varepsilon)^{2}},\frac{y}{(1+t\varepsilon)^{2}},\frac{t}{(1+t\varepsilon)}\right), \\
u^{\varepsilon,v_{2}}(x,y,t)&=u\left(\exp(-\varepsilon)x,\exp(-\varepsilon)y,\exp(-\varepsilon)t\right),\\
u^{\varepsilon,v_{3}}(x,y,t)&=u(x,y,t-\varepsilon), \\
u^{\varepsilon,v_{4}}(x,y,t)&=\exp(\varepsilon) u(x,y, t).
\end{align*}
\end{doublespace}
\subsection{Case 3 ($a=\frac{1}{4}\,,\,d\ne\frac{1}{4}$)}
Since $d\ne\frac{1}{4}$ and from remark \eqref{h,k,l}, using \eqref{b^{2}} we find that  $h=k=0$ and $l_{tt}-\frac{b^{2}}{4}l=0$ in all the following subcases.
 \subsubsection{Subcase 3.1 ($b\ne\pm e\,,\,a=\frac{1}{4},d\ne \frac{1}{4}$)}
 Using \eqref{b^{2}} we find out that $\tau_{t}=0$. Thus 
$$\xi=\sqrt{x}\left(c_{3}\exp(\frac{b}{2}t)+c_{4}\exp(-\frac{b}{2}t)\right),\,\,\gamma=0,\,\,\tau=c_{1},\,\,\phi=\left(2bc_{3}\exp(\frac{b}{2}t)\sqrt{x}+c_{2}\right)u+\mu(x,y,t)$$
{\small
\begin{align*}
    v_{1}&=\frac{\partial}{\partial t}\,\,,\,\,v_{2}=u\frac{\partial}{\partial u}\,\,,\,\, v_{3}=\sqrt{x}\exp(\frac{b}{2}t)\frac{\partial}{\partial x}+2b\exp(\frac{b}{2}t)\sqrt{x}u\frac{\partial}{\partial u}\,\,,\,\,v_{4}=\sqrt{x}\exp(-\frac{b}{2}t)\frac{\partial}{\partial x}.
\end{align*}
}
\begin{longtable}{|c|c|c|c|c|c|} 
\caption{Commutator table($b\ne\pm e\,,\,a=\frac{1}{4},d\ne \frac{1}{4}$).}\\
\hline
$[v_{i},v_{j}]$ & $v_{1}$ & $v_{2}$&$v_{3}$&$v_{4}$& $v_{\mu}$ \endhead
\hline
$v_{1}$&0&0&$\frac{b}{2}v_{3}$&$-\frac{b}{2}v_{4}$&$v_{\mu_{t}}$\\ 
$v_{2}$&0&0&0&0&$-v_{\mu}$\\
$v_{3}$&$-\frac{b}{2}v_{3}$&0&0&$-bv_{2}$&$v_{\mu_{1}}$\\
$v_{4}$&$\frac{b}{2}v_{4}$&0&$bv_{2}$&0&$v_{\mu_{2}}$\\
$v_{\mu}$&$-v_{\mu_{t}}$&$v_{\mu}$&$-v_{\mu^{1}}$&$-v_{\mu^{2}}$&0\\ 
\hline
\end{longtable}
{\small
\begin{align*}
\mu_{1}&=-2\exp(\frac{b}{2}t)\sqrt{x}b\mu+\exp(\frac{b}{2}t)\sqrt{x}\mu_{x}\,\,,\,\,\mu_{2}=\exp(-\frac{b}{2}t)\sqrt{x}\mu_{x}.  
\end{align*}
}
The transformed point  by a one-parameter group $\exp \left(\varepsilon v_{i}\right)(x,y, t, u)= (\tilde{x},\tilde{y}, \tilde{t}, \tilde{u})$ are:
{\small
\begin{doublespace}
\begin{align*}\ds
G_{1}: & (x,y, t+\varepsilon, u), \\
G_{2}: & (x,y, t, \exp(\varepsilon)u), \\
G_{3}: & \left(\frac{\left(\varepsilon \exp(\frac{bt}{2})+2\sqrt{x}\right)^{2}}{4},y, t,u\exp\left(\frac{ \exp(\frac{bt}{2})b\varepsilon\left(\varepsilon \exp(\frac{bt}{2})+4\sqrt{x}\right)}{2}\right)\right), \\
G_{4}: & \left(\frac{\left(2\exp(\frac{bt}{2})\sqrt{x}+\varepsilon\right)^{2}\exp(-bt)}{4},y, t,u\right).
\end{align*}
\end{doublespace}
}
\noindent Since each one parameter group $G_{i}$ is a symmetry group, using Theorem \ref{newSol} if $u(x,y, t)$ is a solution of (\ref{I}), so are the $u^{\varepsilon,v_{i}},\,1\leq i\leq 4$ given by:
{\small
\begin{doublespace}
\begin{align*}\ds
u^{\varepsilon,v_{1}}(x,y,t)&=u(x,y, t-\varepsilon), \\
u^{\varepsilon,v_{2}}(x,y,t)&=\exp(\varepsilon) u(x,y, t), \\
u^{\varepsilon,v_{3}}(x,y,t)&=\exp\left(\frac{ \exp(\frac{bt}{2})b\varepsilon\left(4\sqrt{x}-\varepsilon \exp(\frac{bt}{2})\right)}{2}\right)u\left(\frac{\left(2\sqrt{x}-\epsilon \exp(\frac{bt}{2})\right)^{2}}{4},y,t\right), \\
u^{\varepsilon,v_{4}}(x,y,t)&=u\left(\frac{\left(2\exp(\frac{bt}{2})\sqrt{x}-\epsilon\right)^{2}\exp(-bt)}{4},y,t\right).
\end{align*}
\end{doublespace}
}
\subsubsection{Subcase 3.2 ($b=e\,,\,b\ne0\,,\,a=\frac{1}{4},d\ne\frac{1}{4}$)}\label{paragraph 4.3.2}
{\small
\begin{doublespace}
\begin{align*}\ds
    v_{1}&=b\exp(bt)x\frac{\partial}{\partial x}+b\exp(bt)y\frac{\partial}{\partial y}+\exp(bt)\frac{\partial}{\partial t}+2b\exp(bt)\left(b(x+y)-(\frac{1}{4}+d)\right)u\frac{\partial}{\partial u},\\
     v_{2}&=-b\exp(-bt)x\frac{\partial}{\partial x}-b\exp(-bt)y\frac{\partial}{\partial y}+\exp(-bt)\frac{\partial}{\partial t},\,\,v_{3}=\frac{\partial}{\partial t}-b(\frac{1}{4}+d)u\frac{\partial}{\partial u}\\v_{4}&=\sqrt{x}\exp(\frac{b}{2}t)\frac{\partial}{\partial x}+2b\sqrt{x}\exp(\frac{b}{2}t)u\frac{\partial}{\partial u},\,\,v_{5}=\sqrt{x}\exp(-\frac{b}{2}t)\frac{\partial}{\partial x},\,\, v_{6}=u\frac{\partial}{\partial u}.
\end{align*}
\end{doublespace}
}
{\small
\begin{longtable}{|c|c|c|c|c|c|c|c|} 
\caption{Commutator table($b=e\,,\,b\ne0\,,\,a=\frac{1}{4},d\ne \frac{1}{4}$).}\\
\hline
$[v_{i},v_{j}]$ & $v_{1}$ & $v_{2}$&$v_{3}$&$v_{4}$&$v_{5}$&$v_{6}$& $v_{\mu}$ \endhead
\hline
$v_{1}$&0&$-2bv_{3}$&$-bv_{1}$&0&$-bv_{4}$&0&$v_{\mu_{1}}$\\ 
$v_{2}$&$2bv_{3}$&0&$bv_{2}$&$bv_{5}$&0&0&$v_{\mu_{2}}$\\
$v_{3}$&$bv_{1}$&$-bv_{2}$&0&$\frac{b}{2}v_{4}$&$-\frac{b}{2}v_{5}$&0&$v_{\mu_{3}}$\\
$v_{4}$&0&$-bv_{5}$&$-\frac{b}{2}v_{4}$&0&$-bv_{6}$&0&$v_{\mu_{4}}$\\
$v_{5}$&$bv_{4}$&0&$\frac{b}{2}v_{5}$&$bv_{6}$&0&0&$v_{\mu_{5}}$\\
$v_{6}$&0&0&0&0&0&0&$-v_{\mu}$\\ 
$v_{\mu}$&$-v_{\mu_{1}}$&$-v_{\mu_{2}}$&$-v_{\mu_{3}}$&$-v_{\mu_{4}}$&$-v_{\mu_{5}}$&$v_{\mu}$&0\\ 
\hline
\end{longtable}}
{\small
\begin{doublespace}
\begin{align*}\ds
     \mu_{1}&=b\exp(bt)x\mu_{x}+b\exp(bt)y\mu_{y}+\exp(bt)\mu_{t}+2b\exp(bt)(d+\frac{1}{4}-b(x+y)))\mu,\\
     \mu_{2}&=-b\exp(-bt)x\mu_{x}-b\exp(-bt)y\mu_{y}+\exp(-bt)\mu_{t}\,\,,\,\, \mu_{3}=b(a+d)\mu+\mu_{t},\\
     \mu_{4}&=\sqrt{x}\exp(\frac{b}{2}t)\mu_{x}-2b\exp(\frac{b}{2}t)\sqrt{x}\mu\,\,,\,\,\mu_{5}=\sqrt{x}\exp(-\frac{b}{2}t)\mu_{x}.
\end{align*}
\end{doublespace}
}
\noindent The transformed pointby a one-parameter group $\exp \left(\varepsilon v_{i}\right)(x,y, t, u)=$ $(\tilde{x},\tilde{y}, \tilde{t}, \tilde{u})$ are:
{\small
\begin{doublespace}
\begin{align*}\ds
G_{1}: & ~ \Bigg(\frac{x}{1-b\varepsilon \exp(bt)},\frac{y}{1-b\varepsilon \exp(bt)}, t-\frac{\ln(1-b\varepsilon \exp(bt))}{b},(1-b\varepsilon \exp(bt))^{(2d+\frac{1}{2})}u\exp\left(\frac{2b^{2}\varepsilon \exp(bt)(x+y)}{1-b\varepsilon \exp(bt)}\right)\Bigg), \\
G_{2}: & \left(\frac{x}{1+b\varepsilon \exp(-bt)},\frac{y}{1+b\varepsilon \exp(-bt)}, t+\frac{\ln(1+b\varepsilon \exp(-bt))}{b},u\right), \\
G_{3}: & \left(x,y, t+\varepsilon, \exp(-\varepsilon b(\frac{1}{4}+d))u\right), \\
G_{4}: & \left(\frac{\left(\varepsilon \exp(\frac{bt}{2})+2\sqrt{x}\right)^{2}}{4},y, t,u\exp\left(\frac{  \exp(\frac{bt}{2})b\varepsilon\left(\varepsilon  \exp(\frac{bt}{2})+4\sqrt{x}\right)}{2}\right)\right), \\
G_{5}: & \left(\frac{\left(2\exp(\frac{bt}{2})\sqrt{x}+\varepsilon\right)^{2}\exp(-bt)}{4},y, t,u\right) \\
G_{6}: & \left(x,y,t,\exp(\varepsilon)u\right)
\end{align*}
\end{doublespace}
}
\noindent Since each one parameter group $G_{i}$ is a symmetry group, using Theorem \ref{newSol} if $u(x,y, t)$ is a solution of (\ref{I}), so are the $u^{\varepsilon,v_{i}},\,1\leq i\leq 6$ given by:
{\small
\begin{doublespace}
\begin{align*}\ds
u^{\varepsilon,v_{1}}(x,y,t)&=\frac{\exp\left(\frac{2b^{2}\epsilon \exp(bt)(x+y)}{1+b\epsilon \exp(bt)}\right)}{(1+b\epsilon \exp(bt))^{2(\frac{1}{4}+d)}}
 u\left(\frac{x}{1+\varepsilon b \exp(bt)},\frac{y}{1+\varepsilon b \exp(bt)},t-\frac{\ln(1+b\epsilon \exp(bt))}{b}\right), \\
u^{\varepsilon,v_{2}}(x,y,t)&=u\left(\frac{x}{1-\varepsilon b \exp(-bt)},\frac{y}{1-\varepsilon b \exp(-bt)},t+\frac{\ln(1-\exp(-bt)b\varepsilon)}{b}\right),\\
u^{\varepsilon,v_{3}}(x,y,t)&=\exp\left(-\varepsilon b(\frac{1}{4}+d)\right) u(x,y, t-\varepsilon), \\
u^{\varepsilon,v_{4}}(x,y,t)&=\exp\left(\frac{b\epsilon \exp(\frac{bt}{2})(4\sqrt{x}-\epsilon \exp(\frac{bt}{2}))}{2}\right)u\left(\frac{\left(2\sqrt{x}-\epsilon \exp(\frac{bt}{2})\right)^{2}}{4},y,t\right),  \\
u^{\varepsilon,v_{5}}(x,y,t)&=u\left(\frac{\left(2\exp(\frac{bt}{2})\sqrt{x}-\epsilon\right)^{2}}{4}\exp(-bt),y,t\right), \\
u^{\varepsilon,v_{6}}(x,y,t)&=\exp(\varepsilon) u(x,y,t).
\end{align*}
\end{doublespace}
}
\subsubsection{Subcase 3.3 $(b=-e\,,\,b\ne0\,,\,a=\frac{1}{4},d\ne\frac{1}{4})$}
{\small
\begin{doublespace}
\begin{align*}\ds
    v_{1}&=b\exp(bt)x\frac{\partial}{\partial x}+b\exp(bt)y\frac{\partial}{\partial y}+\exp(bt)\frac{\partial}{\partial t}+\frac{b(4bx-1)}{2}\exp(bt)c_{1}u\frac{\partial}{\partial u},\\
     v_{2}&=-b\exp(-bt)x\frac{\partial}{\partial x}-b\exp(-bt)y\frac{\partial}{\partial y}+\exp(-bt)\frac{\partial}{\partial t}+2b(by+d)\exp(-bt)u\frac{\partial}{\partial u},\\
       v_{3}&=\frac{\partial}{\partial t}-b(\frac{1}{4}-d)u\frac{\partial}{\partial u}\,\,,\,\, v_{4}=\sqrt{x}\exp(\frac{b}{2}t)\frac{\partial}{\partial x}+2b\sqrt{x}\exp(\frac{b}{2}t)u\frac{\partial}{\partial u},\,\, v_{5}=\sqrt{x}\exp(-\frac{b}{2}t)\frac{\partial}{\partial x},\,\,v_{6}=u\frac{\partial}{\partial u}.
\end{align*}
\end{doublespace}
}
{\small
\begin{longtable}{|c|c|c|c|c|c|c|c|} 
\caption{Commutator table($b=-e\,,\,b\ne0\,,\,a=\frac{1}{4},d\ne \frac{1}{4}$).}\\
\hline
 $[v_{i},v_{j}]$& $v_{1}$ & $v_{2}$&$v_{3}$&$v_{4}$&$v_{5}$&$v_{6}$& $v_{\mu}$ \endhead
\hline
$v_{1}$&0&$-2bv_{3}$&$-bv_{1}$&0&$-bv_{4}$&0&$v_{\mu_{1}}$\\ 
$v_{2}$&$2bv_{3}$&0&$bv_{2}$&$bv_{5}$&0&0&$v_{\mu_{2}}$\\
$v_{3}$&$bv_{1}$&$-bv_{2}$&0&$\frac{b}{2}v_{4}$&$-\frac{b}{2}v_{5}$&0&$v_{\mu_{3}}$\\
$v_{4}$&0&$-bv_{5}$&$-\frac{b}{2}v_{4}$&0&$-bv_{6}$&0&$v_{\mu_{4}}$\\
$v_{5}$&$bv_{4}$&0&$\frac{b}{2}v_{5}$&$bv_{6}$&0&0&$v_{\mu_{5}}$\\
$v_{6}$&0&0&0&0&0&0&$-v_{\mu}$\\ 
$v_{\mu}$&$-v_{\mu_{1}}$&$-v_{\mu_{2}}$&$-v_{\mu_{3}}$&$-v_{\mu_{4}}$&$-v_{\mu_{5}}$&$v_{\mu}$&0\\ 
\hline
\end{longtable}}
{\small
\begin{doublespace}
\begin{align*}\ds
    \mu_{1}&=b\exp(bt)x\mu_{x}+b\exp(bt)y\mu_{y}+\exp(bt)\mu_{t}+2b(\frac{1}{4}-bx)\exp(bt)\mu,\\
    \mu_{2}&=-b\exp(-bt)x\mu_{x}-b\exp(-bt)y\mu_{y}+\exp(-bt)\mu_{t}-2b(d+by)\exp(-bt)\mu,\\
    \mu_{3}&=b(\frac{1}{4}-d)\mu+\mu_{t}\,\,,\,\, \mu_{4}=\exp(\frac{bt}{2})\sqrt{x}\mu_{x}-2b\exp(\frac{bt}{2})\sqrt{x}\mu\,\,,\,\,\mu_{5}=\exp(-\frac{bt}{2})\sqrt{x}\mu_{x}.
\end{align*}
\end{doublespace}
}
\noindent The transformed pointby a one-parameter group $\exp \left(\varepsilon v_{i}\right)(x,y, t, u)=$ $(\tilde{x},\tilde{y}, \tilde{t}, \tilde{u})$ are:
{\small
\begin{doublespace}
\begin{align*}\ds
G_{1}: & \left(\frac{x}{1-b\varepsilon \exp(bt)},\frac{y}{1-b\varepsilon \exp(bt)}, t-\frac{\ln(1-b\varepsilon \exp(bt))}{b}, \sqrt{1-b\varepsilon \exp(bt)}u\exp\left(\frac{2b^{2}\varepsilon \exp(bt)x}{1-b\varepsilon \exp(bt)}\right)\right), \\
G_{2}: & ~ \Bigg(\frac{x}{1+b\varepsilon \exp(-bt)},\frac{y}{1+b\varepsilon \exp(-bt)}, t+\frac{\ln(1-b\varepsilon \exp(-bt))}{b},(b\varepsilon +\exp(bt))^{2d}u\exp\left(\frac{-2b\left(dt\exp(bt)+b\varepsilon(dt-y)\right)}{b\varepsilon +\exp(bt)}\right)\Bigg), \\
G_{3}: & \left( x,y,  t+\varepsilon, -b(\frac{1}{4} - d)\varepsilon + u(x, y, t)\right), \\
G_{4}: & \left(\frac{\left(\varepsilon \exp(\frac{bt}{2})+2\sqrt{x}\right)^{2}}{4},y, t,u\exp\left(\frac{ \exp(\frac{bt}{2})b\varepsilon\left(\varepsilon \exp(\frac{bt}{2})+4\sqrt{x}\right)}{2}\right)\right), \\
G_{5}: & \left(\frac{\left(2\exp(\frac{bt}{2})\sqrt{x}+\varepsilon\right)^{2}\exp(-bt)}{4},y, t,u\right), \\
G_{6}: & (x,y,t,u\exp(\varepsilon)). 
\end{align*}
\end{doublespace}
}
\noindent Since each one parameter group $G_{i}$ is a symmetry group, using Theorem \ref{newSol} if $u(x,y, t)$ is a solution of (\ref{I}), so are the $u^{\varepsilon,v_{i}},\,1\leq i\leq 6$ given by:

\begin{doublespace}
\begin{align*}\ds
u^{\varepsilon,v_{1}}(x,y,t)&=\frac{\exp\left(\frac{2b^{2}\varepsilon x \exp(bt) }{1+b\varepsilon \exp(bt)}\right)}{\sqrt{1+b\varepsilon \exp(bt)}}u\left(\frac{x}{1+b\varepsilon \exp(bt)},\frac{y}{1+b\varepsilon \exp(bt)},t-\frac{\ln(1+\varepsilon b \exp(bt))}{b}\right), \\
u^{\varepsilon,v_{2}}(x,y,t)&=\frac{\exp\left(-\frac{b\varepsilon \exp(-bt)y}{1-b\varepsilon \exp(-bt)}\right)}{(1-b\varepsilon \exp(-bt))^{2}}
u\left(\frac{x}{(1-\varepsilon b \exp(-bt))},\frac{y}{(1-\varepsilon b \exp(-bt))},t+\frac{\ln(1+b\varepsilon \exp(-bt))}{b}\right), \\
u^{\varepsilon,v_{3}}(x,y,t)&=\exp\left(-\varepsilon b(\frac{1}{4}-d)\right) u(x,y, t-\varepsilon), \\
u^{\varepsilon,v_{4}}(x,y,t)&=\exp\left(\frac{b\varepsilon \exp(\frac{bt}{2})\left(4\sqrt{x}-\varepsilon \exp(\frac{bt}{2})\right)}{2}\right)u\left(\frac{(2\sqrt{x}-\varepsilon \exp(\frac{bt}{2}))^{2}}{4},y,t\right), \\
u^{\varepsilon,v_{5}}(x,y,t)&=u\left(\frac{(2\exp(\frac{bt}{2})\sqrt{x}-\epsilon)^{2}}{4}\exp(-bt),y,t\right), \\
u^{\varepsilon,v_{6}}(x,y,t)&=\exp(\varepsilon) u(x,y,t).
\end{align*}
\end{doublespace}

\subsubsection{Subcase 3.4$(b=e=0\,,\,a=\frac{1}{4},d\ne\frac{1}{4})$}

\begin{align*}
    v_{1}&=2tx\frac{\partial}{\partial x}+2ty\frac{\partial}{\partial y}+t^{2}\frac{\partial}{\partial t}+(2(x+y)-2(\frac{1}{4}+d)t)u\frac{\partial}{\partial u},\\
    v_{2}&=x\frac{\partial}{\partial x}+y\frac{\partial}{\partial y}+t\frac{\partial}{\partial t},\,\,
    v_{3}=\frac{\partial}{\partial t},\\
    v_{4}&=\sqrt{x}t\frac{\partial}{\partial x}+2\sqrt{x}u\frac{\partial}{\partial u},\,\,v_{5}=\sqrt{x}\frac{\partial}{\partial x},\\v_{6}&=u\frac{\partial}{\partial u}.
\end{align*}

\begin{longtable}{|c|c|c|c|c|c|c|c|} 
\caption{Commutator table($b=e=0\,,\,b\ne0\,,\,a=\frac{1}{4},d\ne \frac{1}{4}$).}\\
\hline
 $[v_{i},v_{j}]$& $v_{1}$ & $v_{2}$&$v_{3}$&$v_{4}$&$v_{5}$&$v_{6}$& $v_{\mu}$ \endhead
\hline
$v_{1}$&0&$-v_{1}$&$-2v_{2}+2(d+\frac{1}{4})v_{6}$&$0$&$-v_{4}$&0&$v_{\mu_{1}}$\\ 
$v_{2}$&$v_{1}$&0&$-v_{3}$&$\frac{1}{2}v_{4}$&$-\frac{1}{2}v_{5}$&0&$v_{\mu_{2}}$\\
$v_{3}$&$2v_{2}-2(d+\frac{1}{4})v_{6}$&$v_{3}$&0&$v_{5}$&$0$&0&$v_{\mu_{3}}$\\
$v_{4}$&$0$&$-\frac{1}{2}v_{4}$&$-v_{5}$&0&$-v_{6}$&0&$v_{\mu_{4}}$\\
$v_{5}$&$v_{4}$&$\frac{1}{2}v_{5}$&$0$&$v_{6}$&0&0&$v_{\mu_{5}}$\\
$v_{6}$&0&0&0&0&0&0&$-v_{\mu}$\\ 
$v_{\mu}$&$-v_{\mu_{1}}$&$-v_{\mu_{2}}$&$-v_{\mu_{3}}$&$-v_{\mu_{4}}$&$-v_{\mu_{5}}$&$v_{\mu}$&0\\ 
\hline
\end{longtable}

\begin{doublespace}
\begin{align*}\ds
     \mu_{1}&=2tx\mu_{x}+2ty\mu_{y}+t^{2}\mu_{t}+(2t(\frac{1}{4}+d)-2(x+y))\mu,\\
     \mu_{2}&=x\mu_{x}+y\mu_{y}+t\mu_{t}\,\,,\,\,\mu_{3}=\mu_{t}\,\,,\,\,\mu_{4}=\sqrt{x}t\mu_{x}-2\sqrt{x}\mu\,\,,\,\,\mu_{5}=\sqrt{x}\mu_{x}.  
\end{align*}
\end{doublespace}

\noindent The transformed pointby a one-parameter group $\exp \left(\varepsilon v_{i}\right)(x,y, t, u)=$ $(\tilde{x},\tilde{y}, \tilde{t}, \tilde{u})$ are:

\begin{doublespace}
\begin{align*}\ds
\nonumber G_{1}: &  \left(\frac{x}{(1-t\varepsilon)^{2}},\frac{y}{(1-t\varepsilon)^{2}},\frac{t}{1-t\varepsilon},u\exp\left(\frac{2\varepsilon(x+y)}{1-t\varepsilon}\right)(1-t\varepsilon)^{\frac{1}{4}+d}\right), \\
\nonumber G_{2}: & (\exp(\varepsilon)x,\exp(\varepsilon)y,\exp(\varepsilon)t,u), \\
\nonumber G_{3}: & ~(x,y,t+\varepsilon,u), \\
\nonumber G_{4}: & \left(x+t\sqrt{x}\varepsilon +\frac{t^{2}\varepsilon^{2}}{4},y,t,u\exp\left(\frac{\varepsilon(t\varepsilon+4\sqrt{x})}{2}\right)\right), \\
\nonumber G_{5}: & \left(\sqrt{x}\varepsilon+\frac{\varepsilon^{2}}{4}+x,y,t,u\right), \\
\nonumber G_{6}: & (x,y,t,\exp(\varepsilon)u).
\end{align*}
\end{doublespace}

\noindent Since each one parameter group $G_{i}$ is a symmetry group, using Theorem \ref{newSol} if $u(x,y, t)$ is a solution of (\ref{I}), so are the $u^{\varepsilon,v_{i}},\,1\leq i\leq 6$ given by: 

\begin{doublespace}
\begin{align*}\ds
u^{\varepsilon,v_{1}}(x,y,t)&=\frac{\exp\left(\frac{2\varepsilon(x+y)}{(1+t\varepsilon)^{3}}\right)}{(1+t\varepsilon)^{\frac{1}{4}+d}}u\left(\frac{x}{(1+t\varepsilon)^{2}},\frac{y}{(1+t\varepsilon)^{2}},\frac{t}{(1+t\varepsilon)}\right), \\
u^{\varepsilon,v_{2}}(x,y,t)&=u\left(\exp(-\varepsilon)x,\exp(-\varepsilon)y,\exp(-\varepsilon)t\right),\\
u^{\varepsilon,v_{3}}(x,y,t)&=u(x,y,t-\varepsilon), \\
u^{\varepsilon,v_{4}}(x,y,t)&=\exp\left(\frac{\varepsilon \left(\varepsilon t+4\sqrt{x-t\sqrt{x}\varepsilon+\frac{t^{2}\varepsilon^{2}}{2}}\right)}{2}\right)u\left(x-t\sqrt{x}\varepsilon+\frac{t^{2}\varepsilon^{2}}{2},y,t\right), \\
u^{\varepsilon,v_{5}}(x,y,t)&=u\left(\sqrt{x}\varepsilon-\frac{\varepsilon^{2}}{4}+x,y,t\right), \\
u^{\varepsilon,v_{6}}(x,y,t)&=\exp(\varepsilon) u(x,y, t).
\end{align*}
\end{doublespace}

\subsection{Case 4 ($a\ne\frac{1}{4}\,,\,d=\frac{1}{4}$)}
This case is symmetric to the last case (case 3), because $a\ne\frac{1}{4}$ and from remark \eqref{h,k,l}, using \eqref{k} we find out that  $h=l=0$ and $k_{tt}-\frac{e^{2}}{4}k=0$ in all the following subcases.
\subsubsection{Subcase 4,1 ($b\ne\pm e\,,\,a\ne\frac{1}{4},d= \frac{1}{4}$)}

\begin{align*}
  v_{1}&=\frac{\partial}{\partial t}\,\,,\,\,v_{2}=u\frac{\partial}{\partial u},\\
    v_{3}&=\sqrt{y}\exp(\frac{e}{2}t)\frac{\partial}{\partial y}+2e\exp(\frac{e}{2}t)\sqrt{y}u\frac{\partial}{\partial u},\\
    v_{4}&=\sqrt{y}\exp(-\frac{e}{2}t)\frac{\partial}{\partial y}.
\end{align*}

\begin{longtable}{|c|c|c|c|c|c|} 
\caption{Commutator table($b\ne\pm e\,,\,a\ne\frac{1}{4},d= \frac{1}{4}$).}\\
\hline
$[v_{i},v_{j}]$ & $v_{1}$ & $v_{2}$&$v_{3}$&$v_{4}$& $v_{\mu}$ \endhead
\hline
$v_{1}$&0&0&$\frac{e}{2}v_{3}$&$-\frac{e}{2}v_{4}$&$v_{\mu_{t}}$\\ 
$v_{2}$&0&0&0&0&$-v_{\mu}$\\
$v_{3}$&$-\frac{e}{2}v_{3}$&0&0&$-ev_{2}$&$v_{\mu_{1}}$\\
$v_{4}$&$\frac{e}{2}v_{4}$&0&$ev_{2}$&0&$v_{\mu^{2}}$\\
$v_{\mu}$&$-v_{\mu_{t}}$&$v_{\mu}$&$-v_{\mu_{1}}$&$-v_{\mu_{2}}$&0\\ 
\hline
\end{longtable}

\begin{doublespace}
\begin{align*}\ds
 \mu_{1}&=-2\exp(\frac{e}{2}t)\sqrt{y}e\mu+\exp(\frac{e}{2}t)\sqrt{y}\mu_{y}\,\,,\,\,\mu_{2}=\exp(-\frac{e}{2}t)\sqrt{y}\mu_{y}.
\end{align*}
\end{doublespace}

\noindent The transformed point by a one-parameter group $\exp \left(\varepsilon v_{i}\right)(x,y, t, u)=$ $(\tilde{x},\tilde{y}, \tilde{t}, \tilde{u})$ are:

\begin{doublespace}
\begin{align*}\ds
G_{1}: & (x,y, t+\varepsilon, u), \\
G_{2}: & (x,y, t, \exp(\varepsilon)u), \\
G_{3}: & \left(x,\frac{\left(\varepsilon \exp(\frac{et}{2})+2\sqrt{y}\right)^{2}}{4}, t,u\exp\left(\frac{\exp(\frac{et}{2})e\varepsilon\left(\varepsilon \exp(\frac{et}{2})+4\sqrt{y}\right)}{2}\right)\right), \\
G_{4}: & \left(x,\frac{\left(2\exp(\frac{et}{2})\sqrt{y}+\varepsilon\right)^{2}\exp(-et)}{4}, t,u\right).
\end{align*}
\end{doublespace}

\noindent Since each one parameter group $G_{i}$ is a symmetry group, using Theorem \ref{newSol} if $u(x,y, t)$ is a solution of (\ref{I}), so are the $u^{\varepsilon,v_{i}},\,1\leq i\leq 4$ given by:

\begin{doublespace}
\begin{align*}\ds
u^{\varepsilon,v_{1}}(x,y,t)&=u(x,y, t-\varepsilon), \\
u^{\varepsilon,v_{2}}(x,y,t)&=\exp(\varepsilon) u(x,y, t), \\
u^{\varepsilon,v_{3}}(x,y,t)&=\exp\left(\frac{\exp(\frac{et}{2})e\varepsilon\left( 4\sqrt{y}-\varepsilon\exp(\frac{et}{2}\right)}{2}\right)u\left(\frac{\left(2\sqrt{y}-\epsilon \exp(\frac{et}{2})\right)^{2}}{4},y,t\right), \\
u^{\varepsilon,v_{4}}(x,y,t)&=u\left(x,\frac{\left(2\exp(\frac{et}{2})\sqrt{y}-\epsilon\right)^{2}\exp(-et)}{4},t\right).
\end{align*}
\end{doublespace}
\subsubsection{Subcase 4.2 ($b=e\,,\,b\ne0\,,\,a\ne\frac{1}{4},d=\frac{1}{4}$)}

\begin{doublespace}
\begin{align*}
    v_{1}&=b\exp(bt)x\frac{\partial}{\partial x}+b\exp(bt)y\frac{\partial}{\partial y}+\exp(bt)\frac{\partial}{\partial t}+2b\exp(bt)\left(b(x+y)-(a+\frac{1}{4})\right)u\frac{\partial}{\partial u},\\
     v_{2}&=-b\exp(-bt)x\frac{\partial}{\partial x}-b\exp(-bt)y\frac{\partial}{\partial y}+\exp(-bt)\frac{\partial}{\partial t},\,\,v_{3}=\frac{\partial}{\partial t}-b(\frac{1}{4}+a)u\frac{\partial}{\partial u}\\
     v_{4}&=\sqrt{y}\exp(\frac{b}{2}t)\frac{\partial}{\partial x}+2b\sqrt{y}\exp(\frac{b}{2}t)u\frac{\partial}{\partial u},\,\, v_{5}=\sqrt{y}\exp(-\frac{b}{2}t)\frac{\partial}{\partial y}\,\,,\,\,v_{6}=u\frac{\partial}{\partial u}.    
\end{align*}
\end{doublespace}

{\small
\begin{longtable}{|c|c|c|c|c|c|c|c|} 
\caption{Commutator table($b=e\,,\,a\ne\frac{1}{4},d= \frac{1}{4}$).}\\
\hline
 $[v_{i},v_{j}]$& $v_{1}$ & $v_{2}$&$v_{3}$&$v_{4}$&$v_{5}$&$v_{6}$& $v_{\mu}$ \endhead
\hline
$v_{1}$&0&$-2bv_{3}$&$-bv_{1}$&0&$-bv_{4}$&0&$v_{\mu_{1}}$\\ 
$v_{2}$&$2bv_{3}$&0&$bv_{2}$&$bv_{5}$&0&0&$v_{\mu_{2}}$\\
$v_{3}$&$bv_{1}$&$-bv_{2}$&0&$\frac{b}{2}v_{4}$&$-\frac{b}{2}v_{5}$&0&$v_{\mu_{3}}$\\
$v_{4}$&0&$-bv_{5}$&$-\frac{b}{2}v_{4}$&0&$-bv_{6}$&0&$v_{\mu_{4}}$\\
$v_{5}$&$bv_{4}$&0&$\frac{b}{2}v_{5}$&$bv_{6}$&0&0&$v_{\mu_{5}}$\\
$v_{6}$&0&0&0&0&0&0&$-v_{\mu}$\\ 
$v_{\mu}$&$-v_{\mu_{1}}$&$-v_{\mu_{2}}$&$-v_{\mu_{3}}$&$-v_{\mu_{4}}$&$-v_{\mu_{5}}$&$v_{\mu}$&0\\ 
\hline
\end{longtable}}

\begin{doublespace}
\begin{align*}
     \mu_{1}&=b\exp(bt)x\mu_{x}+b\exp(bt)y\mu_{y}+\exp(bt)\mu_{t}+2b\exp(bt)(d+\frac{1}{4}-b(x+y)))\mu,\\
      \mu_{2}&=-b\exp(-bt)x\mu_{x}-b\exp(-bt)y\mu_{y}+\exp(-bt)\mu_{t},\\
      \mu_{3}&=b(a+d)\mu+\mu_{t}\,\,,\,\, \mu_{4}=\sqrt{y}\exp(\frac{b}{2}t)\mu_{y}-2b\exp(\frac{b}{2}t)\sqrt{x}\mu\,\,,\,\, \mu_{5}=\sqrt{y}\exp(-\frac{b}{2}t)\mu_{y}.
\end{align*}
\end{doublespace}

\noindent The transformed pointby a one-parameter group $\exp \left(\varepsilon v_{i}\right)(x,y, t, u)=$ $(\tilde{x},\tilde{y}, \tilde{t}, \tilde{u})$ are:
{\small
\begin{doublespace}
\begin{align*}\ds
G_{1}: &~ \Bigg(\frac{x}{1-b\varepsilon \exp(bt)},\frac{y}{1-b\varepsilon \exp(bt)}, t-\frac{\ln(1-b\varepsilon \exp(bt))}{b},
(1-b\varepsilon \exp(bt))^{2a+\frac{1}{2}}u\exp\left(\frac{2b^{2}\varepsilon \exp(bt)(x+y)}{1-b\varepsilon \exp(bt)}\right)\Bigg),\\
G_{2}: & \left(\frac{x}{1+b\varepsilon \exp(-bt)},\frac{y}{1+b\varepsilon \exp(-bt)}, t+\frac{\ln(1+b\varepsilon \exp(-bt))}{b},u\right), \\
G_{3}: & \left(x,y, t+\varepsilon, \exp(-\varepsilon b(\frac{1}{4}+d))u\right), \\
G_{4}: & \left(x,\frac{\left(\varepsilon \exp(\frac{bt}{2})+2\sqrt{y}\right)^{2}}{4}, t,u\exp\left(\frac{ \exp(\frac{bt}{2})b\varepsilon\left(\varepsilon \exp(\frac{bt}{2})+4\sqrt{y}\right)}{2}\right)\right), \\
G_{5}: & \left(x,\frac{\left(2\exp(\frac{bt}{2})\sqrt{y}+\varepsilon\right)^{2}\exp(-bt)}{4}, t,u\right) \\
G_{6}: & \left(x,y,t,\exp(\varepsilon)u\right)
\end{align*}
\end{doublespace}
}
\noindent Since each one parameter group $G_{i}$ is a symmetry group,using Theorem \ref{newSol} if $u(x,y, t)$ is a solution of (\ref{I}), so are the $u^{\varepsilon,v_{i}},\,1\leq i\leq 6$ given by:

\begin{doublespace}
\begin{align*}\ds
u^{\varepsilon,v_{1}}(x,y,t)&=\frac{\exp\left(\frac{2b^{2}\epsilon \exp(bt)(x+y)}{1+b\epsilon \exp(bt)}\right)}{(1+b\epsilon \exp(bt))^{2(\frac{1}{4}+a)}}
u\left(\frac{x}{1+\varepsilon b \exp(bt)},\frac{y}{1+\varepsilon b \exp(bt)},t-\frac{\ln(1+b\epsilon \exp(bt))}{b}\right), \\
u^{\varepsilon,v_{2}}(x,y,t)&=u\left(\frac{x}{1-\varepsilon b \exp(-bt)},\frac{y}{1-\varepsilon b \exp(-bt)},t+\frac{\ln(1-\exp(-bt)b\varepsilon)}{b}\right),\\
u^{\varepsilon,v_{3}}(x,y,t)&=\exp(-\varepsilon b(\frac{1}{4}+a)) f(x,y, t-\varepsilon), \\
u^{\varepsilon,v_{4}}(x,y,t)&=\exp\left(\frac{b\epsilon \exp(\frac{bt}{2})(4\sqrt{y}-\epsilon \exp(\frac{bt}{2}))}{2}\right)u\left(x,\frac{\left(2\sqrt{y}-\epsilon \exp(\frac{bt}{2})\right)^{2}}{4},t\right)  \\
u^{\varepsilon,v_{5}}(x,y,t)&=u\left(x,\frac{\left(2\exp(\frac{bt}{2})\sqrt{y}-\epsilon\right)^{2}}{4}\exp(-bt),t\right), \\
u^{\varepsilon,v_{6}}(x,y,t)&=\exp(\varepsilon) u(x,y,t).
\end{align*}
\end{doublespace}

\subsubsection{ Subcase 4.3$(b=-e\,,\,b\ne0\,,\,a\ne\frac{1}{4},d=\frac{1}{4})$}

\begin{doublespace}
\begin{align*}\ds
    v_{1}&=b\exp(bt)x\frac{\partial}{\partial x}+b\exp(bt)y\frac{\partial}{\partial y}+\exp(bt)\frac{\partial}{\partial t}-2b(-bx+a)\exp(bt)u\frac{\partial}{\partial u},\\
     v_{2}&=-b\exp(-bt)x\frac{\partial}{\partial x}-b\exp(-bt)y\frac{\partial}{\partial y}+\exp(-bt)\frac{\partial}{\partial t}+\frac{b(4by+1)}{2}\exp(-bt)u\frac{\partial}{\partial u},\,\, v_{3}=\frac{\partial}{\partial t}+b(\frac{1}{4}-a)u\frac{\partial}{\partial u}\\
     v_{4}&=\sqrt{y}\exp(-\frac{b}{2}t)\frac{\partial}{\partial x}-2b\sqrt{y}\exp(-\frac{b}{2}t)u\frac{\partial}{\partial u},\,\,
       v_{5}=\sqrt{y}\exp(\frac{b}{2}t)\frac{\partial}{\partial y},\,\,v_{6}=u\frac{\partial}{\partial u}.
\end{align*}
\end{doublespace}

\begin{longtable}{|c|c|c|c|c|c|c|c|} 
\caption{Commutator table($b=-e\,,\,a\ne\frac{1}{4},d= \frac{1}{4}$).}\\
\hline
 $[v_{i},v_{j}]$& $v_{1}$ & $v_{2}$&$v_{3}$&$v_{4}$&$v_{5}$&$v_{6}$& $v_{\mu}$ \endhead
\hline
$v_{1}$&0&$-2bv_{3}$&$-bv_{1}$&$-bv_{5}$&0&0&$v_{\mu_{1}}$\\ 
$v_{2}$&$2bv_{3}$&0&$bv_{2}$&0&$bv_{4}$&0&$v_{\mu_{2}}$\\
$v_{3}$&$bv_{1}$&$-bv_{2}$&0&$-\frac{b}{2}v_{4}$&$\frac{b}{2}v_{5}$&0&$v_{\mu_{3}}$\\
$v_{4}$&$bv_{5}$&0&$\frac{b}{2}v_{4}$&0&$bv_{6}$&0&$v_{\mu_{4}}$\\
$v_{5}$&0&$-bv_{4}$&$-\frac{b}{2}v_{5}$&$-bv_{6}$&0&0&$v_{\mu_{5}}$\\
$v_{6}$&0&0&0&0&0&0&$-v_{\mu}$\\ 
$v_{\mu}$&$-v_{\mu_{1}}$&$-v_{\mu_{2}}$&$-v_{\mu_{3}}$&$-v_{\mu_{4}}$&$-v_{\mu_{5}}$&$v_{\mu}$&0\\ 
\hline
\end{longtable}

\begin{doublespace}
\begin{align*}\ds
\mu_{1}&=b\exp(bt)x\mu_{x}+b\exp(bt)y\mu_{y}+\exp(bt)\mu_{t}+2b(a-bx)\exp(bt)\mu,\\
    \mu_{2}&=-b\exp(-bt)x\mu_{x}-b\exp(-bt)y\mu_{y}+\exp(-bt)\mu_{t}-2b(\frac{1}{4}+by)\exp(-bt)\mu,\\
    \mu_{3}&=-b(\frac{1}{4}-a)\mu+\mu_{t}\,\,,\,\,\mu_{4}=\exp(-\frac{bt}{2})\sqrt{y}\mu_{y}+2b\exp(-\frac{bt}{2})\sqrt{y}\mu\,\,,\,\,\mu_{5}=\exp(\frac{bt}{2})\sqrt{y}\mu_{y}.
\end{align*}
\end{doublespace}

\noindent The transformed pointby a one-parameter group $\exp \left(\varepsilon v_{i}\right)(x,y, t, u)=$ $(\tilde{x},\tilde{y}, \tilde{t}, \tilde{u})$ are:
{\small
\begin{doublespace}
\begin{align*}\ds
G_{1}: & \left(\frac{x}{1-b\varepsilon \exp(bt)},\frac{y}{1-b\varepsilon \exp(bt)}, t-\frac{\ln(1-b\varepsilon \exp(bt))}{b},\frac{1}{1-b\varepsilon \exp(bt)}u\exp\left(\frac{2\exp(bt)b^{2}\varepsilon x}{1-b\varepsilon \exp(bt)}\right)\right), \\
G_{2}: & ~ \Bigg(\frac{x}{1+b\varepsilon \exp(-bt)},\frac{y}{1+b\varepsilon \exp(-bt)}, t+\frac{\ln(1-b\varepsilon \exp(-bt))}{b},\sqrt{1+b\varepsilon \exp(-bt)}u\exp\left(\frac{2b^{2}\varepsilon \exp(-bt)y}{1+b\varepsilon \exp(-bt)}\right)\Bigg),\\
G_{3}: & \left( x,y,  t+\varepsilon, b\left(\frac{b}{4} - ab\right)\varepsilon + u(x, y, t)\right), \\
G_{4}: & \left(x,\frac{\left(\varepsilon \exp(-\frac{bt}{2})+2\sqrt{y}\right)^{2}}{4}, t,u\exp\left(\frac{ -b\exp(-\frac{bt}{2})\varepsilon\left(\varepsilon \exp(-\frac{bt}{2})+4\sqrt{y}\right)}{2}\right)\right), \\
G_{5}: & \left(x,\frac{\left(2\sqrt{y}+\varepsilon \exp(\frac{bt}{2})\right)^{2}}{4},t,u\right), \\
G_{6}: & (x,y,t,u\exp(\varepsilon).
\end{align*}
\end{doublespace}}

\noindent Since each one parameter group $G_{i}$ is a symmetry group, using Theorem \ref{newSol} if $u(x,y, t)$ is a solution of (\ref{I}), so are the $u^{\varepsilon,v_{i}},\,1\leq i\leq 6$ given by:
{\small
\begin{doublespace}
\begin{align*}\ds
 u^{\varepsilon,v_{1}}(x,y,t)&=(1+b\varepsilon \exp(bt))\exp\left(\frac{2b^{2}\varepsilon  \exp(bt) x}{1-b\varepsilon \exp(bt)}\right)u\left(\frac{x}{1+b\varepsilon \exp(bt)},\frac{y}{1+b\varepsilon \exp(bt)},t-\frac{\ln(1+\varepsilon b \exp(bt) )}{b}\right)\\
u^{\varepsilon,v_{2}}(x,y,t)&=\sqrt{1+\frac{b\varepsilon \exp(-bt)}{1+b\varepsilon \exp(-bt)}} \exp\left(\frac{2b^{2}\varepsilon \exp(-bt)y}{(1-b\varepsilon \exp(-bt))(1+2b\varepsilon \exp(-bt))}\right)\\
&u\left(\frac{x}{(1-\varepsilon b \exp(-bt))},\frac{y}{(1-\varepsilon b \exp(-bt))},t+\frac{\ln(1+b\varepsilon \exp(-bt))}{b}\right), \\
u^{\varepsilon,v_{3}}(x,y,t)&=\exp\left(\varepsilon b(\frac{1}{4}-a)\right) u(x,y, t-\varepsilon), \\
u^{\varepsilon,v_{4}}(x,y,t)&=\exp\left(-\frac{b\varepsilon \exp(-\frac{bt}{2})\left(4\sqrt{y}-\varepsilon \exp(-\frac{bt}{2})\right)}{2}\right)u\left(x,\frac{\left(2\sqrt{y}-\varepsilon \exp(-\frac{bt}{2})\right)^{2}}{4},t\right), \\
u^{\varepsilon,v_{5}}(x,y,t)&=u\left(x,\frac{\left(2\sqrt{y}-\epsilon \exp(\frac{bt}{2})\right)^{2}}{4},t\right), \\
u^{\varepsilon,v_{6}}(x,y,t)&=\exp(\varepsilon) u(x,y,t).
\end{align*}
\end{doublespace}
}
\subsubsection{Subcase 4.4$(b=e=0\,,\,a\ne\frac{1}{4},d=\frac{1}{4})$}


\begin{align*}
    v_{1}&=2tx\frac{\partial}{\partial x}+2ty\frac{\partial}{\partial y}+t^{2}\frac{\partial}{\partial t}+(2(x+y)-2(\frac{1}{4}+d)t)u\frac{\partial}{\partial u},\,\,v_{3}=\frac{\partial}{\partial t},\\
    v_{2}&=x\frac{\partial}{\partial x}+y\frac{\partial}{\partial y}+t\frac{\partial}{\partial t},\,\,v_{4}=\sqrt{y}t\frac{\partial}{\partial x}+2\sqrt{y}u\frac{\partial}{\partial u},\,\,v_{5}=\sqrt{y}\frac{\partial}{\partial y}\,\,,\,\,v_{6}=u\frac{\partial}{\partial u}. 
\end{align*}

\begin{longtable}{|c|c|c|c|c|c|c|c|} 
\caption{Commutator table($b=e=0\,,\,a\ne\frac{1}{4},d= \frac{1}{4}$).}\\
\hline
 $[v_{i},v_{j}]$& $v_{1}$ & $v_{2}$&$v_{3}$&$v_{4}$&$v_{5}$&$v_{6}$& $v_{\mu}$ \endhead
\hline
$v_{1}$&0&$-v_{1}$&$-2v_{2}+2(d+\frac{1}{4})v_{6}$&$0$&$-v_{4}$&0&$v_{\mu_{1}}$\\ 
$v_{2}$&$v_{1}$&0&$-v_{3}$&$\frac{1}{2}v_{4}$&$-\frac{1}{2}v_{5}$&0&$v_{\mu_{2}}$\\
$v_{3}$&$2v_{2}-2(d+\frac{1}{4})v_{6}$&$-bv_{2}$&0&$v_{5}$&$0$&0&$v_{\mu_{3}}$\\
$v_{4}$&$0$&$-\frac{1}{2}v_{4}$&$-v_{5}$&0&$-v_{6}$&0&$v_{\mu_{4}}$\\
$v_{5}$&$v_{4}$&$\frac{1}{2}v_{5}$&$0$&$-v_{6}$&0&0&$v_{\mu_{5}}$\\
$v_{6}$&0&0&0&0&0&0&$-v_{\mu}$\\ 
$v_{\mu}$&$-v_{\mu_{1}}$&$-v_{\mu_{2}}$&$-v_{\mu_{3}}$&$-v_{\mu_{4}}$&$-v_{\mu_{5}}$&$v_{\mu}$&0\\ 
\hline
\end{longtable}
{\small
\begin{doublespace}
\begin{align*}\ds
     \mu_{1}&=2tx\mu_{x}+2ty\mu_{y}+t^{2}\mu_{t}+(2t(\frac{1}{4}+d)-2(x+y))\mu,\\
     \mu_{2}&=x\mu_{x}+y\mu_{y}+t\mu_{t}\,\,,\,\,\mu_{3}=\mu_{t}\,\,,\,\,\mu_{4}=\sqrt{y}t\mu_{y}-2\sqrt{y}\mu\,\,,\,\,\mu_{5}=\sqrt{y}\mu_{y}.
\end{align*}
\end{doublespace}
}

\noindent The transformed pointby a one-parameter group $\exp \left(\varepsilon v_{i}\right)(x,y, t, u)=$ $(\tilde{x},\tilde{y}, \tilde{t}, \tilde{u})$ are:
{\small
\begin{doublespace}
\begin{align*}\ds
\nonumber G_{1}: & \left(\frac{x}{(1-t\varepsilon)^{2}},\frac{y}{(1-t\varepsilon)^{2}},\frac{t}{1-t\varepsilon},u\exp\left(\frac{2\varepsilon(x+y)}{1-t\varepsilon}\right)(1-t\varepsilon)^{a+\frac{1}{4}}\right), \\
\nonumber G_{2}: & (\exp(\varepsilon)x,\exp(\varepsilon)y,\exp(\varepsilon)t,u), \\
\nonumber G_{3}: & ~(x,y,t+\varepsilon,u), \\
\nonumber G_{4}: & \left(x,y+t\sqrt{y}\varepsilon +\frac{t^{2}\varepsilon^{2}}{4},t,u\exp\left(\frac{\varepsilon(t\varepsilon+4\sqrt{y})}{2}\right)\right), \\
\nonumber G_{5}: & \left(x,\sqrt{y}\varepsilon+\frac{\varepsilon^{2}}{4}+y,t,u\right), \\
\nonumber G_{6}: & (x,y,t,\exp(\varepsilon)u).
\end{align*}
\end{doublespace}
}
\noindent Since each one parameter group $G_{i}$ is a symmetry group, using Theorem \ref{newSol} if $u(x,y, t)$ is a solution of (\ref{I}), so are the $u^{\varepsilon,v_{i}},\,1\leq i\leq 6$ given by: 
{\small
\begin{doublespace}
\begin{align*}\ds
u^{\varepsilon,v_{1}}(x,y,t)&=\frac{\exp\left(\frac{2\varepsilon(x+y)}{(1+t\varepsilon)^{3}}\right)}{(1+t\varepsilon)^{a+\frac{1}{4}}}u\left(\frac{x}{(1+t\varepsilon)^{2}},\frac{y}{(1+t\varepsilon)^{2}},\frac{t}{(1+t\varepsilon)}\right), \\
u^{\varepsilon,v_{2}}(x,y,t)&=u\left(\exp(-\varepsilon)x,\exp(-\varepsilon)y,\exp(-\varepsilon)t\right),\\
u^{\varepsilon,v_{3}}(x,y,t)&=u(x,y,t-\varepsilon), \\
u^{\varepsilon,v_{4}}(x,y,t)&=\exp\left(\frac{\varepsilon \left(\varepsilon t+4\sqrt{y-t\sqrt{y}\varepsilon+\frac{t^{2}\varepsilon^{2}}{2}}\right)}{2}\right)u\left(x,y-t\sqrt{y}\varepsilon+\frac{t^{2}\varepsilon^{2}}{2},t\right), \\
u^{\varepsilon,v_{5}}(x,y,t)&=u\left(x,-\sqrt{y}\varepsilon+\frac{\varepsilon^{2}}{4}+y,t\right), \\
u^{\varepsilon,v_{6}}(x,y,t)&=\exp(\varepsilon) u(x,y, t).
\end{align*}
\end{doublespace}
}
\section{The structure of the Lie algebras }
Recall $\mathcal{V}_{a,b,d,e}$ defined in Remark \ref{nu }, the algebra generated by $v_{1},...,v_{r}\,\,\text{and}\,\, v_{\mu},\,\text{where}\\ r\in \{2,4,6,9\}$. Similarly, $\mathcal{L}_{a,b,d,e}$ represent the finite- dimensional algebra generated by\\ $v_{1},...,v_{r}\,\,r\in \{2,4,6,9\}$. Additionally, $\mathcal{H}$ is the infinite-dimensional ideal of $\mathcal{V}_{a,b,d,e}$ consisting of the $v_{\mu}$ where $\mu$ is a solution of \eqref{I}. \\$[\mathcal{V}_{a,b,d,e},\mathcal{H}]\subset \mathcal{H}$ therefore  $\mathcal{V}_{a,b,d,e}=\mathcal{L}_{a,b,d,e}\ltimes \mathcal{H}$.\\
Now, we proceed to analyze the structure of $\mathcal{L}_{a,b,d,e}$ in each case.

\subsection{Case 1 ($a=\frac{1}{4}\,,\,d=\frac{1}{4}$)}
\subsubsection{Subcase 1.1 $(b\ne\pm e\,,\,a=d=\frac{1}{4})$}
 We have $\mathcal{Z}(\mathcal{L}_{a,b,d,e})=<v_{2}>$ thus $\mathcal{L}_{a,b,d,e}/\mathcal{Z} (\mathcal{L}_{a,b,d,e})\cong\mathbb{R}v_{1}\ltimes \mathcal{A}$ \, where $\mathcal{A}=<\Bar{v_{3}},\Bar{v_{4}},\Bar{v_{5}},\Bar{v_{6}}>$ is an abelian Lie algebra of dimension 4.
\begin{itemize}
     \item [$\blacksquare$]{if $b=0 \,,\, e\ne 0$:}(\ref{case1.1})\\
      From the commutator table, the center of $\mathcal{L}_{a,b,d,e}\,\, \text{is}\,\, \mathcal{Z} (\mathcal{L}_{a,b,d,e})=<v_{2},v_{3},v_{5}>$, which is an abelian Lie algebra of dimension 3. Then 

$$ \frac{\mathcal{L}_{a,b,d,e}}{{\mathcal{Z} (\mathcal{L}_{a,b,d,e})}}\cong \frac{\frac{\mathcal{L}_{a,b,d,e}}{\mathbb{R}v_{2}}}{\frac{\mathcal{Z} (\mathcal{L}_{a,b,d,e})}{\mathbb{R}v_{2}}}=\frac{\mathbb{R}v_{1}\ltimes <\Bar{v_{3}},\Bar{v_{4}},\Bar{v_{5}},\Bar{v_{6}}>}{<\Bar{v_{3}},\Bar{v_{5}}>}\cong \mathbb{R}\Bar{\Bar{v_{1}}}\ltimes <\Bar{\Bar{v_{4}}},\Bar{\Bar{v_{6}}}>\cong L(3,2,x=-1)\cong \mathfrak{iso}_{2}$$ which is the Poincaré algebra, using the isomorphism defined below 
      
     {\small $$
\phi:\left\{\begin{aligned}
e_{1}&\longrightarrow \Bar{\Bar{v_{4}}},\\
e_{2}&\longrightarrow \Bar{\Bar{v_{6}}},\\
e_{3}&\longrightarrow \frac{2}{e}\Bar{\Bar{v_{1}}}.
\end{aligned}\right.
$$}
where $\{e_{1},e_{2},e_{3}\}$ are the basis of $\mathfrak{iso}_{2}$, using the notation of \cite[p.~16]{Bowers}.
     \item [$\blacksquare$]{if $b\ne0 \,,\, e= 0$:}(\ref{case1.1})\\
    This case is symmetric to the last one. We have that, $\mathcal{Z} (\mathcal{L}_{a,b,d,e})=<\Bar{v_{2}},\Bar{v_{4}},\Bar{v_{6}}>$, then

$$ \frac{\mathcal{L}_{a,b,d,e}}{{\mathcal{Z} (\mathcal{L}_{a,b,d,e})}}\cong \frac{\frac{\mathcal{L}_{a,b,d,e}}{\mathbb{R}v_{2}}}{\frac{\mathcal{Z} (\mathcal{L}_{a,b,d,e})}{\mathbb{R}v_{2}}}=\frac{\mathbb{R}v_{1}\ltimes <\Bar{v_{3}},\Bar{v_{4}},\Bar{v_{5}},\Bar{v_{6}}>}{<\Bar{v_{4}},\Bar{v_{6}}>}\cong \mathbb{R}\Bar{\Bar{v_{1}}}\ltimes <\Bar{\Bar{v_{3}}},\Bar{\Bar{v_{5}}}>\cong L(3,2,x=-1)\cong \mathfrak{iso}_{2}$$ using the isomorphism defined below 
{\small  $$  
\psi:\left\{\begin{aligned}
e_{1}&\longrightarrow \Bar{\Bar{v_{5}}},\\
e_{2}&\longrightarrow \Bar{\Bar{v_{6}}},\\
e_{3}&\longrightarrow \frac{2}{b}\Bar{\Bar{v_{1}}}.
\end{aligned}\right.
$$}
where $\{e_{1},e_{2},e_{3}\}$ are the basis of $\mathfrak{iso}_{2}$, using the notation of \cite[p.~16]{Bowers}.

\end{itemize}
\subsubsection{Subcase 1.2$(b=e\,,\,b\ne 0\,\,a=\frac{1}{4}\,,\,d=\frac{1}{4})$}
    We have $J=<v_{4},v_{5},v_{6},v_{7},v_{8},v_{9}>$ is an ideal of $\mathcal{L}_{a,b,d,e}$ of dimension 6, then\\
    $\mathcal{L}_{a,b,d,e}/J=<\Bar{v_{1}},\Bar{v_{2}},\Bar{v_{3}}>\cong\mathfrak{sl}_{2}$ via the isomorphism 
   {\small \begin{subequations}\label{eq:system}
\begin{empheq}[left=\empheqlbrace]{align}\nonumber
\Bar{v_{1}}&\longrightarrow -bE,\\
\Bar{v_{2}}&\longrightarrow bF, \label{iso2.1}\\\nonumber
\Bar{v_{3}}&\longrightarrow \frac{b}{2}H,
\end{empheq}
\end{subequations}}

where $\{E,F,H\}$ is the canonical basis of $\mathfrak{sl}_{2}$, 
    then: $$\mathcal{L}_{a,b,d,e}\cong \mathfrak{sl}_{2}\ltimes J.$$
\subsubsection{Subcase 1.3$(b=-e\,,\,b\ne 0\,\,a=\frac{1}{4}\,,\,d=\frac{1}{4})$}
Similiar to the last case, $$\mathcal{L}_{a,b,d,e}\cong \mathfrak{sl}_{2}\ltimes J$$
\subsubsection{Subcase 1.4$(b=e=0\,,\,a=\frac{1}{4}\,,\,d=\frac{1}{4})$}

 We have that $J=<\Bar{v_{4}},\Bar{v_{5}},\Bar{v_{6}},\Bar{v_{7}},\Bar{v_{8}},\Bar{v_{9}}>$ is an ideal of $\mathcal{L}_{a,b,d,e}$ of dimension 6, and, $\mathcal{L}_{a,b,d,e}/J\cong \mathfrak{sl}_{2}$
via the isomorphism defined as follows:
 {\small\begin{subequations}\label{eq:system}
\begin{empheq}[left=\empheqlbrace]{align}\nonumber
\Bar{v_{1}}&\longrightarrow E,\\
\Bar{v_{2}}&\longrightarrow \frac{H}{2},\label{iso3.1} \\
\Bar{v_{3}}&\longrightarrow -F\nonumber,
\end{empheq}
\end{subequations}}
Thus, $$\mathcal{L}_{a,b,d,e}\cong \mathfrak{sl}_{2}\ltimes J. $$
\subsection{Case 2 ($a\ne\frac{1}{4}\,,\,d\ne\frac{1}{4}$)}
\subsubsection{Subcase 2.1 ($b\ne\pm e\,,\,a\ne\frac{1}{4},d\ne \frac{1}{4}$)}
From the commutator table, we deduce that 
$\mathcal{L}_{a,b,d,e}=<v_{1},v_{2}>$ is an abelian Lie algebra of dimension 2.
\subsubsection{Subcase 2.2$(b= e\,,\,b\ne 0\,\,a\ne\frac{1}{4}\,,\,d\ne\frac{1}{4})$}
We have $I=<v_{1},v_{2},v_{3}>\cong \mathfrak{sl}_{2}$ via the isomorphism (\ref{iso2.1}), and $\mathcal{Z}(\mathcal{L}_{a,b,d,e})=<v_{4}>$. Thus:
$$\mathcal{L}_{a,b,d,e}\cong \mathbb{R}v_{4}\times \mathfrak{sl}_{2}\cong \mathcal{N}.$$
\subsubsection{Subcase 2.3$(b=-e\,,\,b\ne 0\,\,a\ne\frac{1}{4}\,,\,d\ne\frac{1}{4})$}
Same as last case: $$\mathcal{L}_{a,b,d,e}\cong \mathbb{R}v_{4}\times \mathfrak{sl}_{2}\cong \mathcal{N}.$$
\subsubsection{Subcase 2.4$(b=e=0\,\,a\ne\frac{1}{4}\,,\,d\ne\frac{1}{4})$}
Again we have that: $$\mathcal{L}_{a,b,d,e}\cong\mathbb{R}v_{4}\times \mathfrak{sl}_{2}\cong \mathcal{N}.$$
 \subsection{Case 3 ($a=\frac{1}{4}\,,\,d\ne\frac{1}{4}$)}
 \subsubsection{Subcase 3.1 $(b\ne\pm e,\,\,a=\frac{1}{4},d\ne \frac{1}{4})$}
Suppose that $b\ne 0$, because the case $b=0\,\,\mathcal{L}_{a,b,d,e}$ is an abelian Lie algebra of dimension 4.
From the commutator table $\mathcal{Z}(\mathcal{L}_{a,b,d,e})=\mathbb{R}v_{2}$. Thus: $$\mathcal{L}_{a,b,d,e}/\mathbb{R}v_{2}\cong L(3,2,x=-1)\cong\mathfrak{iso}_{2},$$using the notations of \cite[p.~16]{Bowers}.\label{5.3.1}
\subsubsection{Subcase 3.2 $(b= e\,,\,b\ne 0\,,\,a=\frac{1}{4}\,d\ne\frac{1}{4})$}\label{5.3.2}

\noindent Using the table in paragraph \ref{paragraph 4.3.2}, one sees that the linear mapping from $\mathcal{L}_{a,b,d,e}$ to $\mathcal{M}$ defined by:  
\small \begin{subequations}\label{eq:system}
\begin{empheq}[left=\empheqlbrace]{align}
v_{1}&\longrightarrow E,\\
v_{2}&\longrightarrow -b^{2}F \label{iso2.1}\\
v_{3}&\longrightarrow \frac{b}{2}H,\\ 
v_{4}&\longrightarrow X,\\
v_{5}&\longrightarrow -bY,\\
v_{6}&\longrightarrow Z.
\end{empheq}
\end{subequations}

It is a morphism of Lie algebra. Therefore, $\mathcal{L}_{a,b,d,e}$ is isomorphic to $\mathcal{M}$, Hence to $\mathcal{H}_{0,0}$, used in \cite{lescot2014solving}.

\subsubsection{Subcase 3.3 $(b=-e\,,\,b\ne 0\,,\,a=\frac{1}{4}\,d\ne\frac{1}{4})$}
Similar to the previous case:
$$\mathcal{L}_{a,b,d,e}\cong\mathcal{M}\cong \mathcal{H}_{0,0}.$$
\subsubsection{Subcase 3.4$(b=e=0\,,\,a=\frac{1}{4}\,,\,d\ne\frac{1}{4})$}\label{5.3.3}
Again, in this case: $$\mathcal{L}_{a,b,d,e}\cong\mathcal{M}\cong \mathcal{H}_{0,0}.$$
 \subsection{Case 4 ($a\ne\frac{1}{4}\,,\,d=\frac{1}{4}$)}
 This case is symmetric to case 3.
 \subsubsection{Subcase 4.1 $(b\ne\pm e\,\,a\ne\frac{1}{4},d=\frac{1}{4})$}
We assume $e\ne 0$, because the case $e=0$ makes $\mathcal{L}_{a,b,d,e}$ an abelian algebra dimension 4. From the commutator table,
$$\mathcal{L}_{a,b,d,e}/\mathbb{R}v_{2}\cong \mathfrak{iso}_{2}$$
\subsubsection{Subcase 4.2$(b=e\,,\,b\ne 0\,\,a\ne\frac{1}{4}\,,d=\frac{1}{4})$}
 $$\mathcal{L}_{a,b,d,e}\cong\mathcal{M}\cong \mathcal{H}_{0,0}.$$
\subsubsection{Subcase 4.3$(b=-e\,,\,b\ne 0\,\,a\ne\frac{1}{4}\,,d=\frac{1}{4})$}
$$\mathcal{L}_{a,b,d,e}\cong\mathcal{M}\cong \mathcal{H}_{0,0}.$$
\subsubsection{Subcase 4.4$(b=e=0\,\,a\ne\frac{1}{4}\,,\,d=\frac{1}{4})$}
$$\mathcal{L}_{a,b,d,e}\cong\mathcal{M}\cong \mathcal{H}_{0,0}.$$
\section{Isomorphism between Lie algebras.}
In each case we designate the respective Lie algebras as follows:
{\scriptsize
\begin{align*}
&\mathcal{L}^{0}_{0} : b\ne \pm e \,\,\text{and}\,\, a\ne\frac{1}{4}\,\,d\ne\frac{1}{4},\hspace{0.2cm}\mathcal{L}^{1}_{0} : b=e \,\,\text{and}\,\, a\ne\frac{1}{4}\,\,d\ne\frac{1}{4},\hspace{0.2cm}\mathcal{L}^{2}_{0} :b=-e \,\,\text{and}\,\, a\ne\frac{1}{4}\,\,d\ne\frac{1}{4},\hspace{0.2cm}\mathcal{L}^{3}_{0} :b=e=0 \,\,\text{and}\,\, a\ne\frac{1}{4}\,\,d\ne\frac{1}{4},\\
&\mathcal{L}^{0}_{1} :  b\ne \pm e \,\,\text{and}\,\, a=\frac{1}{4}\,\,d\ne\frac{1}{4},\hspace{0.2cm}\mathcal{L}^{1}_{1} :  b=e \,\,\text{and}\,\, a=\frac{1}{4}\,\,d\ne\frac{1}{4},\hspace{0.2cm}\mathcal{L}^{2}_{1} :  b=-e \,\,\text{and}\,\, a=\frac{1}{4}\,\,d\ne\frac{1}{4},\hspace{0.2cm}\mathcal{L}^{3}_{1} :  b=e=0 \,\,\text{and}\,\, a=\frac{1}{4}\,\,d\ne\frac{1}{4},\\
&\mathcal{L}^{0}_{2} : b\ne \pm e \,\,\text{and}\,\, a\ne\frac{1}{4}\,\,d=\frac{1}{4},\hspace{0.2cm}\mathcal{L}^{1}_{2} : b=e \,\,\text{and}\,\, a\ne\frac{1}{4}\,\,d=\frac{1}{4},\hspace{0.2cm}\mathcal{L}^{2}_{2} : b=-e \,\,\text{and}\,\, a\ne\frac{1}{4}\,\,d=\frac{1}{4},\hspace{0.2cm}\mathcal{L}^{3}_{2} : b=e=0 \,\,\text{and}\,\, a\ne\frac{1}{4}\,\,d=\frac{1}{4},\\
&\mathcal{L}^{0}_{3}: b\ne \pm e \,\,\text{and}\,\, a=\frac{1}{4}\,\,d=\frac{1}{4},\hspace{0.2cm}\mathcal{L}^{1}_{3}: b=e \,\,\text{and}\,\, a=\frac{1}{4}\,\,d=\frac{1}{4}\hspace{0.2cm}\mathcal{L}^{2}_{3}: b=-e \,\,\text{and}\,\, a=\frac{1}{4}\,\,d=\frac{1}{4},\hspace{0.2cm}\mathcal{L}^{3}_{3}: b=e=0 \,\,\text{and}\,\, a=\frac{1}{4}\,\,d=\frac{1}{4}.
\end{align*}}
 In all cases we noticed that
\begin{align*}
     &\mathcal{L}^{0}_{0}\subset \mathcal{L}^{0}_{1}\subset \mathcal{L}^{0}_{3},\quad \mathcal{L}^{1}_{0}\subset \mathcal{L}^{1}_{1}\subset \mathcal{L}^{1}_{3},\quad\mathcal{L}^{2}_{0}\subset \mathcal{L}^{2}_{1}\subset \mathcal{L}^{2}_{3},\\
      &\mathcal{L}^{0}_{0}\subset \mathcal{L}^{0}_{2}\subset \mathcal{L}^{0}_{3},\quad\mathcal{L}^{1}_{0}\subset \mathcal{L}^{1}_{2}\subset \mathcal{L}^{1}_{3},\quad \mathcal{L}^{2}_{0}\subset \mathcal{L}^{2}_{2}\subset \mathcal{L}^{2}_{3}.     
\end{align*}
 $$ \text{ and}\quad\mathcal{L}^{1}_{0}\cong \mathcal{L}^{2}_{0}\cong  \mathcal{L}^{3}_{0},\hspace{0.3cm}\mathcal{L}^{1}_{1}\cong \mathcal{L}^{2}_{1}\cong  \mathcal{L}^{3}_{1}\cong \mathcal{L}^{1}_{2}\cong \mathcal{L}^{2}_{2}\cong  \mathcal{L}^{3}_{2},\hspace{0.3cm}\mathcal{L}^{0}_{1}\cong  \mathcal{L}^{0}_{2},\hspace{0.3cm}\mathcal{L}^{1}_{3}\cong  \mathcal{L}^{2}_{3} \cong  \mathcal{L}^{3}_{3}.$$
We summarize our results in the following table :
{\doublespacing
\begin{longtable}{|c|c|c|c|}
\caption{Summarized results.}\\
\hline
 \multirow{2}{*}{Cases}&  \multirow{2}{*}{Subcases} & Finite dimensional  & \multirow{2}{*}{isomorphism }
 \\
 & & Lie algebra of symmetry & between Lie algebras  \\\hline
\multirow{4}{*}{$a=\frac{1}{4}\,d=\frac{1}{4}$}&$b\ne \pm e$&$\mathcal{L}^{0}_{3}$ : 6 vector fields& $\mathcal{L}^{0}_{3}\cong \mathfrak{iso}_{2}$ \\ \cline{2-3} 
&$b=e\,,\,b\ne 0$&$\mathcal{L}^{1}_{3}$ : 9&$\mathcal{L}^{1}_{3}\cong \mathcal{L}^{2}_{3}\cong \mathcal{L}^{3}_{3}\cong \mathfrak{sl}_{2}\ltimes J  $ \\\cline{2-3} &  $b=-e\,,\,b\ne 0$&$\mathcal{L}^{2}_{3}$ : 9&\\\cline{2-3}&$b=e=0$&$\mathcal{L}^{3}_{3}$ : 9&\\\hline
&$b\ne \pm e$&$\mathcal{L}^{0}_{0}$ : 2 vector fields&\\ \cline{2-3}
$a\ne\frac{1}{4}\,d\ne\frac{1}{4}$&$b=e\,,\,b\ne 0$&$\mathcal{L}^{1}_{0}$ : 4& $\mathcal{L}^{1}_{0}\cong \mathcal{L}^{2}_{0}\cong \mathcal{L}^{3}_{0}\cong \mathbb{R}v_{4}\times\mathfrak{sl}_{2}$\\\cline{2-3}
&$b=-e\,,\,b\ne 0$&$\mathcal{L}^{2}_{0}$ : 4&\\\cline{2-3}&$b=e=0$&$\mathcal{L}^{3}_{0}$ : 4&\\\hline
 &$b\ne \pm e$&$\mathcal{L}^{0}_{1}$ : 4 vector fields &$\mathcal{L}^{0}_{1}\cong \mathfrak{iso}_{2}$\\ \cline{2-3}
$a=\frac{1}{4}\,d\ne\frac{1}{4}$&$b=e\,,\,b\ne 0 $&$\mathcal{L}^{1}_{1}$ : 6& $\mathcal{L}^{1}_{1}\cong \mathcal{L}^{2}_{1}\cong \mathcal{L}^{3}_{1}\cong \mathfrak{sl}_{2}\ltimes \mathfrak{h}_{3}$\\\cline{2-3}&$b=-e\,,\,b\ne 0$&$\mathcal{L}^{2}_{1}$ : 6&\\\cline{2-3}&$b=e=0$&$\mathcal{L}^{3}_{1}$ : 6&\\\hline
 &$b\ne \pm e$&$\mathcal{L}^{0}_{2}$ : 4 vector fields&$\mathcal{L}^{0}_{2}\cong \mathfrak{iso}_{2}$\\ \cline{2-3}
$a\ne\frac{1}{4}\,d=\frac{1}{4}$&$b=e\,,\,b\ne 0$&$\mathcal{L}^{1}_{2}$ : 6& $\mathcal{L}^{1}_{2}\cong \mathcal{L}^{2}_{2}\cong \mathcal{L}^{3}_{2}\cong \mathfrak{sl}_{2}\ltimes \mathfrak{h}_{3}$\\\cline{2-3}&$b=-e\,,\,b\ne 0$&$\mathcal{L}^{2}_{2}$ : 6&\\\cline{2-3}&$b=e=0$&$\mathcal{L}^{3}_{2}$ : 6&\\
\hline

\end{longtable}
}

\begin{remark}
  From the beginning, we observed that the PDE class (\ref{KBELS}) admits the permutation of $x$ and $y$ as an equivalence transformation(symmetry group). Therefore, we might have anticipated that $\mathcal{L}_{a,b,d,e}$ and $ \mathcal{L}_{d,e,a,b} $ would be isomorphic.
\end{remark}
\section{Conclusion}
We examined the Lie point symmetries of the Kolmogorov backward equation \eqref{KBELS} by, analyzing Lie algebras and their structures. Furthermore, we can extend this analysis to compute the Lie symmetries of stochastic differential equation \eqref{LS} associated with partial differential equation \eqref{KBELS}. These symmetries are referred to as W-symmetries or random Lie point symmetries following the notation introduced by Gaeta \cite{Giuseppe-Gaeta_2000}. For more comprehensive information, one can refer to Gaeta's works \cite{Giuseppe-Gaeta_1999,Giuseppe-Gaeta_2000} as well as the contributions by Kozlov \cite{kozlov2010symmetries}, who haves made additional advancements in this theory.
For example, let us consider the Longstaff-Schwartz model \eqref{LS} in the case $(b=e=0)$,\, and $(a=d=\frac{1}{4})$. Using Theorem 3.5 in \cite{kozlov2020symmetries}, we find that the SDE has the following symmetries:
$$ 
    v_{2}=x\frac{\partial}{\partial x}+y\frac{\partial}{\partial y}+\frac{\partial}{\partial t}\,,\,
     v_{4}=\sqrt{x}\frac{\partial}{\partial x}\,,\,
       v_{6}=\sqrt{y}\frac{\partial}{\partial y}\,,\,
       v_{8}=\frac{\partial}{\partial t}.
$$
Exponentiating the symmetry $v_{4}$, gives the Lie group of transformation  associated to $v_{4}$ which is\\
$(\Tilde{x},\Tilde{y},\Tilde{t})=((\sqrt{x}+\frac{\varepsilon}{2})^{2},y,t)$, using Theorem 3.5 in \cite{kozlov2020symmetries}, $\Tilde{X_{t}}=(\sqrt{X_{t}}+\frac{\varepsilon}{2})^{2}$ must be a solution of the SDE 
$d X_{t} =\frac{1}{4} d t+\sqrt{X_{t}} d B^{1}_{t}$. Using Ito formula, one finds that,
\begin{align*}
    d\Tilde{X_{t}}&=\left(1+\frac{\varepsilon}{2\sqrt{X_{t}}}\right)(a d t+\sqrt{X_{t}} d B^{1}_{t})-\frac{1}{8}\left(\varepsilon X_{t}^{-\frac{1}{2}}\right)dt\\
    &=adt+\frac{\varepsilon}{2\sqrt{X_{t}}}(a-\frac{1}{4})+(\sqrt{X_{t}}+\frac{\varepsilon}{2})dB^{1}_{t}\\
     d\Tilde{X_{t}}&=\frac{1}{4}dt+\sqrt{\tilde{X_{t}}}dB^{1}_{t} 
\end{align*}
Whence, $\Tilde{X_{t}}$ is a solution of the SDE.

\section*{Acknowledgement}
This paper and the research behind it would not have been possible without the exceptional support of Mr Paul LESCOT. His enthusiasm, knowledge and exacting attention to details have been an inspiration and kept my work on track.
\bibliography{biblio.bib}

\bibliographystyle{plain}

\end{document}